\RequirePackage{fix-cm}
\documentclass[smallcondensed]{svjour3}     % onecolumn (ditto)
\smartqed  % flush right qed marks, e.g. at end of proof
\usepackage{graphicx}
\usepackage{latexsym}
\usepackage{amsmath}
\usepackage{amssymb}
\usepackage{float}

\newcommand\dt{{\partial_{t}}}
\newcommand\Dt{{\Delta t}}
\newcommand\ex{{\rm ex}}
\newcommand\FD{{\rm FD}}

\newcommand\vp{\varphi}
\newcommand\ve{\varepsilon}

\newcommand\cB{\mathcal{B}}
\newcommand\cE{\mathcal{E}}
\newcommand\cK{\mathcal{K}}

\newcommand\cP{\mathcal{P}}

\newcommand\NN{\mathbb{N}}
\newcommand\RR{\mathbb{R}}
\newcommand\ZZ{\mathbb{Z}}

\newcommand\Nr {N_{\rm r}}
\newcommand\Dtr{\Dt_{\rm r}}

\newcommand\fin{{f^0}}
\newcommand\rmd{\, {\rm d}} %for integrals

\providecommand{\abs}[1]{\lvert#1\rvert}
\providecommand{\bigabs}[1]{\big|#1\big|}
\providecommand{\Bigabs}[1]{\Big|#1\Big|}

\providecommand{\norm}[1]{\lVert#1\rVert}
\providecommand{\Bignorm}[1]{\Big\|#1\Big\|}

\DeclareMathOperator{\Div}{div}
\DeclareMathOperator{\supp}{supp}
\DeclareMathOperator{\erf}{erf}
\DeclareMathOperator{\curl}{curl}

\begin{document}

\title{Smooth particle methods without smoothing%\thanks{Grants or other notes
%about the article that should go on the front page should be
%placed here. General acknowledgments should be placed at the end of the article.}
}
%\subtitle{Do you have a subtitle?\\ If so, write it here}

%\titlerunning{Short form of title}        % if too long for running head

\author{Martin Campos Pinto
%         \and
%        Second Author %etc.
}

%\authorrunning{Short form of author list} % if too long for running head

\institute{M. Campos Pinto \at
              CNRS, UMR 7598, Laboratoire Jacques-Louis Lions, F-75005, Paris, France \\
              UPMC Univ Paris 06, UMR 7598, Laboratoire Jacques-Louis Lions, F-75005, Paris, France \\
              Tel.: +33 1 44 27 91 54\\
              Fax:  +33 1 44 27 72 00\\
              \email{campos@ann.jussieu.fr}           %  \\
    %             \emph{Present address:} of F. Author  %  if needed
    %           \and
    %           S. Author \at
    %              second address
}
    %\institute{F. Author \at
    %              first address \\
    %              Tel.: +123-45-678910\\
    %              Fax: +123-45-678910\\
    %              \email{fauthor@example.com}           %  \\
    %%             \emph{Present address:} of F. Author  %  if needed
    %           \and
    %           S. Author \at
    %              second address
    %}

\date{Received: date / Accepted: date}
% The correct dates will be entered by the editor

\maketitle

\begin{abstract}
We present a new class of particle methods with deformable shapes that converge in the uniform norm without requiring remappings,
extended overlapping or vanishing moments for the particles. 
The crux of the method is to use polynomial expansions of the backward characteristic flow to transport the numerical particles 
with improved accuracy: in the first order case the method consists of representing the transported density with 
linearly-transformed particles, the second order version computes quadratically-transformed particles, and so on.
For programming purposes we provide explicit implementations of the resulting LTP and QTP schemes
that only involve pointwise evaluations of the forward characteristic flow, and also come with local indicators for 
the accuracy of the corresponding transport scheme.
Numerical tests using different transport problems demonstrate the accuracy of the proposed methods compared to 
standard particle schemes, and establish their robustness with respect to the remapping period. In particular, it is 
shown that QTP particles can be transported without remappings on very long periods of times, without 
hampering the accuracy of the numerical solutions.
Finally, a dynamic criterion is proposed to automatically select the time steps where the particles should be remapped. 
The strategy does not require additional inter-particle communications, and it is validated by numerical experiments.

% \PACS{PACS code1 \and PACS code2 \and more}

%   MSC (2000 and 2010): 
%   76M28 = Particle methods and lattice-gas methods
%   35F10 = Initial value problems for linear first-order PDE, linear evolution equations
%   65M12 = Stability and convergence of numerical methods

\end{abstract}

%%% %%% %%% %%% %%% %%% %%% %%% %%% %%% %%% %%% %%% %%% %%% %%% %%% %%% %%% %%% %%% %%% %%% %%% %%% %%% %%% %%% %%% %%% %%% %%% %%% 
%%% %%% %%% %%% %%% %%% %%% %%% %%% %%% %%% %%% %%% %%% %%% %%% %%% %%% %%% %%% %%% %%% %%% %%% %%% %%% %%% %%% %%% %%% %%% %%% %%% 

\section{Introduction}
\label{sec:intro}

Efficient and simple particle methods are a very popular tool for the numerical simulation 
of transport equations involved in many physical problems, ranging from 
fluid dynamics \cite{Chorin.1973.jfm,Cottet.Koumoutsakos.2000.cup} to kinetic (e.g., Vlasov) equations \cite{Hockney.Eastwood.1988.tf,Birdsall.Langdon.2005.tf}.
However, particle methods also suffer from weak convergence properties that cause difficulties in many practical cases.
Specifically, it is known that they only converge in a strong sense when the particles present an extended overlapping,
that is, when the number of overlapping particles tends to infinity as the mesh size $h$ of their initialization grid tends to 0, see e.g.
\cite{Beale.Majda.1982b.mcomp,Raviart.1985.lnm}.
Moreover, convergence rates are known to be suboptimal and to require vanishing moments
for the particle shape functions, which prevents high orders to be achieved with positive shapes.
In practice, extended particle overlapping is expensive and it involves an additional parameter to be optimized, such as 
the exponent $q < 1$ for which the particles radius behaves like $h^q \gg h$.
In Particle-In-Cell (PIC) codes for instance, taking $q < 1$ typically leads to increasing the number of particles per cell faster than the number of cells,
since the latter determines the effective radius of the particles \cite{Birdsall.Langdon.2005.tf}. 
In Smoothed Particle Hydrodynamics (SPH) schemes it amounts to increasing
the number of neighbors, i.e., interacting particles \cite{Rasio.2000.ptps}.
Because of these issues, accurate results often require numerically intensive runs, and 
in many cases the simulations do not meet the conditions of strong convergence.
Significant oscillations are then produced, which are sometimes seen as a statistical noise that hampers interpretation of results,
and can further cause large scale errors.

To suppress noise, many methods (like the redeposition scheme introduced by Denavit~\cite{Denavit.1972.jcp} for plasma simulation 
and recently revisited as a Forward Semi-Lagrangian scheme (FSL), see e.g.~\cite{Nair.Scroggs.Semazzi.2003.jcp,Cotter.Frank.Reich.2007.rms,Crouseilles.Respaud.Sonnendrucker.2009.cpc}) 
use periodic remappings, i.e., particle re-initializations that smooth out the evolution. 
However, frequent remappings introduce numerical dissipation which in many cases contradicts 
the benefit of using non-dissipative particle schemes, and to reduce this effect many works have been devoted to the design
of accurate, low dissipative remapping schemes, including adaptive multi-resolution techniques,
see e.g.~\cite{Koumoutsakos.1997.jcp,Bergdorf.Koumoutsakos.2006.mms,Wang.Miller.Colella.2011.sicomp,VanRees.Leonard.Pullin.Koumoutsakos.2011.jcp,Magni.Cottet.2012.jcp}.
To improve the accuracy of the density evaluation without introducing unwanted smoothing, 
some authors have studied alternative methods sometimes referred to as ``pseudo-particle'' methods. 
In Beale's method~\cite{Beale.1988.cfd-minn} for instance, new weights are iteratively computed 
from the positions of the particles to evaluate the density by an approximate interpolation technique,
but particles are never remapped. And in Strain's method~\cite{Strain.1996.jcp} local interpolation formulas 
are dynamically constructed on cells that contain a given number of neighboring particles (in 
\cite{Cohen.Perthame.2000.simath} a similar idea is further investigated to establish optimal error estimates).

In this article, we present a new class of particle methods where polynomial expansions of the characteristic flow are 
locally computed to transport the numerical particles with improved accuracy. Specifically this amounts to transforming
the particle shape functions with polynomial mappings which coefficients involve local derivatives of the backward flow:
the first order version is a linearly-transformed particle (LTP) method, the second order case is a quadratically-transformed 
particle (QTP) scheme, and so on. 

On a theoretical level we prove that the resulting method converges in the uniform norm at a potential rate that depends 
on the smoothness of the flow, provided that the initial particle shape functions are Lipschitz. In particular, 
the proof does not require remappings, extended overlapping or vanishing moments for the particles.

On a practical level we provide explicit implementations of the LTP and QTP schemes that comply with the algorithmic structure of 
most existing particle codes, in the sense that they only involve pointwise evaluations of the forward characteristic flow.
%Moreover, these come with local error indicators for the transport scheme 

Numerical tests using different transport problems eventually demonstrate the accuracy of the proposed methods compared to 
standard particle schemes, and their robustness with respect to the remapping periods. In particular, it is 
shown that QTP particles can be transported without remappings on very long periods of times, without 
hampering the accuracy of the numerical solutions.
Finally, we propose a dynamic criterion to automatically select the time steps where the particles should be remapped. Our
strategy does not require additional inter-particle communications, and is validated by numerical experiments.

We note that 
%the method can be viewed as a modified version of Hou's method \cite{Hou.1990.sinum} 
%where instead of using a global deformation mapping, each particle is transported by the linearized 
%flow around its trajectory. Moreover, 
variants of this approach have been studied by some authors, 
at least to the first order. In a theoretical study indeed, Cohen and Perthame~\cite{Cohen.Perthame.2000.simath} 
observed that by transporting the particles with the linearized flow around their trajectories one would 
obtain a convergent method (in $L^1$) with particles scaled with their initialization grid, and no remappings.
And more recently, Alard and Colombi~\cite{Alard.Colombi.2005.ras} have 
described a similar scheme obtained by evolving Gaussian particles with locally affine force fields 
in PIC simulations of the Vlasov-Poisson system.

The outline of the article is as follows. In Section~\ref{sec:overview} we begin with a rapid overview of the main
particle methods and introduce some notations. In Section~\ref{sec:pmws} we present the new class of methods 
with polynomial transformations, and establish a priori convergence theorems. 
In Section~\ref{sec:discrete} we provide explicit implementations for the LTP and QTP schemes that only 
involve pointwise evaluations of the forward flow, and numerical results are eventually presented in Section~\ref{sec:num},
together with a dynamic criterion to select the remapping time steps.

%% %%% %%% %%% %%% %%% %%% %%% %%% %%% %%% %%% %%% %%% %%% %%% %%% %%% %%% %%% %%% %%% %%% %%% %%% %%% %%% %%% %%% %%% %%% %%% %%% 
%%% %%% %%% %%% %%% %%% %%% %%% %%% %%% %%% %%% %%% %%% %%% %%% %%% %%% %%% %%% %%% %%% %%% %%% %%% %%% %%% %%% %%% %%% %%% %%% %%% 

\section{A brief review of particle methods}
\label{sec:overview}

To introduce the notations and state our main results we begin with a brief review of particle methods.
Following \cite{Raviart.1985.lnm,Cohen.Perthame.2000.simath} we consider the linear %$d$-dimensional
transport equation
\begin{equation}
\label{transport}
\dt f(t, x) + u(t, x) \cdot \nabla f(t, x) = 0, \qquad t \in [0,T], \quad x \in \RR^d
\end{equation}
associated with an initial data $\fin: \RR^d \to \RR$, a final time $T$
and a velocity field $u: [0,T] \times \RR^d \to \RR^d$.
In fluid problems for instance we have $d=2,3$, while in kinetic formulations
$\RR^d$ is a phase space with $d \le 6$ and $u$ is a generalized velocity field 
with components of velocity and acceleration.
We assume that $u$ is smooth enough (e.g., Lipschitz) for the characteristic trajectories $X(t) = X(t;t_0,x_0)$, solutions to
\begin{equation}
\label{traj}
X'(t) = u(t,X(t)), \qquad X(t_0) = x_0,
\end{equation}
to be defined on $[0,T]$ for all $x_0 \in \RR^d$ and $t_0 \in [0,T]$.
In particular, the characteristic flow $F_{t_0,t}: x_0 \mapsto X(t)$ is invertible 
and satisfies $(F_{t_0,t})^{-1} = F_{t,t_0}$. Solutions to \eqref{transport} read then
\begin{equation} \label{ex-transp}
f(t,x) = f(t_0, (F_{t_0,t})^{-1}(x)) \qquad \text{ for } ~ t_0, t \in [0,T], ~  x \in \RR^d. 
\end{equation}
For simplicity, we restrict ourselves to incompressible fields satisfying $\Div u =0$. In this case the flow is measure preserving 
in the sense that its Jacobian matrix $J_{F_{t_0,t}}(x) = \big(\partial_{j} (F_{t_0,t})_{i}\big)_{1\le i,j \le d}$ satisfies 
$$
\det \big(J_{F_{t_0,t}}(x)\big) = 1 \qquad \text{ for } ~ t_0, t \in [0,T], ~ x \in \RR^d. 
$$

%%% %%% %%% %%% %%% %%% %%% %%% %%% %%%
%%% %%% %%% %%% %%% %%% %%% %%% %%% %%%
\subsection{The traditional smoothed particle method (TSP)}
\label{sec:tsp}
In the standard ``academic'' particle method \cite{Raviart.1985.lnm}, numerical solutions are 
typically computed as follows: considering deterministic initializations for simplicity,
the initial data $\fin$ is first approximated by a collection of particles on a regular (say, cartesian) grid of step $h>0$,
$$
f_{h,\epsilon}^0(x) := \sum_{k \in \ZZ^d} w_k(\fin) \vp_\epsilon (x-x^0_k)
\qquad \text{ with } \quad x^0_k := h k
$$
and with weights typically defined as
\begin{equation}
\label{stand-wk}
w_k(\fin) := \int_{x_k^0 + \left[-\frac h2,\frac h2\right]^d} \fin(x) \rmd x
\qquad \text{ or } \qquad
w_k(\fin) := h^d \fin(x_k^0).
\end{equation}
Here $\vp_\epsilon = \epsilon^{-d} \vp(\cdot / \epsilon)$ is a particle shape function with radius proportional to $\epsilon$, 
usually seen as a smooth approximation of the Dirac measure obtained by scaling a compactly supported ``cut-off''
function $\vp$ for which common choices include B-splines and smoothing kernels with vanishing moments, 
see e.g. \cite{Koumoutsakos.1997.jcp,Cottet.Koumoutsakos.2000.cup}.
Particle centers are then pushed forward at each time step $t^n = n\Dt$ by following a numerical flow $F^n$
which approximates the exact $F_{t^n,t^{n+1}}$, %and the weights are kept constant, 
leading to
$$
f_{h,\epsilon}^{n+1}(x) := \sum_{k \in \ZZ^d} w_k(\fin) \vp_\epsilon (x-x^{n+1}_k) \approx f(t^{n+1},x) 
\quad \text{ with } \quad x_k^{n+1} := F^n(x_k^{n}).
$$
In the classical error analysis \cite{Beale.Majda.1982b.mcomp,Raviart.1985.lnm}, the above process is seen as 
(i) an approximation (in the distribution sense) of the initial data by a collection of weighted Dirac measures,
(ii) the exact transport of the Dirac particles along the flow, and 
(iii) the smoothing of the resulting distribution $\sum_{k} w_k(\fin) \delta_{x^n_k}$ with the convolution kernel $\vp_\epsilon$.
The classical error estimate reads then as follows: 
if for some prescribed integers $m>0$ and $r>0$, the cut-off $\vp$ has $m$-th order smoothness and satisfies 
a moment condition of order $r$, namely if $\int \vp = 1$, $\int \abs{y}^r \abs{\vp(y)} \rmd y < \infty$ and
$$
\int y_1^{s_1} \ldots y_d^{s_d} \vp(y_1, \ldots, y_d) \rmd y = 0 
\quad \text{ for } ~ s \in \NN^d ~ \text{ with } ~  1 \le s_1 + \cdots + s_d \le r-1,
$$
then there exists a constant $C$ independent of $\fin, h$ or $\epsilon$, such that for all $n \le T/\Dt$ we have 
\begin{equation}
\label{stand-est}
\norm{f(t^n)-f_{h,\epsilon}^n}_{L^\mu} \le C \Big( \epsilon^r \norm{\fin}_{W^{r,\mu}} + (h / \epsilon)^m \norm{\fin}_{W^{m,\mu}} \Big)
\end{equation}
with $1 \le \mu \le \infty$. More recently, Cohen and Perthame \cite{Cohen.Perthame.2000.simath} observed that defining the weights as
\begin{equation}
\label{impr-wk}
w_k(\fin) := \int_{\RR^d} \fin(x) \tilde \vp_h(x-x_k^0) \rmd x
\end{equation}
with a weighting function $\tilde \vp_h = \tilde \vp(\cdot/h)$ derived from a continuous and compactly supported $\tilde \vp$ 
such that 
$$
\sum_{k\in\ZZ^d} k_1^{s_1} \ldots k_d^{s_d} \tilde \vp(y-k) = y_1^{s_1} \ldots y_d^{s_d}
\quad \text{ for } ~ s \in \NN^d ~ \text{ with } ~  0 \le s_1 + \cdots + s_d \le m-1,
$$
one has the improved estimate
\begin{equation}
\label{impr-est}
\norm{f(t^n)-f_{h,\epsilon}^n}_{L^\mu} \le C \Big( \epsilon^r \norm{\fin}_{W^{r,\mu}} + (h / \epsilon)^m \norm{\fin}_{L^\mu} \Big)
\end{equation}
with a new constant that is again independent of $\fin, h$ or $\epsilon$.
Note that \eqref{impr-est} is better than \eqref{stand-est} in that $m$ is not constrained by the smoothness of $\fin$,
which allows to reach higher convergence rates. Indeed balancing the error terms in the above estimates suggests to take 
$\epsilon \sim h^q$ with $q = {\frac {m}{m+r}}$, yielding a convergence in $h^q = h^{\frac{rm}{m+r}}$. 
In particular, if $\fin \in W^{s,\mu}$ for some integer $s$ then the best possible rate with standard weights is only 
$h^{s/2}\norm{\fin}_{W^{s,\mu}}$, obtained with $m=r=s$. With the improved weights instead, one can take a higher value for $m$ and obtain
estimates close to $h^{s}\norm{\fin}_{W^{s,\mu}}$. Moreover, the latter approach also allows to improve (i.e., reduce) 
the particle overlapping, since the corresponding exponents are $q = \frac 12$ and $\frac {m}{m+r} \approx 1$, respectively.
In either case, we see from the terms $\epsilon^r \sim h^{qr}$ in the estimates that 
extended particle overlapping does not only make the simulations more expensive, it also deteriorates their convergence order.

%%% %%% %%% %%% %%% %%% %%% %%% %%% %%%
%%% %%% %%% %%% %%% %%% %%% %%% %%% %%%
\subsection{The forward semi-Lagrangian scheme (FSL)}
\label{sec:fsl}

In forward semi-Lagrang\-ian schemes (also called remapped or remeshed particle methods), extended overlapping is usually not 
required and particles have the same scale than their initialization grid, i.e., $\epsilon = h$. Instead
they are periodically remapped on a regular grid, say every $\Nr$ time steps. Thus, letting 
$$
A_h : g \mapsto \sum_{k \in \ZZ^d} w_k(g) \vp_h (x-x_k^0)  %, ~~~ k \in \ZZ^d
$$
be an approximation operator defined on the grid, and denoting %where $x_k^0 = hk$ and with weights computed, e.g., as in \eqref{stand-wk}, 
$$
T^n_h: \vp_h(\cdot-x^{n}_k) \mapsto \vp_h(\cdot-F^n(x^{n}_k))
$$
the fixed-shape particle transport operator, FSL schemes take the generic form
%$$
%f^n_h(x) = \sum_{k \in \ZZ^d} w_k(\fin) \vp_\epsilon (x-x^n_k) \approx f(t_n,x), \qquad x_k^n = F^n(x_k^{n-1}), ~~~ k \in \ZZ^d,
%$$
\begin{equation}
\label{FSL}
f^{n+1}_h = \sum_{k\in \ZZ^d} \tilde w_{k}^{n} \vp_{h}(\cdot - F^n(\tilde x^{n}_k)) 
= T^n \tilde f^n_h
\quad \text{ with } \quad
\tilde f^n_h :=
    \begin{cases}
    ~ A_h f^n_h ~  & \text{ if } n \in \Nr \NN
    \\
    ~ f^n_h   & \text{ otherwise.}
    \end{cases}
\end{equation} 
%Although there does not seem to be a comlete error analysis of such schemes, 
In practice these schemes give satisfactory results when the remapping operators are well designed, see e.g. 
\cite{Bergdorf.Koumoutsakos.2006.mms,Wang.Miller.Colella.2011.sicomp,VanRees.Leonard.Pullin.Koumoutsakos.2011.jcp,Magni.Cottet.2012.jcp} 
for recent references.

%%% %%% %%% %%% %%% %%% %%% %%% %%% %%%
%%% %%% %%% %%% %%% %%% %%% %%% %%% %%%
%\subsection{The linearly-transformed particle method (LTP)}
%\label{sec:ltp}
%

\subsection{Particle transport with polynomial transformations}
%\label{sec:tp}

In this article we develop a lesser-known approach where particles are subject to polynomial transformations
that approximate the backward flow involved in the exact transport, see Equation~\eqref{ex-transp}. 
At the first order for instance, this amounts to formally defining linearly-transformed particles as
\begin{equation}
\label{formal-ltp}
\vp^n_{h,k}(x) := \vp_h\big(D^n_k(x-x^n_k)\big)
%
%  \quad  \text{ where } 
%\quad \begin{cases} 
%    x_k(t) := F_{0,t}(x_k^0) 
%    \\
%    J^n_k(t) := J_{F_{t,0}}(x_k(t)). %= J_{F_{0,t}}(x_k)^{-1}
%\end{cases}
\end{equation}
where in addition to pushing forward the particle centers as before, one also needs to compute $d\times d$ 
deformation matrices $D^n_k$, $k \in \ZZ^d$, that approximate the local Jacobian matrices of the backward flow.
Similarly, at the second order the particle shape functions are transformed with local quadratic mappings
which coefficients involve the derivatives of the backward flow, and so on.
Because occasional remappings are usually needed for accurate solutions we may write our particle methods in a form 
similar to \eqref{FSL}. However, as will be shown below the resulting schemes converge in the uniform norm 
without remapping the particles, and in practice this leads to improved convergence compared to fixed-shape
particle schemes, with significantly larger remapping periods $\Dtr = \Nr \Dt$.

%but particles are now associated with invertible 
%$d \times d$ deformation matrices $D_k^n$ representing backward Jacobian matrices at $x^n_k$. Thus, 
%numerical solutions read
%\begin{equation}
%\label{vpn-ltp}
%f^n_h(x) = \sum_{k \in \ZZ^d} w^n_k \vp^n_{h,k}(x)  := \sum_{k \in \ZZ^d} w^n_k \vp_h(D^n_k(x-x^n_k))
%\end{equation}
%and transporting the particles consists in updating the deformation matrices $D^n_k$ together with the particle centers $x^n_k$, 
%initialized as $D^0_k := I_d$ and $x^0_k := hk$, respectively. 
% We may summarize our findings as follows.

%\paragraph{Main results}
%The formal LTP method \eqref{formal-ltp} converges with order 1 in $L^\infty$, and arbitrary orders can be reached 
%with proper polynomial deformations which coefficients involve the derivatives of the backward flow (see Theorem~\ref{theo. cv}).
%On the numerical side, an explicit implementation of the LTP transport operator based on finite-difference approximations of the 
%forward Jacobian is also shown to converge with order 1 in $L^\infty$ with no remappings required (see Theorem~\ref{theo. ltp}). 
%In practice this leads to improved convergence compared to both TSP and FSL schemes, with lower remapping frequencies than the latter
%(see Section~\ref{sec:unif}).

In the sequel it will be convenient to use the maximum norm $\norm{x}_\infty := \max_i \abs{x_i}$
for vectors and the associated $\norm{A}_\infty := \max_i \sum_j \abs{A_{i,j}}$ for matrices.
For functions in Sobolev spaces $W^{m,\infty}(\omega)$ with $\omega \subset \RR^d$ and integer index $m > 0$, 
we will use the semi-norm
\begin{equation}\label{norm}
\abs{v}_{m,\omega} 
:= \max_{i} 
\Big\{\sum_{l_1=1}^d \cdots \sum_{l_m=1}^d %\partial_{l_1, \cdots, l_m} 
\norm{\partial_{l_1} \cdots \partial_{l_m} v_i }_{L^\infty(\omega)}\Big\},
\end{equation}
and for conciseness we will drop the domain when it is the whole space.

%%% %%% %%% %%% %%% %%% %%% %%% %%% %%% %%% %%% %%% %%% %%% %%% %%% %%% %%% %%% %%% %%% %%% %%% %%% %%% %%% %%% %%% %%% %%% %%% %%% 
%%% %%% %%% %%% %%% %%% %%% %%% %%% %%% %%% %%% %%% %%% %%% %%% %%% %%% %%% %%% %%% %%% %%% %%% %%% %%% %%% %%% %%% %%% %%% %%% %%% 

\section{Particle methods without smoothing}
\label{sec:pmws} % was {sec:nosmooth}

The particle methods that we present in this work deviate from the ``smoothed'' approaches described in Sections~\ref{sec:tsp}
and \ref{sec:fsl} in the following ways.
\begin{itemize}
\item[$\bullet$]
    Convergence (including high-order) is proved without resorting to a smoothing kernel argument.
    Instead it relies on local expansions of the flow. In particular,
    \begin{itemize}
    \item[-] 
        particles may have a radius proportional to the mesh-size $h$ of their initialization grid, so that
        their overlapping is uniformly bounded with respect to $h$~;
    \item[-]         
        many particle shape functions can be used, including heterogeneous particle collections such as standard finite element bases,
        as no vanishing moments or high-order smoothness is required~;
    \item[-] 
        when initialized or remapped, particle weights can be computed with standard approximation schemes.
    \end{itemize}
\item[$\bullet$]
    Remappings are not required for convergence as in FSL schemes. In practice they may improve the results, but at a frequency
    that is significantly lower than what is needed with FSL schemes.
\end{itemize}

Specifically, we will address the problem of transporting a collection of particles
\begin{equation} \label{f0h}
f^0_h = \sum_{k \in \ZZ^d} w_{k} \vp^0_{h,k} \approx \fin
\end{equation}
with a particle-wise operator $\bar T^n_h$ that approximates the exact transport associated with the flow $\bar F^n_\ex = F_{0,t^n}$,
\begin{equation}
\label{Tex}
\bar T^n_\ex : \vp^0_{h,k} \mapsto \vp^0_{h,k}(\bar B^n_\ex(\cdot)) \quad \text{ where } \quad \bar B^n_\ex = (\bar F^n_\ex)^{-1}
\end{equation}
(here and below we use a bar to distinguish flows and transport operators defined on the global time interval $[0,t^n]$ from those 
on the single time step $[t^n,t^{n+1}]$).
For simplicity we consider homogeneous particles based on a reference shape function $\vp$ 
that is supported inside the $d$-dimensional cube $B_{\ell^\infty}(0,\rho^0)$ (see Section~\ref{sec:shapes} for examples). 
Initially our particles are centered on the cartesian nodes $x^0_k := hk$,
$k \in \ZZ^d $, they are scaled to the grid and normalized in $L^1$, i.e.,
\begin{equation}
\label{vphk}
\vp_{h,k}^0(x) := \vp_{h}(x - x^0_k) := h^{-d} \vp (h^{-1} x - k).
\end{equation}
In particular, the particles are initially supported in small $d$-dimensional cubes
\begin{equation} \label{supp-vp0}
\Sigma^0_{h,k}  :=  \supp(\vp_{h,k}^0)  \subset  B_{\ell^\infty}(x^0_k,h \rho^0)
\end{equation}
and consistent with \eqref{f0h} and \eqref{vphk} we may assume that the weights satisfy
\begin{equation}
\label{wk-bound}
\abs{w_k} \le c_w h^d \norm{\fin}_{L^\infty}, \qquad k \in \ZZ^d.
\end{equation}

%\todo{check outline}
%In Section~\ref{sec:bspline} we first review one local approximation scheme with B-splines.
%In Section~\ref{sec:polytrans} we introduce particle transport operators $T_{(r)}$ 
%with polynomial shape transformations that converge without smoothing arguments,
%and in Section~\ref{sec:ltp-fd} we describe one explicit implementation of the first-order operator $T_{(1)}$ . 
%In Section~\ref{sec:local} we give a practical tool for localizing particles with linearly transformed supports, and in Section~\ref{sec:error}
%we establish rigorous error estimates that prove both the uniform convergence and the uniformly bounded overlapping of the deformed particles.
%

%%% %%% %%% %%% %%% %%% %%% %%% %%% %%%
%%% %%% %%% %%% %%% %%% %%% %%% %%% %%%
\subsection{Two preliminary estimates}
\label{sec:pre-est}

The general form considered in this paper for particle transport operators is
\begin{equation} \label{T-with-bound}
\bar T^n_h \vp^0_{h,k} (x) := \chi_{\Sigma^n_{h,k}}(x) \vp^0_{h,k}(\bar B^n_{h,k}(x)).
\end{equation}
Here $\bar B^n_{h,k} \approx \bar B^n_\ex$ is the approximated backward flow for the $k$-th particle,
$\chi$ is the set characteristic function and $\Sigma^n_{h,k}$ is an a priori bound for the particle support.
A priori bounds are needed when the domains $(\bar B^n_{h,k})^{-1}(\Sigma^0_{h,k})$ are not easily computable,
or when they are very large compared with $\bar F^n_\ex(\Sigma^0_{h,k})$. For linearly-transformed particles this will not be the case 
and the particle transport will have a simpler definition, see \eqref{T-no-bound} below.

It will be useful to state a preliminary estimate for the convergence of such particle methods, 
based on the overlapping constant
\begin{equation} \label{over-Theta}
\Theta^n(h) := \sup_{x \in \RR^d} \#\big( \{k \in \ZZ^d : x \in \Sigma^n_{h,k} \} \big)
\end{equation}
and on the backward flow error
\begin{equation} \label{eB}
e^n_B(h) := \sup_{k \in \ZZ^d} \norm{\bar B^n_{h,k}-\bar B^n_\ex}_{L^\infty(\Sigma^n_{h,k})}.
\end{equation}

\begin{lemma}
\label{lem:conv-f-gal}
If the exactly transported particles vanish outside the domains $\Sigma^n_{h,k}$,
\begin{equation} \label{hyp-domain}
\bar F^n_\ex(\Sigma^0_{h,k}) \subset \Sigma^n_{h,k}, \qquad k \in \ZZ^d,
\end{equation}
then the approximate transport operator \eqref{T-with-bound} satisfies
$$
\norm{(\bar T^n_h - \bar T^n_\ex) f^0_h}_{L^\infty} \lesssim \frac{e^n_B(h) \Theta^n(h) }{h} \norm{f^0}_{L^\infty}
$$
with a constant independent of $h$, $f^0$, $u$ and $n$.
\end{lemma}

\begin{proof}
For $x \in \RR^d$, we let
%\begin{equation} \label{bar-cK}
$\cK^n_h(x) := \{k \in \ZZ^d : x \in \Sigma^n_{h,k} \}$
%\end{equation}
and infer from \eqref{hyp-domain} that 
$$
\bar T^n_h \vp^0_{h,k}(x) - \vp^0_{h,k}(\bar B^n_\ex (x))
=
\begin{cases}
    0  
        &\text{ if } ~~ k \notin  \cK^n_h(x)
    \\
    \vp^0_{h,k}(\bar B^n_{h,k} (x)) - \vp^0_{h,k}(\bar B^n_\ex (x))
        &\text{ otherwise. }
\end{cases}
$$
It follows that
\begin{align*}
\abs{(\bar T^n_h - \bar T^n_\ex) f^0_h(x)}
    &= \Bigabs{\sum_{k \in \ZZ^d} w_k \bar T^n_h \vp^0_{h,k}(x) - f^0_h(\bar B^n_\ex (x))} 
\\
    &= \Bigabs{ \sum_{k \in \cK^n_h(x)} w_k \big( \vp^0_{h,k}(\bar B^n_{h,k} (x)) - \vp^0_{h,k}(\bar B^n_\ex (x))\big)}
\\
    &\le \sum_{k \in \cK^n_h(x)} 
            \abs{w_k} \abs{\vp^0_{h,k}}_1 \norm{\bar B^n_{h,k} - \bar B^n_\ex}_{L^\infty(\Sigma^n_{h,k})}
\\
    &\le C \norm{f^0}_{L^\infty} h^{-1} \Theta^n(h) e^n_B(h),
\end{align*}
where we have used the bound \eqref{wk-bound} and the scaling $\abs{\vp^0_{h,k}}_1 \sim h^{-1-d}$.
\qed
\end{proof}

In the case where one can define approximate forward flows $\bar F^n_{h,k} := (\bar B^n_{h,k})^{-1}$ that are easily computable 
(for instance if the $\bar B^n_{h,k}$ are affine mappings), the particle transport operator can be simplified 
as we no longer need a priori bounds for the particle supports. Indeed, letting
\begin{equation} \label{T-no-bound}
\bar T^n_h \vp^0_{h,k} (x) := \vp^0_{h,k}(\bar B^n_{h,k}(x))
\end{equation}
readily defines particles supported in the domains $\bar F^n_{h,k}(\Sigma^0_{h,k})$, $k \in \ZZ^d$.
And if we now {\em denote}
\begin{equation} \label{bSigma}
\Sigma^n_{h,k} := \bar F^n_\ex(\Sigma^0_{h,k}) \cup \bar F^n_{h,k}(\Sigma^0_{h,k}),
\end{equation}
then the definitions \eqref{T-no-bound} and \eqref{T-with-bound} are equivalent,
moreover Assumption~\eqref{hyp-domain} is readily fulfilled.
It is then possible to establish an a priori bound for the corresponding overlapping constant \eqref{over-Theta} 
and the transport error, that is either based on the above backward flow error \eqref{eB}
or on the forward flow error
\begin{equation} \label{eF}
e^n_F(h) := \sup_{k \in \ZZ^d} \norm{\bar F^n_{h,k}-\bar F^n_\ex}_{L^\infty(\Sigma^0_{h,k})}.
\end{equation}

\begin{lemma}
\label{lem:conv-f-aff}
The approximate transport operator \eqref{T-no-bound} satisfies
\begin{equation}\label{conv-f-aff}
\norm{(\bar T^n_h - \bar T^n_\ex) f^0_h}_{L^\infty} \lesssim \Big(1+\frac{e^n_B(h)}{h}\Big)^d \frac{e^n_B(h)}{h} \norm{f^0}_{L^\infty}
\end{equation}
with a constant independent of $h$, $f^0$, $u$ and $n$. Moreover, if the exact and approximate backward flows satisfy
uniform Lipschitz estimates
\begin{equation}\label{Bhk-Lip}
\abs{\bar B^n_\ex}_{1},  \sup_{k \in \ZZ^d} \abs{\bar B^n_{h,k}}_{1,\Sigma^n_{h,k}} \le C,
\end{equation}
then the transport error is also controlled by the forward flow error,
\begin{equation}\label{conv-f-F}
\norm{(\bar T^n_h - \bar T^n_\ex) f^0_h}_{L^\infty} \lesssim \Big(1+\frac{e^n_F(h)}{h}\Big)^d \frac{e^n_F(h)}{h} \norm{f^0}_{L^\infty}.
\end{equation}
\end{lemma}

\begin{proof}
For $x \in \RR^d$, we now denote 
$$
\cK^n_\ex(x) := \{k \in \ZZ^d : x \in \bar F^n_\ex(\Sigma^0_{h,k})\}
\quad \text{ and } \quad 
\cK^n_h(x) := \{k \in \ZZ^d : x \in \bar F^n_{h,k}(\Sigma^0_{h,k})\}.
$$
The cardinality of $\cK^n_\ex(x)$ is readily bounded by the overlapping of the initial supports 
$\Sigma^0_{h,k}$: from \eqref{supp-vp0} we find indeed
$$
\#\big(\cK^n_h(x)\big) \le (2 \rho^0)^d.
$$
Moreover, for $k \in \cK^n_h(x)$ we write
$$
\norm{hk - \bar B^n_\ex(x)}_\infty \le \norm{hk - \bar B^n_{h,k}(x)}_\infty + \norm{\bar B^n_{h,k}(x) - \bar B^n_\ex(x)}_\infty
< h\rho^0 + e^n_B(h)
$$
and since $\cK^n_h(x)$ is a subset of $\ZZ^d$, the above bound yields
$$
\#\big(\cK^n_h(x)\big) \le \big(2 (\rho^0 + h^{-1}e^n_B(h))\big)^d.
$$
It follows that the overlapping constant \eqref{over-Theta} is bounded by
\begin{equation} \label{Theta-bound}
\Theta^n(h) \le \sup_{x\in\RR^d} \big(\#(\cK^n_h(x)) + \#(\cK^n_\ex(x)) \big) \lesssim (1 + h^{-1} e^n_B(h))^d,
\end{equation}
so that Lemma~\ref{lem:conv-f-gal} gives the first estimate \eqref{conv-f-aff}.

In order to derive an estimate based on the forward flow error we next write for $k \in \cK^n_\ex(x)$ that
\begin{align*}
\norm{\bar B^n_{h,k}(x) - \bar B^n_\ex(x)}_\infty 
    &\le \norm{\bar B^n_{h,k}\big(\bar F^n_\ex(\bar B^n_\ex (x))\big) - \bar B^n_{h,k}\big(\bar F^n_{h,k}(\bar B^n_\ex (x))\big)}_\infty
\\
    &\le \abs{\bar B^n_{h,k}}_{1,\Sigma^n_{h,k}} \norm{\bar F^n_\ex(\bar B^n_\ex(x)) - \bar F^n_{h,k}(\bar B^n_\ex(x))}_\infty
\\
    &\le \abs{\bar B^n_{h,k}}_{1,\Sigma^n_{h,k}} e^n_F(h)
\end{align*}
where we have used that $\bar B^n_\ex(x) \in \Sigma^0_{h,k}$ in the last two inequalities.
Similarly, for $k \in \cK^n(x)$ we write
\begin{align*}
\norm{\bar B^n_{h,k}(x) - \bar B^n_\ex(x)}_\infty 
    &\le \norm{\bar B^n_{h,k}(x) - \bar B^n_\ex\big(\bar F^n_{h,k}(\bar B^n_{h,k}(x))\big)}_\infty
\\
    &\le \abs{\bar B^n_\ex}_{1} \norm{\bar F^n_\ex(\bar B^n_{h,k}(x)) - \bar F^n_{h,k}(\bar B^n_{h,k}(x))}_\infty
\\
    &\le \abs{\bar B^n_\ex}_{1} e^n_F(h)
\end{align*}
where we have now used that $\bar B^n_{h,k}(x) \in \Sigma^0_{h,k}$. According to \eqref{eB}
we thus have
$$
e^n_B(h) \le  \max\big\{\abs{\bar B^n_\ex}_{1}, \sup_{k \in \ZZ^d} \abs{\bar B^n_{h,k}}_{1,\Sigma^n_{h,k}}\big\} \, e^n_F(h)
$$ 
which gives the desired estimate.
\qed
\end{proof}

%%% %%% %%% %%% %%% %%% %%% %%% %%% %%%
%%% %%% %%% %%% %%% %%% %%% %%% %%% %%%
\subsection{Particle approximations using a cartesian grid}
\label{sec:shapes}

Several choices can be made for the reference shape function $\vp$ used in the particle definition \eqref{vphk}.
One standard option consists of taking the interpolating kernel $M'_4$ introduced by Monaghan \cite{Monaghan.1985.jcp}, i.e.,
\begin{equation} \label{vp-M4}
\vp(x) = \prod_{i=1}^d M'_4(x_i)
\quad \text{ with } \quad
M'_4(x_i) = \begin{cases}
    1 - 5 \frac{\abs{x_i}^2}{2} + 3 \frac{\abs{x_i}^3}{2} ~ &\text{ if } 0 \le \abs{x_i} \le 1 
    \\
    \frac 12 (2 - \abs{x_i})^2 (1 - \abs{x_i}) ~ &\text{ if } 1 \le \abs{x_i} \le 2 
    \\
    0  \quad &\text{ otherwise.}
\end{cases}
\end{equation}
%and of defining $\vp$ with a tensor product,
%\begin{equation} \label{vp-M4}
%\vp(x) := \prod_{i=1}^d M'_4(x_i).
%\end{equation}
Particles are then initialized and remapped on the grid by a simple call to the approximation operator
\begin{equation}
\label{Ah-M4}
A_h : g \mapsto \sum_{k \in \ZZ^d} w_{k}(g) \vp^0_{h,k}
~~  \text { with } ~~ 
w_{k}(g) := h^{d} g(x^0_k).
\end{equation}
Here $A_h$ is an interpolation that reproduces second-degree polynomials, hence it is third order accurate.

Another option consists of using cardinal B-splines, defined recursively with
$$
\cB_0 := \chi_{[-\tfrac 12, \tfrac 12]} ~~ \text{ and } ~~
\cB_p(x) := (\cB_{p-1} * \cB_0) (x) = \int_{x-\tfrac 12}^{x+\tfrac 12} \cB_{p-1} \quad \text{ for } ~~ p \ge 1,
$$
%Clearly $\cB_p$ is a piecewise polynomial of degree $p$, it is globally $\cC^{p-1}$ 
%and is supported on $[-\rho^0,\rho^0]$ with $\rho^0  =\frac{p+1}{2}$. Moreover its integer translations 
%span the space of cardinal splines of degree $p$, see, e.g., \cite{de-Boor.1978.sv}.
so that $\cB_1(x) = \max\{1-\abs{x},0\}$ is the traditional ``hat-function'', $\cB_3$ is the centered cubic B-spline, and so on.
For the reference particle shape function we then take the tensorized B-spline,
\begin{equation}
\label{vp-Bp}
\vp(x) := \prod_{i=1}^d \cB_p(x_i) \quad \text{ supported on } ~~ \supp(\vp) =
B_{\ell^\infty}(0,\rho^0)
~~ \text{ with } ~ \rho^0 := \tfrac{p+1}{2}.
\end{equation}
For the initialization and remappings we can then use standard approximation schemes that rely on the fact that the span 
of their integer translates \eqref{vphk} contains the space $\cP_p$ of polynomials 
with coordinate degree less or equal to $p$.
Specifically, we can use the quasi-interpolation schemes described by 
\cite{Chui.Diamond.1990.numat} and \cite{Unser.Daubechies.1997.ieee-tsp}, 
where high-order B-spline approximations are locally obtained by pointwise evaluations of the target function. 
The resulting approximation reads
\begin{equation}
\label{Ah}
A_h : g \mapsto \sum_{k \in \ZZ^d} w_{k}(g) \vp^0_{h,k}
~~  \text { with } ~~ 
w_{k}(g) := h^{d} \sum_{\norm{l}_\infty \le m_p} a_{l} \, g(x^0_{k+l}),
\quad 
a_{l} := \prod_{1 \le i \le d}a_{l_i}
\end{equation}
with symmetric coefficients $a_l=a_{-l}$ computed with the iterative algorithm in \cite[Section 6]{Chui.Diamond.1990.numat}: 
for the first odd orders we obtain
\begin{itemize}
\item[$\bullet$]
    $m_p = 0$ and $a_0 = 1$ for $p=1$,
\item[$\bullet$]
    $m_p = 1$ and $(a_0,a_1) = (\frac 86, -\frac 16)$ for $p=3$,
\item[$\bullet$]
    $m_p = 4$ and $(a_0,a_1,a_2,a_3,a_4) = (\frac{503}{288}, -\frac{1469}{3600},\frac{7}{225},\frac{13}{3600},\frac{1}{14400})$ for $p=5$.
\end{itemize}
The resulting approximation reproduces polynomials of order $p$, hence it is of order $p+1$: for all $q \le p$ we have
\begin{equation}
\label{Ae-errestim}
\norm{A_h g - g}_{L^\infty} \le  c_A h^{q+1} \abs{g}_{q+1}
%\quad \text{ with } \quad
%c_A = (\norm{A_h}_{{\mathcal L}(L^\infty)} +1) \frac{(m_p+\rho^0)^{q+1}}{(q+1)!}.
\end{equation}
for some constant $c_A$ that is independent of $h$.

In Section~\ref{sec:num} we will show numerical results obtained by applying our transport schemes to either $M'_4$ particles defined
as in \eqref{vp-M4}, or cubic spline particles defined as in \eqref{vp-Bp} with $p=3$. We already note that the latter choice
results in remappings that are fourth-order, but are also more numerically dissipative due to the wider stencil.

\subsection{Particle transport with polynomial shape transformations}
\label{sec:polytrans} % was {sec:hot}

In this section we assume that the exact flow $\bar F^n_\ex$ is known and can be applied exactly, as well as its derivatives
(in Section~\ref{sec:discrete} this assumption will be relaxed). 
Thus, in the traditional method the particles keep their shape and are simply translated with
\begin{equation}
\label{T0}
\bar T^n_{h,(0)}\vp_{h,k}^0 (x) = \vp^0_{h,k}(\bar B^n_{(0),k}(x)) = \vp_h(x - \bar F^n_\ex(x^0_{k})),
%\bar x_{k} := \bar F^n_\ex(x^0_{k}).
\end{equation}
which corresponds to approximating the exact backward flow with
$$
\bar B^n_{(0),k}(x) := x - \bar F^n_\ex(x^0_{k}) + x^0_{k}.
$$
For point (Dirac) particles this operator coincides with \eqref{Tex} and is exact. For finite-size particles however, 
the method does not converge in general.
Assume indeed that $\vp$ is the hat function, and consider the smooth 2d problem where $\fin = 1$ 
and $u(t,x) = (-x_2,x_1)$ over the time interval $[0,\tfrac{\pi}{4}]$. 
Then any reasonable initialization will give $w_k = h^2$, hence, $f^0_h(x) = 1$, and clearly the exact final solution is $f(\tfrac \pi 4,x) = 1$.
Now, at the final time the particle centers will have rotated of $T=\tfrac{\pi}{4}$, therefore every particle with $\abs{k_1}+\abs{k_2}=1$ 
will be centered on $(\cos(\theta+\tfrac{\pi}{4}),\sin(\theta+\tfrac{\pi}{4}))$ with $\theta \in \tfrac {\pi}{2} \NN$, and hence contributes
to $x=0$ with $\bar T^n_{h,(0)} \vp^0_{h,k}(0) = h^{-2}(1-\frac{1}{\sqrt 2})^2$, in addition to $\vp_{h,0}$ which does not move. 
Since the other particles do not contribute to $x=0$, the final error satisfies
$$
\norm{(\bar T^n_{h,(0)} - \bar T^n_\ex)f_h^0}_{L^\infty} \ge  \abs{\bar T^n_{h,(0)} f^0_h(0) - 1} = 2({\sqrt 2}-1)^2 , \qquad \text{ regardless of $h$}.
$$ 
To improve the accuracy of the transport operator, the error estimates in Section~\ref{sec:pre-est} suggest to use higher-order approximations
of the backward flow. Letting indeed
\begin{equation}
\label{phi}
\phi_{k}(s) = \phi_{k}(s;x):= (\bar B^n_\ex-I)(\bar F^n_\ex(x^0_{k}) + s(x - \bar F^n_\ex(x^0_{k}))), % = (I-F)(x'+s(\bar x_\gamma - x))
\end{equation}
we see that the approximation $\bar B^n_{(0),k}(x) \approx \bar B^n_\ex(x)$
corresponds to the lowest-order expansion $\phi_{k}(0) \approx \phi_{k}(1)$.
We may then consider $r$-th degree expansions, $r \ge 1$,
%replacing $\vp_{h}(x- \bar x_{k})$ by $\vp_h\big( \Phi_{k,(r)} (x) \big)$, $r \ge 1$, where
%\begin{equation}
%\label{Phi}
%\Phi_{k,(r)}(x) := x - \bar x_{k} + \phi'_{k}(0) + \cdots + \frac{1}{r!}\phi^{(r)}_{k}(0)
%\approx \Phi_{k,\ex}(x) := F^{-1}(x)-x_{k}^0
%\end{equation}
\begin{equation}
\label{Bhk-r}
\bar B^n_{(r),k}(x) := x - \bar F^n_\ex(x^0_{k}) + x^0_k + \phi'_{k}(0) + \cdots + \frac{1}{r!}\phi^{(r)}_{k}(0)
\approx \bar B^n_\ex(x).
\end{equation}
%is formally an $r$-th order approximation. 
Here we could have used the alternate $\tilde \phi_{k}(s) := (I-\bar F^n_\ex)(x_k^0+s(\bar B^n_\ex(x) - x_k^0))$ 
since $\tilde \phi_{k}(1)-\tilde \phi_{k}(0) = \phi_{k}(1)-\phi_{k}(0)$, but we observe that the form of \eqref{phi} gives
\begin{equation}
\label{der-phi}
\phi^{(r)}_{k}(s) = \!\! \sum_{l_1 \cdots l_r=1}^d  
\Big[\partial_{l_1} \cdots \partial_{l_r} (\bar B^n_\ex-I)(\bar F^n_\ex(x^0_{k}) + s(x-\bar F^n_\ex(x^0_{k})))
\prod_{i=1}^r (x-\bar F^n_\ex(x^0_{k}))_{l_i} \Big]
\end{equation}
so that $\bar B^n_{(r),k}$ is a polynomial mapping which coefficients involve derivatives of $\bar B^n_\ex$ 
at $\bar F^n_\ex(x^0_{k})$, which can be written in terms of the derivatives of $\bar F^n_\ex$ at $x^0_k$.
Moreover, \eqref{der-phi} allows to specify the accuracy of the Taylor expansions \eqref{Bhk-r}.
Indeed for every $x$ in a localized domain $\omega \subset B_{\ell^\infty}(\bar F^n_\ex(x^0_{k}),h \rho)$ with $\rho > 0$, we have
\begin{equation}    
\label{acc-Taylor}
\norm{\bar B^n_{(r),k}(x)-\bar B^n_\ex(x)}_\infty 
=  \Bignorm{\int_0^1 \frac{(1-s)^r}{r!} \phi_{k}^{(r+1)}(s)\rmd s}_\infty
\le h^{r+1}\frac{\rho^{r+1}}{(r+1)!} \abs{\bar B^n_\ex}_{r+1,\langle \omega \rangle}
%\end{array}
\end{equation}
where $\langle \omega \rangle$ denotes the convex hull of $\omega$.
%\newline \newline
\bigskip

For $r=1$, observing that $J_{\bar B^n_\ex}(\bar F^n_\ex(x_{k}^0)) = (J_{\bar F^n_\ex}(x^0_{k}))^{-1}$ we obtain
\begin{equation} \label{Bhk-1}
\bar B^n_{(1),k}(x) = x^0_k + (J^n_k)^{-1}(x-\bar F^n_\ex(x^0_k))
\quad \text{ with } \quad 
J^n_k := J_{\bar F^n_\ex}(x^0_{k}),
\end{equation}
so that the linearly-transformed particle (LTP) transport operator is
\begin{equation} \label{T1}
\bar T^n_{h,(1)}\vp^0_{h,k}(x) := \vp^0_{h,k}( \bar B^n_{(1),k}(x) ) = \vp_h \big( (J^n_k)^{-1}(x-\bar F^n_\ex(x^0_{k})) \big).
\end{equation}
We observe that this corresponds to using for the $k$-th particle the {\em exact} transport operator 
$\bar T^n_\ex[\bar F^n_{(1),k}]$ associated with the linearized flow at $x^0_k$,
\begin{equation} \label{Fhk-1}
\bar F^n_{(1),k}(\hat x) := (\bar B^n_{(1),k})^{-1}(\hat x) = \bar F^n_\ex(x^0_k) + J^n_k(\hat x-x^0_k).
\end{equation} 
We are thus in position to use Lemma~\ref{lem:conv-f-aff}.

\begin{theorem} \label{th:conv-f-T1}
The LTP transport operator \eqref{T1} satisfies
\begin{equation}\label{conv-f-T1}
\norm{(\bar T^n_{h,(1)} - \bar T^n_\ex) A_h f^0}_{L^\infty} \lesssim h c^{n}_{F} \big(1+h c^{n}_{F}\big)^d \norm{f^0}_{L^\infty}
\end{equation}
with $c^{n}_{F} := \abs{\bar F^n_\ex}_1^2 \abs{\bar B^n_\ex}_2$, and an unspecified constant depending only on $p$ and $d$.
\end{theorem}

\begin{proof}
Applying Lemma~\ref{lem:conv-f-aff}, we obtain
\begin{equation} \label{errest-eB-T1}
\norm{(\bar T^n_{h,(1)} - \bar T^n_\ex) A_h f^0}_{L^\infty} \lesssim \Big(1+\frac{e_{B,(1)}^n(h)}{h}\Big)^d \frac{e_{B,(1)}^n(h)}{h} \norm{f^0}_{L^\infty}
\end{equation}
where we have set
\begin{equation} \label{eB-1}
e_{B,(1)}^n(h) := \sup_{k \in \ZZ^d} \norm{\bar B^n_{(1),k}-\bar B^n_\ex}_{L^\infty( \Sigma^n_{h,k})}
\end{equation}
and $\Sigma^n_{h,k} := \bar F^n_{(1),k}(\Sigma^0_{h,k}) \cup \bar F^n_\ex(\Sigma^0_{h,k})$.
Next from \eqref{supp-vp0} one easily sees that the transported particles are supported on
$$
%\supp(\bar T^n_{h,(1)}\vp^0_{h,k}) = 
\bar F^n_{(1),k}(\Sigma^0_{h,k}) 
= \bar F^n_\ex(x^0_{k}) +  J^n_k(B_{\ell^\infty}(0,h \rho^0)) \subset B_{\ell^\infty}(\bar F^n_\ex(x_k^0), h \rho^0 \norm{J^n_k}_\infty)
$$
%with $\rho^0 \norm{J^n_k}_\infty = \rho^0 \norm{J_{\bar F^n_\ex}(x^0_{k})}_\infty \le \rho^0 \abs{\bar F^n_\ex}_{1}$.
Moreover, the supports of the exactly transported particles satisfy
\begin{equation} \label{supp-Tex}
%\supp(\bar T^n_\ex \vp^0_{h,k}) = 
\bar F^n_\ex(\Sigma^0_{h,k}) = \bar F^n_\ex(B_{\ell^\infty}(x_k^0,h \rho^0)) 
\subset B_{\ell^\infty}(\bar F^n_\ex(x_k^0), h \rho^0 \abs{\bar F^n_\ex}_{1}).
\end{equation}
Using next $\norm{J^n_k}_\infty \le \abs{\bar F^n_\ex}_{1}$ we obtain
% the set $\Sigma^n_{h,k}$ is inside $B_{\ell^\infty}(\bar F^n_\ex(x_k^0), h \rho^0 \abs{\bar F^n_\ex}_{1})$,
$\Sigma^n_{h,k} \subset B_{\ell^\infty}(\bar F^n_\ex(x_k^0), h \rho^0 \abs{\bar F^n_\ex}_{1})$
and estimate \eqref{acc-Taylor} yields
\begin{equation} \label{estim-eB-1}
e_{B,(1)}^n(h) \le h^2 \frac{(\rho^0\abs{\bar F^n_\ex}_1)^2}{2} \abs{\bar B^n_\ex}_2 \lesssim h^2 c^{n}_{F}
\end{equation}
with a constant depending only on $p$. This completes the proof.
\qed
\end{proof}

\begin{remark}
From \eqref{Theta-bound} we also see that the particles transported with the LTP operator \eqref{T1} 
have a bounded overlapping constant 
\begin{equation}
\sup_{x \in \RR^d} \#\big( \{k \in \ZZ^d : \bar T^n_{h,(1)}\vp^0_{h,k}(x) \neq 0 \} \big)  \lesssim (1 + hc^{n}_{F})^d
\end{equation}
with again $c^{n}_{F} := \abs{\bar F^n_\ex}_1^2 \abs{\bar B^n_\ex}_2$ and a constant depending only on $p$ and $d$.
\end{remark}
    
For $r > 1$ the approximate flow \eqref{Bhk-r} is in general not invertible,
hence there is no guarantee that the support of $\vp^0_{h,k}(\bar B^n_{(r),k}(\cdot))$ is contained in a ball of radius $Ch$.
In our simulations we have observed that this feature caused strong oscillations in the numerical density,
so that it was indeed necessary to restrict the transported particles on a priori domains as in \eqref{T-with-bound}.
Now, although for practical efficiency the particle supports need to be defined with care, as will be seen in 
Section~\ref{sec:qtp}, for proving convergence rates one may simply define them by transporting 
small extensions of the initial supports with the linearized forward flow. Specifically, we consider (using a tilde to 
distinguish this set from the one in \eqref{eB-1})
\begin{equation} \label{bS-n}
\tilde \Sigma_{h,k}^n := \bar F^n_{(1),k}(B_{\ell^\infty}(x^0_k,h \tilde \rho^n_{h,k}))
~~  \text{ with } ~~
\tilde \rho^n_{h,k} := \rho^0 + \frac 1h \norm{\bar B^n_{(1),k}-\bar B^n_\ex}_{L^\infty(\bar F^n_\ex(\Sigma^0_{h,k}))}.
\end{equation}

\begin{theorem} \label{th:conv-f-Tr}
For $r \ge 1$, the $r$-th order transport operator defined by 
\begin{equation} \label{Tr}
\bar T^n_{h,(r)} \vp^0_{h,k}(x) := \chi_{\tilde \Sigma^n_{h,k}}(x) \vp^0_{h,k}(\bar B^n_{(r),k}(x))
\end{equation}
satisfies
\begin{equation}\label{conv-f-Tr}
\norm{(\bar T^n_{h,(r)} - \bar T^n_\ex) A_h f^0}_{L^\infty} \lesssim h^r \tilde c^{n}_{F,(r)} (1+h\tilde c^{n}_{F,(1)})^d \norm{f^0}_{L^\infty}
\end{equation}
with an unspecified constant independent of $h$, $f^0$, $u$ and $n$, and where for $r \ge 1$ we have set
$
\tilde c^{n}_{F,(r)} = (\tilde \rho^n\abs{\bar F^n_\ex}_1)^{r+1}\abs{\bar B^n_\ex}_{r+1}
$
and $\tilde\rho^n = \rho^0 (1+\tfrac{h\rho^0}{2} \abs{\bar F^n_\ex}_{1}^2 \abs{\bar B^n_\ex}_2)$.
\end{theorem}

\begin{proof}
We first check that the support \eqref{bS-n} contains $\bar F^n_\ex(\Sigma^0_{h,k})$. To do so we take $x = \bar F^n_\ex(\hat x)$ with
$\hat x \in \Sigma^0_{h,k}$ and write
$$
\norm{\bar B^n_{(1),k}(x)-x^0_k}_\infty \le 
\norm{(\bar B^n_{(1),k}-\bar B^n_\ex)(x)}_\infty + \norm{\hat x-x^0_k}_\infty
\le h \tilde\rho^n_{h,k}.
$$
This shows that $x \in \tilde \Sigma^n_{h,k}$, hence Assumption~\eqref{hyp-domain} is fulfilled indeed and Lemma~\ref{lem:conv-f-gal} 
applies: For some constant independent of $h$, $f^0$, $u$ and $n$, we have
\begin{equation} \label{errest-eB*-Tr}
\norm{(\bar T^n_{h,(r)} - \bar T^n_\ex) A_h f^0}_{L^\infty} \lesssim \frac{\tilde e_{B,(r)}^{n}(h) \tilde \Theta^n(h) }{h} \norm{f^0}_{L^\infty}
\end{equation}
with an overlapping constant defined similarly as in \eqref{over-Theta}, i.e.,
\begin{equation} \label{over-Theta-*}
\tilde \Theta^n(h) := \sup_{x \in \RR^d} \#\big( \{k \in \ZZ^d : x \in \tilde \Sigma^n_{h,k} \} \big)
\end{equation}
and a backward flow error defined similarly as in \eqref{eB}, i.e.,
$$ %\begin{equation} \label{eB-r*}
\tilde e_{B,(r)}^{n}(h)  := \sup_{k \in \ZZ^d} \norm{\bar B^n_{(r),k}-\bar B^n_\ex}_{L^\infty(\tilde \Sigma^n_{h,k})}.
$$ %\end{equation}
To further bound the overlapping constant we proceed similarly as in the proof of Lemma~\ref{lem:conv-f-aff}:
given $x$ and $k$ such that $x \in \tilde \Sigma^n_{h,k}$, we write
\begin{align*}
\norm{hk - \bar B^n_\ex(x)}_\infty  
    &\le \norm{x^0_k - \bar B^n_{(1),k}(x)}_\infty + \norm{\bar B^n_{(1),k}(x) - \bar B^n_\ex(x)}_\infty
\\
    &< h \tilde\rho^n_{h,k} + \tilde e_{B,(1)}^{n}(h)
    \le h \rho^0 + 2\tilde e_{B,(1)}^{n}(h),
\end{align*}
and using that $k \in \ZZ^d$ we find
\begin{equation} \label{errest-Theta*}
\tilde \Theta^n(h) \le \Big( 2\rho^0 + 4 \frac{\tilde e_{B,(1)}^{n}(h)}{h}\Big)^d.
\end{equation}
It remains to estimate the flow errors.
From \eqref{bS-n} and \eqref{eB-1}-\eqref{estim-eB-1} we first derive
\begin{equation} \label{bound-rho*hk}
\tilde\rho^n_{h,k} %= \rho^0 + \tfrac {\norm{\bar B^n_{(1),k}-\bar B^n_\ex}_{L^\infty(\bar F^n_\ex(\Sigma^0_{h,k}))}}{h}
\le \rho^0 + \frac{e_{B,(1)}^n(h)}{h} 
\le \rho^0 + \frac{h}{2} (\rho^0\abs{\bar F^n_\ex}_1)^2 \abs{\bar B^n_\ex}_2  = \tilde\rho^n.
\end{equation}
We then observe that any $x \in \tilde \Sigma^n_{h,k}$ reads 
$x = \bar F^n_{(1),k}(\hat x) = \bar F^n_\ex(x^0_k) + J^n_k(\hat x-x^0_k)$ for some $\hat x \in B_{\ell^\infty}(x^0_k,h\tilde\rho^n_{h,k})$. 
Using \eqref{bound-rho*hk} and $\norm{J^n_k}_\infty \le \abs{\bar F^n_\ex}_1$ this gives
$$
\tilde \Sigma^n_{h,k} \subset B_{\ell^\infty}(\bar F^n_\ex(x^0_k), h \tilde\rho^n_{h,k}\norm{J^n_k}_\infty)
\subset B_{\ell^\infty}(\bar F^n_\ex(x^0_k), h \tilde\rho^n\abs{\bar F^n_\ex}_1).
$$
Thus we can apply \eqref{acc-Taylor} with $\omega = \tilde \Sigma^n_{h,k}$, and take the supremum over $k$: this yields
\begin{equation} \label{acc-Br-bS*}
\tilde e_{B,(r)}^{n}(h) \lesssim h^{r+1}(\tilde\rho^n\abs{\bar F^n_\ex}_1)^{r+1} \abs{\bar B^n_\ex}_{r+1} = h^{r+1} \tilde c^{n}_{F,(r)}
\end{equation}
which completes the proof.
\qed
\end{proof}

\begin{remark}[extension to heterogeneous bases]
\label{rem. heterogeneous}
As previously pointed out, since our analysis does not rely on a smoothing kernel argument it readily extends
to cases where the ``particles'' $\vp^0_{h,k}$ are not derived from a reference function $\vp$.
For instance one could use the same method to transport continuous piecewise polynomial basis functions
defined on an unstructured mesh of $\RR^d$, and the same results would apply under the usual shape regularity 
and quasi-uniformity assumptions.
\end{remark}

%%% %%% %%% %%% %%% %%% %%% %%% %%% %%% %%% %%% %%% %%% %%% %%% %%% %%% %%% %%% %%% %%% %%% %%% %%% %%% %%% %%% %%% %%% %%% %%% %%% 
%%% %%% %%% %%% %%% %%% %%% %%% %%% %%% %%% %%% %%% %%% %%% %%% %%% %%% %%% %%% %%% %%% %%% %%% %%% %%% %%% %%% %%% %%% %%% %%% %%% 

\section{Fully discrete schemes for the LTP and QTP methods}
\label{sec:discrete}

In this section we describe fully discrete implementations for the linear ($r=1$) and quadratic ($r=2$) transport operators 
$\bar T^n_{h,(r)}$. 
To comply with the algorithmic structure of standard particle codes, our schemes only involve pointwise evaluations of a given
forward flow $F^n$ which may either be the exact flow on the time step $F^n_\ex := F_{t^n,t^{n+1}}$, or more likely
some approximation to it.
In practice $F^n$ will typically be given by the algorithm used to push the particles in an existing code.

For each method we describe two implementations, the benefits of which may depend on the context. The first one is an 
``incremental'' approach where each particle carries a numerical approximation of the local derivatives of the backward flow 
$\bar B^n_\ex := (F_{0,t^n})^{-1}$. The method then consists of updating these derivatives at each time step with a few calls to 
the numerical flow $F^n$.

The second implementation follows a ``direct'' approach where the particles carry no flow derivatives but a small number of auxiliary 
markers $x^n_{k,\ell}$, $\ell \in I_{\FD, (r)}$. These markers play the role of local Finite-Difference stencils associated with
the particle trajectories $x^n_k$: they are initialized as cartesian patches
\begin{equation}\label{markers}
x^0_{k,\ell} := x^0_k + \ell h', \qquad \ell \in I_{\FD, (r)} \subset \ZZ^d
\end{equation}
and they are pushed forward in time, 
$$
x^{n+1}_{k,\ell} := F^n(x^n_{k,\ell}).
$$
Here $h'$ is an ad-hoc resolution that may be taken close to $h$ or smaller 
(for instance if the derivatives of $F^n$ are good approximations to those of $F^n_\ex$), without affecting the cost of the method. 
With this approach the local derivatives of the backward flow must be evaluated whenever needed from the positions of the markers,
which should only represent a mild overhead since no inter-particle communications are involved. 
Most importantly, the method readily provides local indicators for the backward flow error: 
indeed one may estimate the quantity \eqref{eB} by
\begin{equation} \label{hat-eB}
\hat e_{B,(r)}^{n}(h) 
:= \sup_{k\in\ZZ^d} 
%\hat e^n_{B,k}(h) := 
\Big\{\max_{\ell \in I_{\FD,(r)}} \norm{ \bar B^n_{(r),k}(x^n_{k,\ell}) - x^0_{k,\ell} }_\infty \Big\}.
\end{equation}

%%% %%% %%% %%% %%% %%% %%% %%% %%% %%%
%%% %%% %%% %%% %%% %%% %%% %%% %%% %%%
\subsection{Finite Difference implementation of the LTP method}
\label{sec:ltp}

The linearly-transformed particle (LTP) approximation is based on first-order expansions of the exact flow around the particle 
trajectories. In a time-discrete setting this consists of computing, following \eqref{Bhk-1}-\eqref{T1},
\begin{equation} \label{ltp}
\vp^n_{h,k}(x) = \bar T^n_{h,(1)} \vp^0_{h,k}(x) := \vp_h \big( D^n_k(x-x^n_k) \big)
\quad \text{ with } \quad 
\begin{cases} 
    x^n_k \approx \bar F^n_\ex(x^0_{k})    
    \\
    D^n_k \approx J_{\bar B^n_\ex}(x^n_{k}),
\end{cases}
\end{equation}
which corresponds to transporting the particles forward with %defining $\bar T^n_{h,(1)} \vp^0_{h,k}(x) = \vp^0_{h,k}(\bar B^{n,(1)}_{h,k})$ with
\begin{equation} \label{numflow-ltp}
\bar F^{n}_{(1),k} = (\bar B^{n}_{(1),k})^{-1}, \qquad
\bar B^{n}_{(1),k}(x) := x^0_k + D^n_k(x-x^n_k), \qquad k \in \ZZ^d.
\end{equation}
In the {\em direct} approach we push forward local stencils of the form \eqref{markers}, namely 
$$
%I_{\FD, (1)} := \{\pm e_j : j = 1, \ldots, d \}  %  ------    CENTERED
I_{\FD, (1)} := \{\alpha e_j : \alpha \in \{0, 1\}, ~ j = 1, \ldots, d \}
$$
and we approximate the Jacobian matrix of the forward flow with 
$$
%\bar J^n_k := \Big( \frac{(x^n_{k,e_j}-x^n_{k,0})_i}{h'} \Big)_{1\le i,j\le d} \approx J_{\bar F^n_\ex}(x^0_{k}).
\bar J^n_k := \left( \frac{(x^n_{k,e_j}-x^n_{k,0})_i}{h'} \right)_{1\le i,j\le d} \approx J_{\bar F^n_\ex}(x^0_{k}).
$$
Using the formal identity $J_{\bar F^n}(x^0_{k}) J_{\bar B^n}(x^n_{k}) = I_d$, we finally set 
$$
D^n_k := (\bar J^n_k)^{-1}.
$$
In the {\em incremental} approach the matrix $D^n_k$ is stored in place of the phase space markers. Since $\bar B^0_\ex = I$
it is initialized with $D^0_k := I_d$ and updated using the formal identity 
$J_{\bar B^{n+1}}(x^{n+1}_{k}) = J_{\bar B^n}(x^n_{k}) J_{B^n}(x^{n+1}_{k})$.
Using a local stencil we compute 
$$
J^n_k := \left( \frac{F^n_i(x^n_{k}+h'e_j) - F^n_i(x^n_{k}-h'e_j)}{2h'} \right)_{1\le i,j\le d} \approx J_{F^n_\ex}(x^n_{k})
$$
and then set 
$$
D^{n+1}_k := D^n_k (J^n_k)^{-1}.
$$

%%% %%% %%% %%% %%% %%% %%% %%% %%% %%%
%%% %%% %%% %%% %%% %%% %%% %%% %%% %%%
\subsection{Finite Difference implementation of the QTP method}
\label{sec:qtp}

The quadratically-transformed particle (QTP) approximation is based on second-order expansions of the exact flow around the particle 
trajectories. In a time-discrete setting this consists of computing, following \eqref{Tr} and \eqref{Bhk-r},
$$ %\begin{equation} \label{qtp}
\vp^n_{h,k}(x) = \bar T^n_{h,(2)} \vp^0_{h,k}(x) = \chi_{\Sigma^n_{h,k}}(x) 
        \vp_h \big( D^n_k(x-x^n_k) + \tfrac 12 (x-x^n_k)^t Q^n_k (x-x^n_k) \big),
$$ %\end{equation}
which corresponds to taking as numerical backward flows the quadratic mappings
\begin{equation} \label{numflow-qtp}
\bar B^{n}_{(2),k}(x) := x^0_k + D^n_k(x-x^n_k) + \tfrac 12 (x-x^n_k)^t Q^n_k (x-x^n_k),
\qquad k \in \ZZ^d
\end{equation}
and we recall that in general these quadratic mappings are not invertible:
for this reason we restrict the particles to a-priori supports $\Sigma^n_{h,k}$ that will be specified below.
As for the particle centers $x^n_k$ and deformation matrices $D^n_k$ they are defined and computed as in Section~\ref{sec:ltp} above.
Finally the quadratic deformation terms 
$$
(x-x^n_k)^t Q^n_k (x-x^n_k) = \Big( \sum_{j_1,j_2 = 1}^d (Q^n_k)_i (x-x^n_k)_{j_1} (x-x^n_k)_{j_2}\Big)_{1 \le i \le d}
$$ 
involve approximations of the Hessian matrices of the backward flow,
$$
(Q^n_k)_i \approx H_{(\bar B^n_\ex)_i}(x^n_k) := \Big( \partial_{j_1,j_2} {(\bar B^n_\ex)_i}(x^n_k)\Big)_{1 \le j_1,j_2 \le d}
$$
that are computed as follows. In the {\em direct} approach we push forward local stencils of the form \eqref{markers} with 
$$
I_{\FD, (2)} := \{\alpha_1 e_{j_1} + \alpha_2 e_{j_2} : \alpha_1, \alpha_2 \in \{0, 1\}, ~ j_1, j_2 = 1, \ldots, d \}
$$
and first compute approximate forward Hessian matrices with finite differences,
$$
(\bar H^n_k)_i := \Big( (h')^{-2} \sum_{\alpha_1,\alpha_2 = 0}^1 
                            (-1)^{\alpha_1+\alpha_2} \big(x^n_{k,\alpha_1 e_{j_1}+\alpha_2 e_{j_2}}\big)_i
                    \Big)_{1\le j_1,j_2\le d} 
               \approx H_{(\bar F^n_\ex)_i}(x^0_{k}).
$$
Differentiating twice (at $x^0_k$) the formal identity $0 = \bar B^n \bar F^n$ we then obtain
\begin{equation} \label{d2_bBn_bFn}
0 = \big(J_{\bar F^n}(x^0_{k})\big)^t H_{(\bar B^n)_i}(x^n_k) J_{\bar F^n}(x^0_{k}) 
    + \sum_{j=1}^d \big(J_{\bar B^n}(x^n_{k})\big)_{i,j} H_{(\bar F^n)_j}(x^0_k)
\end{equation}
so that we finally set 
\begin{equation}\label{Q-direct}
(Q^n_k)_i := - (D^n_k)^t \Big( \sum_{j=1}^d (D^n_{k})_{i,j} (\bar H^n_k)_j \Big) D^n_k.
\end{equation}

In the {\em incremental} approach the $d \times d$ matrices $(Q^n_k)_i$ are stored in place of the markers. Since $\bar B^0_\ex = I$
they are initialized with $(Q^0_k)_i := 0$. To update them we differentiate twice (at $x^{n+1}_k$) the formal identity 
$\bar B^{n+1} = \bar B^{n} B^n$: this gives
\begin{multline}\label{d2_bBn_Bn}
H_{(\bar B^{n+1})_i}(x^{n+1}_k) = \big(J_{B^n}(x^{n+1}_{k})\big)^t H_{(\bar B^n)_i}(x^n_k) J_{B^n}(x^{n+1}_{k})
            \\
            + \sum_{j=1}^d \big(J_{\bar B^n}(x^n_{k})\big)_{i,j} H_{(B^n)_j}(x^{n+1}_k)
\end{multline}
and we observe that rewriting \eqref{d2_bBn_bFn} on the local time step yields a formula expressing the local backward 
Hessian matrices in terms of the forward ones.
Thus we first compute finite difference approximations of the matrices $H_{(F^n_\ex)_i}(x^n_{k})$,
$$
(H^n_k)_i := \Big( (h')^{-2} \sum_{\alpha_1,\alpha_2 = 0}^1 
                                    (-1)^{\alpha_1+\alpha_2} (F^n)_i(x^n_{k} + h'(\alpha_1 e_{j_1}+\alpha_2 e_{j_2}))
                              \Big)_{1\le j_1,j_2\le d} 
$$
%\begin{align*}
%(H^n_k)_i &:= \Big( (h')^{-2} \sum_{\alpha_1,\alpha_2 = 0}^1 
%                                    (-1)^{\alpha_1+\alpha_2} (F^n)_i(x^n_{k} + h'(\alpha_1 e_{j_1}+\alpha_2 e_{j_2}))
%                              \Big)_{1\le j_1,j_2\le d} 
%        \\
%        &\approx H_{(F^n_\ex)_i}(x^n_{k})
%\end{align*}
then we approximate the local backward Hessian matrices similarly as in \eqref{Q-direct},
$$
(\check H^n_k)_i := - (\check J^n_k)^{t} \Big( \sum_{j=1}^d (\check J^n_{k})_{i,j} (H^n_k)_j \Big) \check J^n_k \approx H_{(B^n_\ex)_i}(x^{n+1}_{k})
$$
(here $\check J^n_k = (J^n_k)^{-1}$ denotes the approximate local backward Jacobian), and finally using \eqref{d2_bBn_Bn} we let 
\begin{equation}\label{Q-increment}
\begin{split}
(Q^{n+1}_k)_i & := (\check J^n_k)^{t} (Q^n_k)_i \check J^n_k + \sum_{j=1}^d (D^n_k)_{i,j} (\check H^n_k)_j
    \\
    & = (\check J^n_k)^{t} \Big( (Q^n_k)_i -  \sum_{j,j'=1}^d (D^n_k)_{i,j} (\check J^n_k)_{j,j'} (H^n_k)_{j'} \Big) \check J^n_k.
\end{split}
\end{equation}

It remains to define the a-priori particle support. From Theorem~\ref{th:conv-f-Tr} we know that setting 
$$
\Sigma^n_{h,k} := \bar F^n_{(1),k}(B_{\ell^\infty}(x^0_k,h \tilde \rho^n_{h,k}))
~~  \text{ with } ~~
\tilde \rho^n_{h,k} := \rho^0 + \frac 1h \norm{\bar B^n_{(1),k}-\bar B^n_\ex}_{L^\infty(\bar F^n_\ex(\Sigma^0_{h,k}))}
%\bar F^{n,(1)}_{h,k}(B_{\ell^\infty}(x^0_k,h \rho^0 \norm{D^n_k}_\infty \abs{(D^n_k)^{-1}}_\infty))
$$
is enough to guarantee the second order convergence of the method, and in the direct approach we can use
the local indicators in \eqref{hat-eB} to estimate the terms $\tilde \rho^n_{h,k}$. In practice however we 
have observed that these supports were too large, yielding strong oscillations in the solutions. 
Now, since the latter are likely caused by the non-invertibility of the quadratic flow $\bar B^{n}_{(2),k}$, 
we decided to further restrict the particle supports to the regions where this flow is locally
invertible. Specifically, we define the a-priori particle supports as
$$
\Sigma^n_{h,k} := \{ x \in \bar F^n_{(1),k}(B_{\ell^\infty}(x^0_k,h \tilde \rho^n_{h,k})) :
    \det(J_{\bar B^{n}_{(2),k}}(x)) > 0 \}.
$$
Since $J_{\bar B^{n}_{(2),k}}(x) = D^n_k + \big(\sum_{j'} (Q^n_k)_{i,j,j'}(x-x^n_k)_{j'}\big)_{1 \le i,j \le d}$ 
this strategy is easy to implement. In the next section we will see that it results in a very robust numerical method.

%%% %%% %%% %%% %%% %%% %%% %%% %%% %%% %%% %%% %%% %%% %%% %%% %%% %%% %%% %%% %%% %%% %%% %%% %%% %%% %%% %%% %%% %%% %%% %%% %%% 
%%% %%% %%% %%% %%% %%% %%% %%% %%% %%% %%% %%% %%% %%% %%% %%% %%% %%% %%% %%% %%% %%% %%% %%% %%% %%% %%% %%% %%% %%% %%% %%% %%% 

\section{Numerical experiments}
\label{sec:num}

In this section we test the proposed LTP and QTP schemes and compare them with the standard TSP and FSL methods described in 
Section~\ref{sec:overview}, using either $M'_4$ or cubic B-spline particles as described in Section~\ref{sec:shapes}. 
To assess the robustness of the method with respect to the velocity field and the initial data,
we use several 2d test-cases summarized in Table~\ref{tab:test-cases}. For the velocity fields we consider
\begin{itemize} 
\item[$\bullet$]
    the reversible ``swirling'' velocity field proposed by LeVeque \cite{LeVeque.1996.sinum},
    $$ %\begin{equation} \label{u_swirl} % was {LV_field}
    u_{\rm SW}(t,x; T)
    %           :=   \cos\big(\frac {\pi t}{T}\big)
    %                \begin{pmatrix}
    %                        {-}\sin^2(\pi x_1) \sin(2\pi x_2)
    %                \\
    %                \phantom{-}\sin^2(\pi x_2) \sin(2\pi x_1)
    %                \end{pmatrix}        
           := \frac 1 \pi \cos\Big(\frac {\pi t}{T}\Big) \curl \phi_{\rm SW}(x)
    %    ~~~ \text{with} ~~~ \phi_{\rm SW}(x) := - \sin^2(\pi x_1)\sin^2(\pi x_2)
    $$ %\end{equation}    
    with $\phi_{\rm SW}(x) := - \sin^2(\pi x_1)\sin^2(\pi x_2)$~;
\item[$\bullet$]
    another reversible velocity field emulating a Raylegh-Benard convection cell,
    $$ 
    %\begin{equation} \label{u_rb}
    u_{\rm RB}(t,x; T)
           := \cos\Big(\frac {\pi t}{T}\Big) \curl \phi_{\rm RB} (x)
%           \quad \text{with} \quad  
%           \phi_{\rm RB}(x) := \big(x_1-\tfrac 12\big)(x_1-x_1^2)(x_2-x_2^2) 
    $$ %\end{equation}
    with $\phi_{\rm RB}(x) := \big(x_1-\tfrac 12\big)(x_1-x_1^2)(x_2-x_2^2)$~;
\item[$\bullet$]
    and finally a constant non-linear rotation field derived from Example~2 in \cite{Bokanowski.Garcke.Griebel.Klompmaker.2013.jsc}, 
    $$ %\begin{equation} \label{u_nlr} 
    u_{\rm NLR}(x)
           :=   
            \alpha(x)
            \begin{pmatrix}
                      \tfrac 12 - x_2
            \\ \noalign{\smallskip}
            x_1 - \tfrac 12
            \end{pmatrix}
    \quad \text{with} \quad \alpha(x) := \Big(1-\frac{\norm{x-(\tfrac 12, \tfrac 12)}_2}{0.4}\Big)^3_+ .
    $$ %\end{equation}
%    with $\alpha(x) = (1-2.5\norm{x-(\tfrac 12, \tfrac 12)}_2)^3_+$.
\end{itemize}
Here the time-dependency in $u_{\rm SW}$ and $u_{\rm RB}$ yield reversible problems: at $t = T/2$ the solutions 
reach a state of maximum stretching, and they revert to their initial value at $t=T$.
As for the non-linear rotation field $u_{\rm NLR}$, it is associated with the exact backward flow
$$
\bar B^n_\ex(x) =  
                    \begin{pmatrix} 
                    \tfrac 12
                    \\ \noalign{\smallskip}
                    \tfrac 12
                    \end{pmatrix}
                    +
                    \begin{pmatrix} 
                    \cos(\alpha(x)t^n) &  \sin(\alpha(x)t^n)
                    \\ \noalign{\smallskip}
                    - \sin(\alpha(x)t^n) & \cos(\alpha(x)t^n)
                    \end{pmatrix}
                    \begin{pmatrix} 
                    x_1-\tfrac 12
                    \\ \noalign{\smallskip}
                    x_2-\tfrac 12
                    \end{pmatrix},
$$
and the exact solutions are given by $f(t^n,x) = f^0(\bar B^n_\ex(x))$.
In addition to the above velocity fields we consider the following initial data:
\begin{itemize} 
\item[$\bullet$]
    smooth humps of approximate radius $0.2$ given by
    $$ %\begin{equation} \label{f_hump}
    f^0_{\rm hump}(x;\bar x)
           :=  \frac 12 \Big( 1 + \erf\Big( \tfrac 13 (11- 100\norm{x-\bar x}_{2}) \Big)\Big)              
    $$ %\end{equation}
    and centered on $\bar x = (0.5,0.4)$ or $(0.5,0.7)$, depending on the cases~;
\item[$\bullet$]
    a cone of radius $0.15$, 
    $$ %\begin{equation}   \label{f_cone}
    f^0_{\rm cone}(x;\bar x)
            := \big( 1 - \tfrac{20}{3}\norm{x-\bar x}_{2}\big)_+
    $$ %\end{equation}
    centered on $\bar x = (0.5,0.25)$~;
\item[$\bullet$]
    and finally for the non-linear rotation field $u_{\rm NLR}$ we take an initial data corresponding to Example~2 from
    \cite{Bokanowski.Garcke.Griebel.Klompmaker.2013.jsc}, i.e.,
    $$ %\begin{equation}  \label{f_nlr} 
    f^0(x)
           :=  x_2 - \tfrac 12.
    $$ %\end{equation}
\end{itemize}    
By combining the above values we obtain the four test-cases defined in Table~\ref{tab:test-cases},
and accurate solutions are shown in Figures~\ref{fig:cc-SW-cone}-\ref{fig:cc-NLR} for the purpose of illustration.
In Table~\ref{tab:test-cases} we also give the respective time steps $\Dt$ used in the time integration of the particle trajectories. 
In every case indeed, the numerical flow $F^n$ is computed with a RK4 scheme, and the time steps have been 
taken small enough to have no significant effect on the final accuracy. It happens that in every case we
have $\Dt = T/100$, but this is unintended.

    %%%%%%%%%%%%%%%%%%%%%%%%%%%%%%%%%%%%%%%%%%%%%%%%%%%%%%%%%%%%%%%%%%%%%%%%%%%%%%%%%%%%%%%%%%%%%%%%%%%%%%%%%%%%%%%%%%%%%%%%
    % details of the smooth hump used in the SW-hump (previously called RL-2) test case:
    %   
    %        erf_blob_r_min = set_param('erf_blob_r_min',0.02,params)
    %        erf_blob_r_max = set_param('erf_blob_r_max',0.2,params)
    %        # blob center 
    %        erf_blob_x0 = set_param('erf_blob_x0',0.5,params)
    %        erf_blob_x1 = set_param('erf_blob_x1',0.7,params)
    %
    %        c1 = 6./(get_param('erf_blob_r_max',params)-get_param('erf_blob_r_min',params))
    %        af90.erf_blob_c1 = c1
    %        af90.erf_blob_c2 = 3./c1 + get_param('erf_blob_r_min',params)
    %
    %        f0 = 0.5*erf(erf_blob_c1*(erf_blob_c2-sqrt((x0-erf_blob_x0)**2+(x1-erf_blob_x1)**2))) + 0.5
    %
    % python commands to get the constants: 
    %
    % r_max = 0.2
    % r_min = 0.02
    % c1 = 6./(r_max-r_min)
    % c2 = 3./c1 + r_min
    % print "b = c1 = ", c1
    % print "a = c1*c2 = ", c1*c2
    %
    %%%%%%%%%%%%%%%%%%%%%%%%%%%%%%%%%%%%%%%%%%%%%%%%%%%%%%%%%%%%%%%%%%%%%%%%%%%%%%%%%%%%%%%%%%%%%%%%%%%%%%%%%%%%%%%%%%%%%%%%

% For tables use
\begin{table}[!h]
%\begin{center}
% table caption is above the table
\caption{Definition of the benchmark test-cases}
\label{tab:test-cases}       % Give a unique label
% For LaTeX tables use
\begin{tabular}{lllll}
\hline\noalign{\smallskip}
name & $u(t,x)$ & $f^0(x)$ & $T$ & $\Dt$ 
\\
\noalign{\smallskip}\hline\noalign{\smallskip}
SW-cone
    & $u_{\rm SW}(t,x; T)$
        &   $f^0_{\rm cone}(x; \bar x)$ with $\bar x = (0.5,0.25)$
             & 5
                & 0.05 
\\
\noalign{\smallskip}\hline\noalign{\smallskip}
SW-hump 
    &   $u_{\rm SW}(t,x; T)$
        &   $f^0_{\rm hump}(x; \bar x)$ with $\bar x = (0.5,0.7)$
            & 5
                & 0.05 
\\
\noalign{\smallskip}\hline\noalign{\smallskip}
RB-hump
    & $u_{\rm RB}(t,x; T)$
        &   $f^0_{\rm hump}(x; \bar x)$ with $\bar x = (0.5,0.4)$
            & 3 
                & 0.03 
\\
\noalign{\smallskip}\hline\noalign{\smallskip}
NLR 
    & $u_{\rm NLR}(x)$ 
        &  $x_2 - \tfrac 12$ 
            & 50
                & 0.5
\\
\noalign{\smallskip}\hline
\end{tabular}
%\end{center}
\end{table}

%%% %%% %%% %%% %%% %%% %%% %%% %%% %%%
%%% %%% %%% %%% %%% %%% %%% %%% %%% %%%
\subsection{Numerical convergence rates}

\begin{figure} [!ht]
\begin{center}
\begin{tabular}{cc}
\hspace{30pt}
\includegraphics[height=0.3\textwidth]{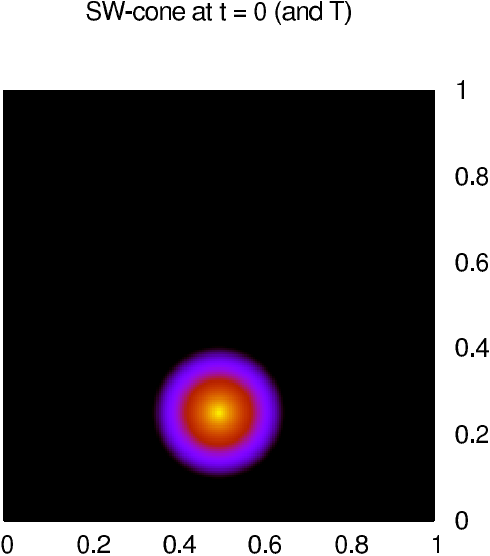}
& % \hspace{5pt} 
\hspace{20pt}
\includegraphics[height=0.3\textwidth]{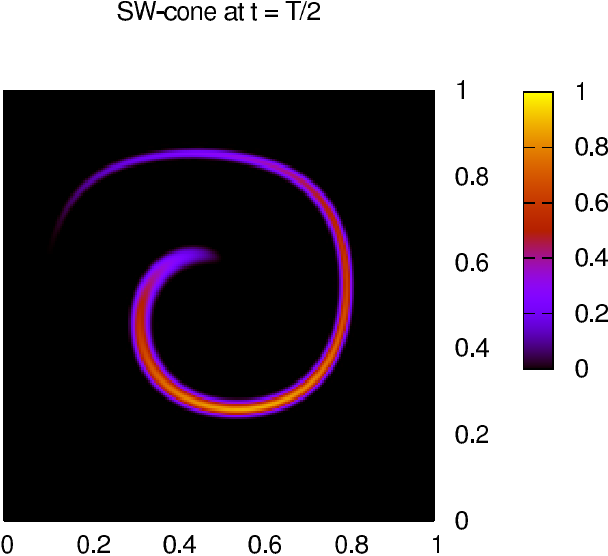}
\vspace{5pt}
\\
\includegraphics[width=0.45\textwidth]{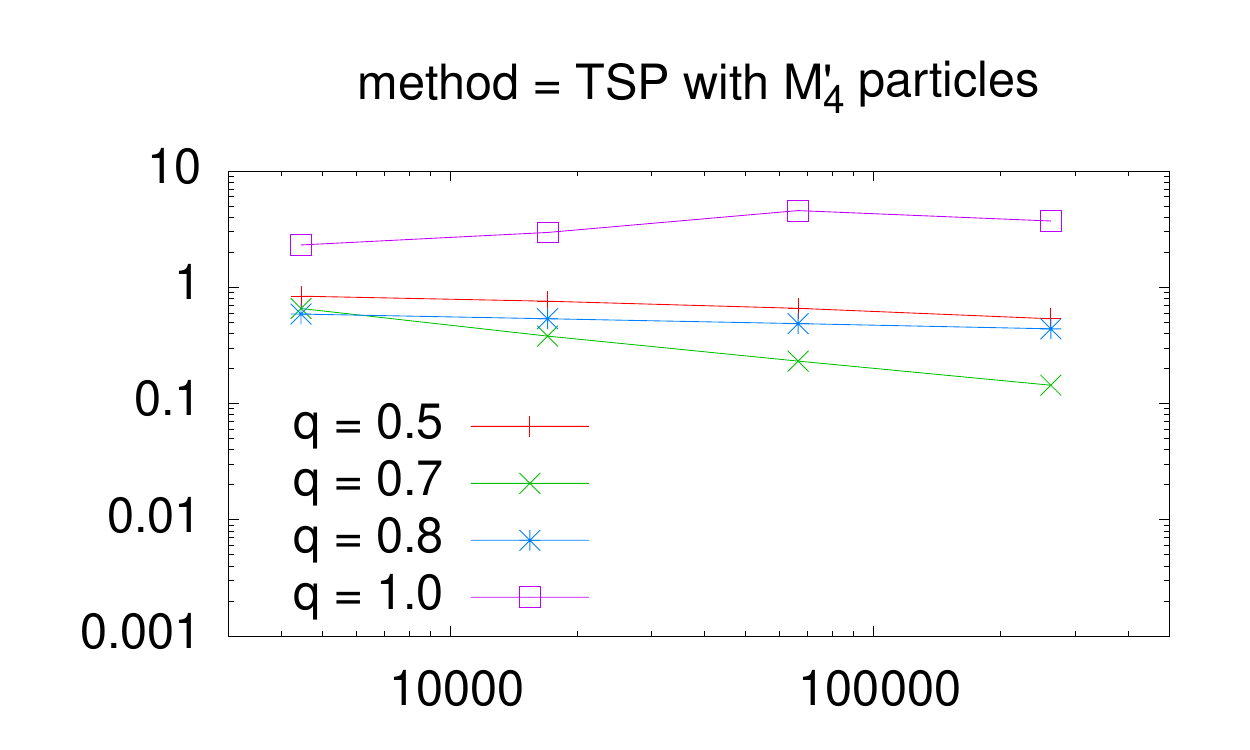}
& % \hspace{5pt} 
\includegraphics[width=0.45\textwidth]{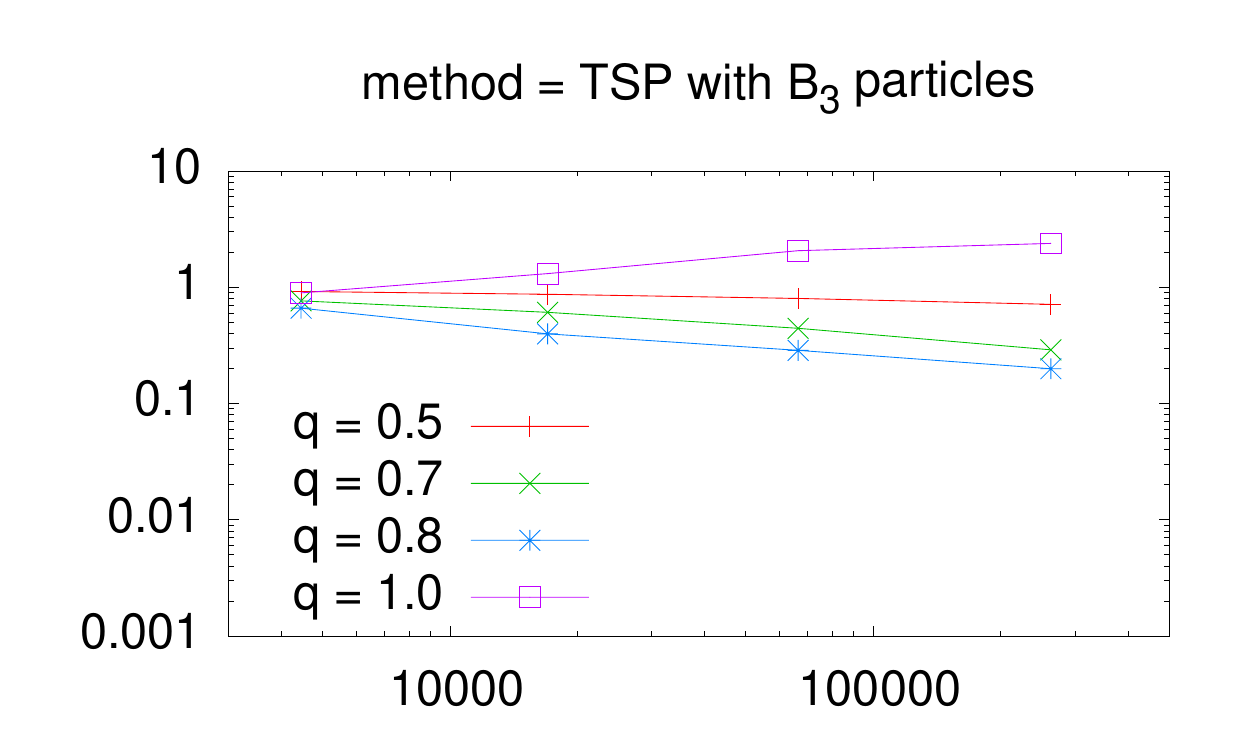}
\vspace{5pt}
\\
\includegraphics[width=0.45\textwidth]{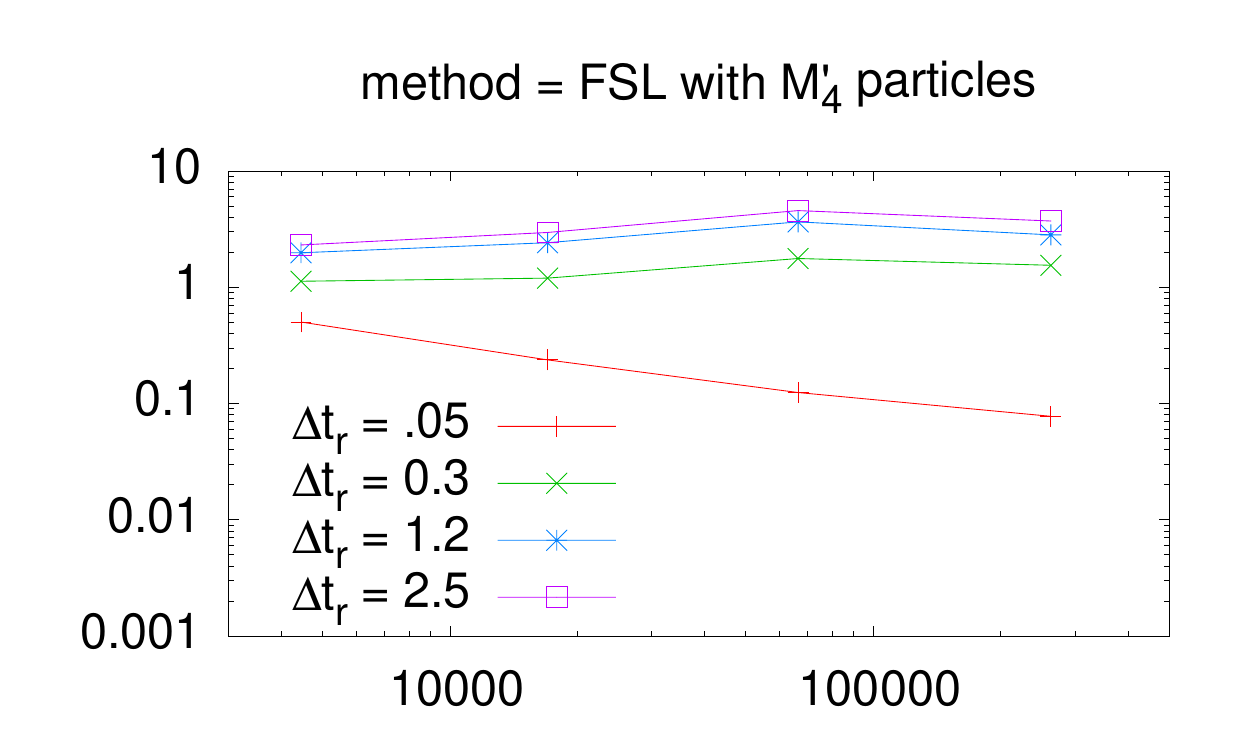}
& % \hspace{5pt} 
\includegraphics[width=0.45\textwidth]{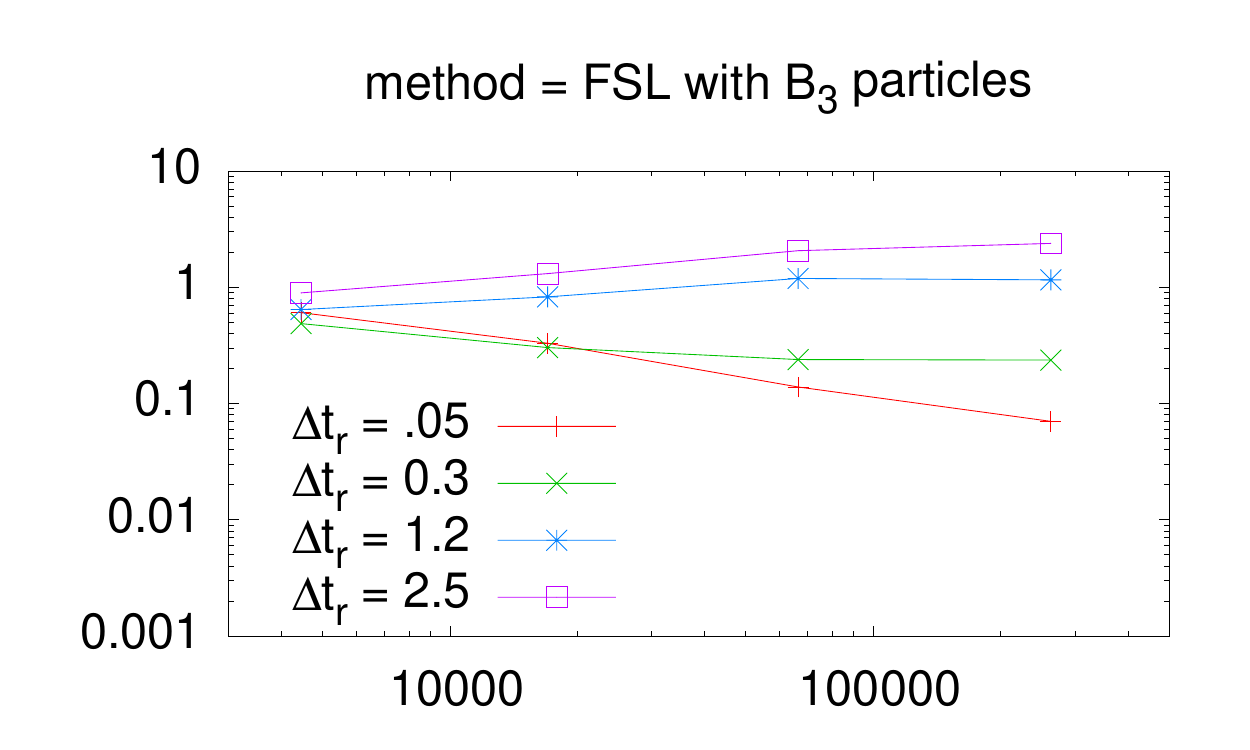}
\vspace{5pt}
\\
\includegraphics[width=0.45\textwidth]{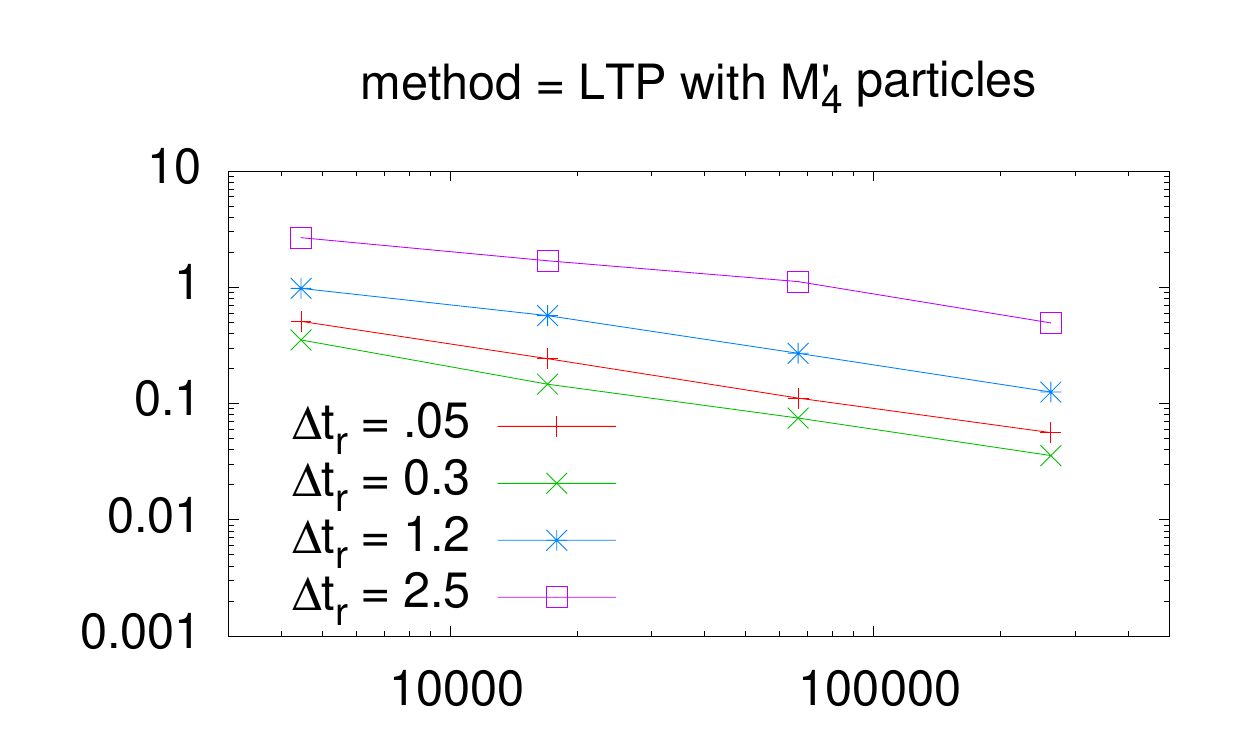}
& % \hspace{5pt} 
\includegraphics[width=0.45\textwidth]{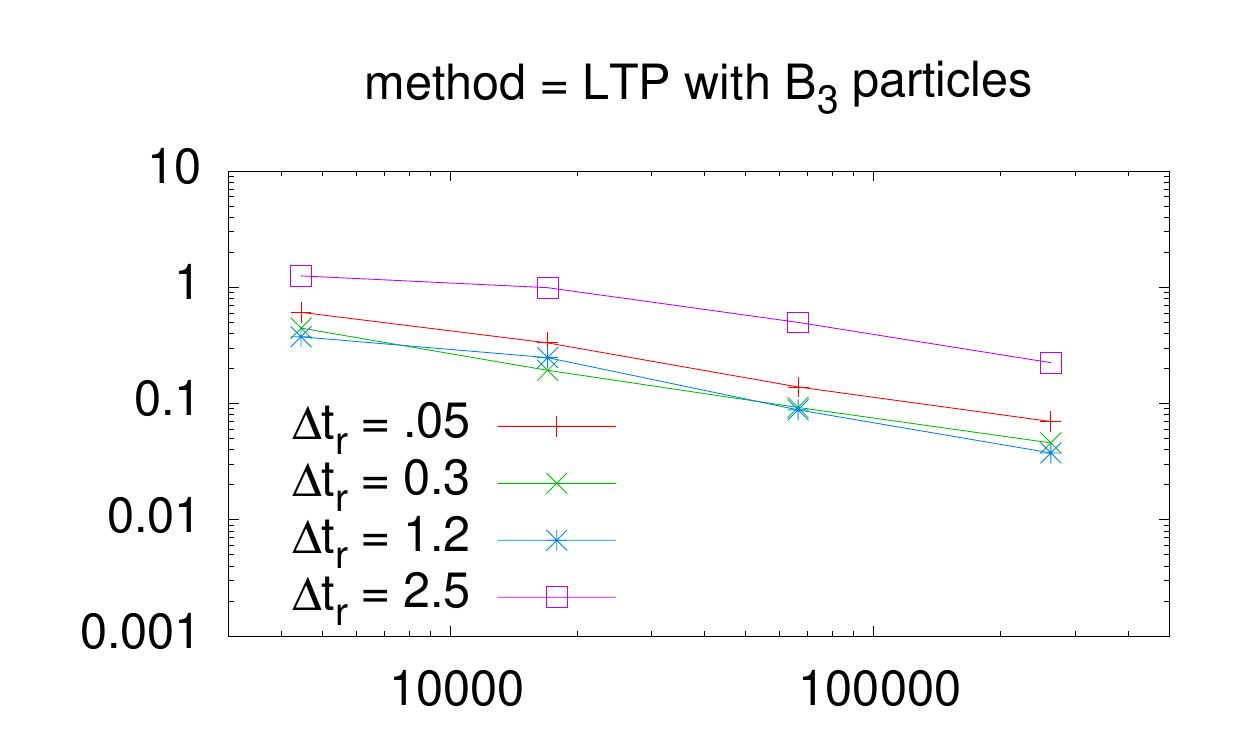}
\vspace{5pt}
\\
\includegraphics[width=0.45\textwidth]{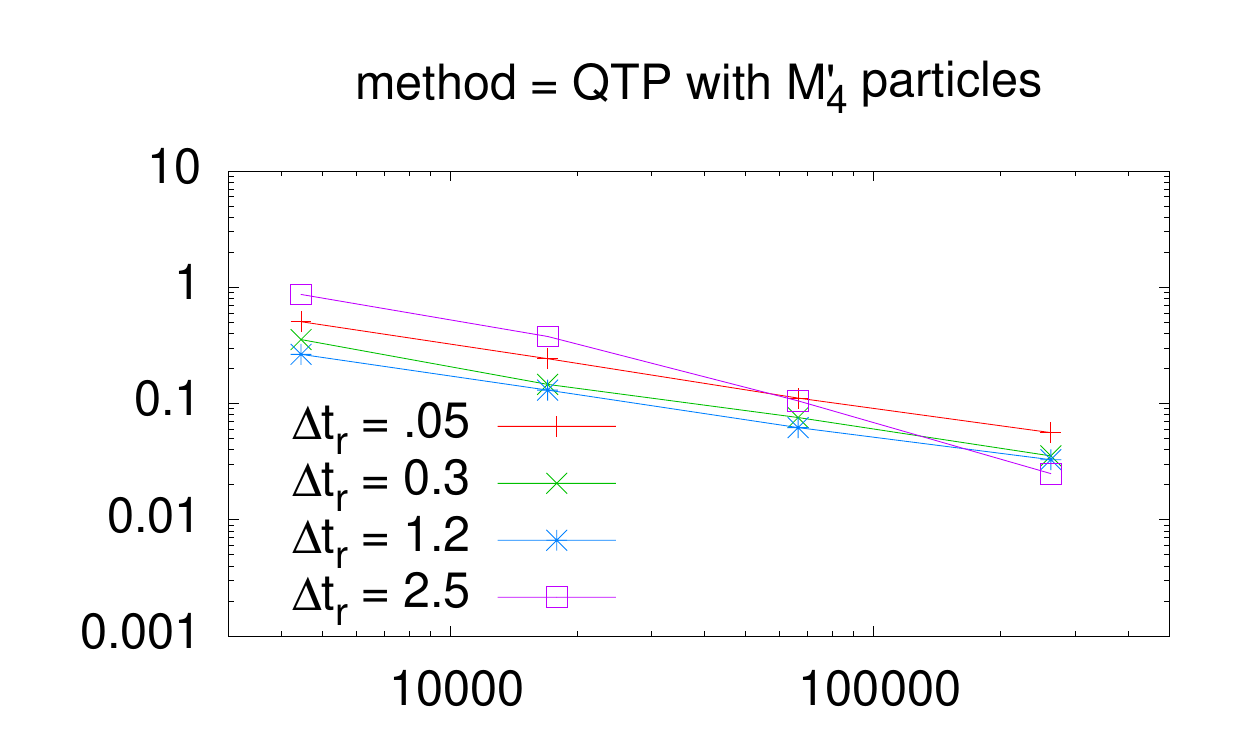}
& % \hspace{5pt} 
\includegraphics[width=0.45\textwidth]{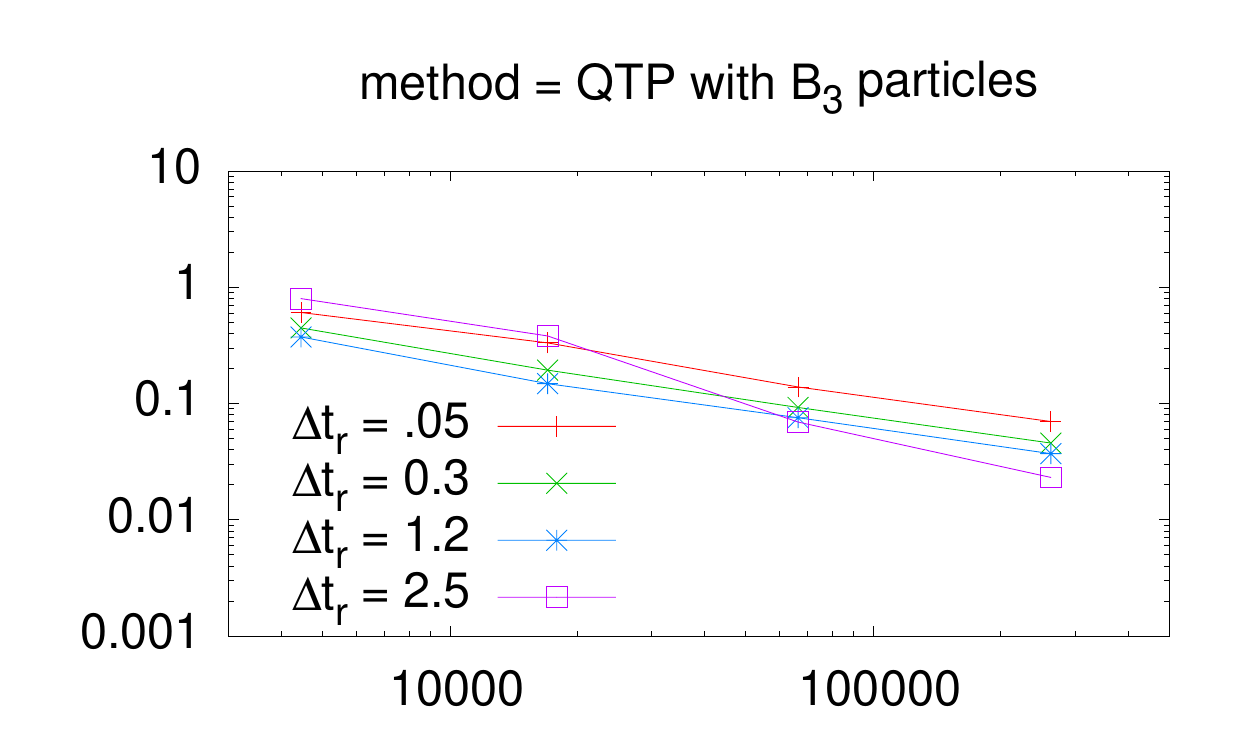}
\end{tabular}
  \caption{(Color) Convergence curves (relative $L^\infty$ errors at $t=T$ vs. average number of active particles)
  for the reversible test case SW-cone defined in Table~\ref{tab:test-cases}, solved with the different methods (see text for details). 
  The first row shows the profile of the exact solution:
  the initial (and final) density $f^0=f(T)$ is on the left, whereas the intermediate solution $f(T/2)$ 
  (with maximum stretching) is on the right.
  }
  \label{fig:cc-SW-cone}
 \end{center}
\end{figure}

\begin{figure} [!ht]
\begin{center}
\begin{tabular}{cc}
\hspace{30pt}
\includegraphics[height=0.3\textwidth]{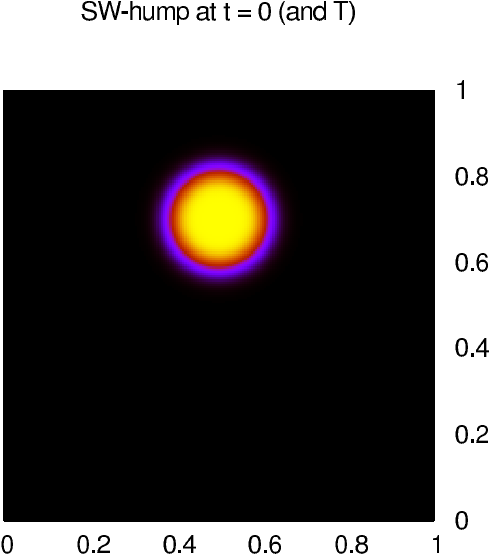}
& % \hspace{5pt} 
\hspace{20pt}
\includegraphics[height=0.3\textwidth]{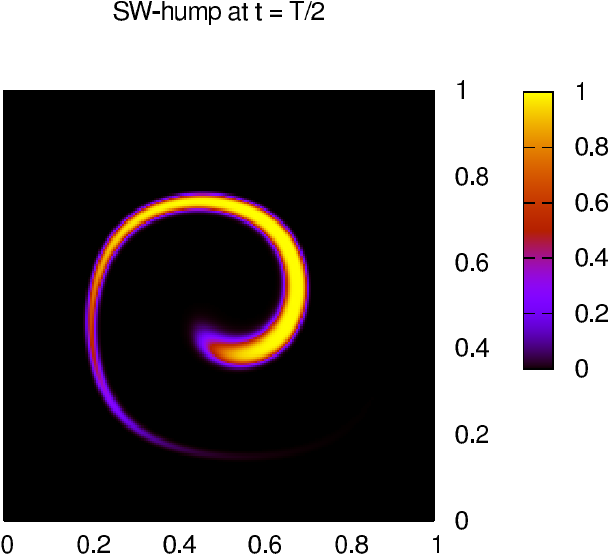}
\vspace{5pt}
\\
\includegraphics[width=0.45\textwidth]{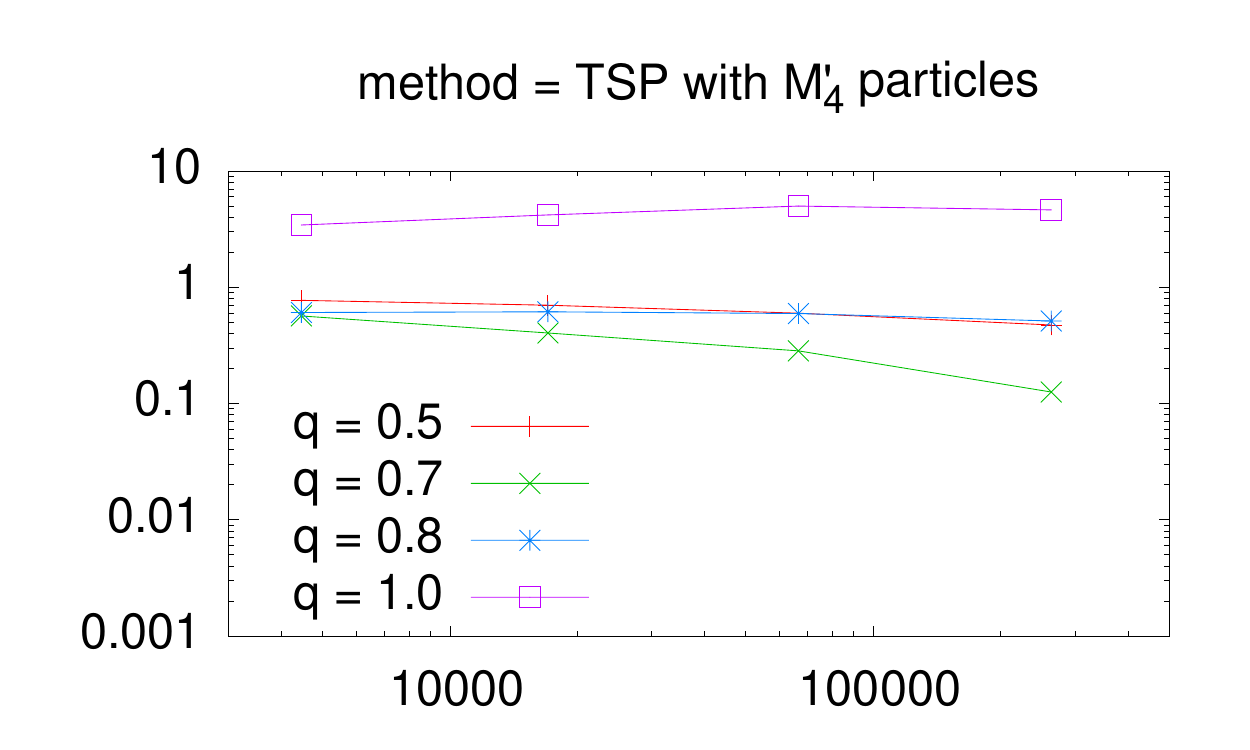}
& % \hspace{5pt} 
\includegraphics[width=0.45\textwidth]{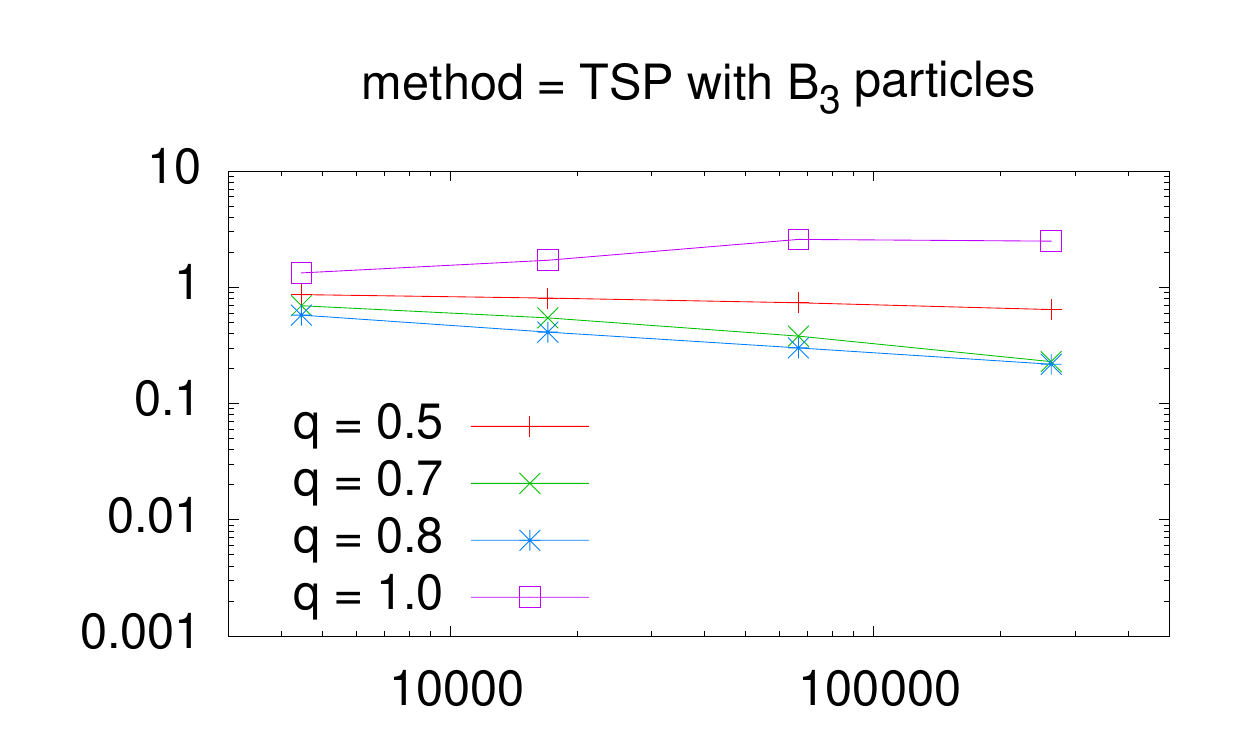}
\vspace{5pt}
\\
\includegraphics[width=0.45\textwidth]{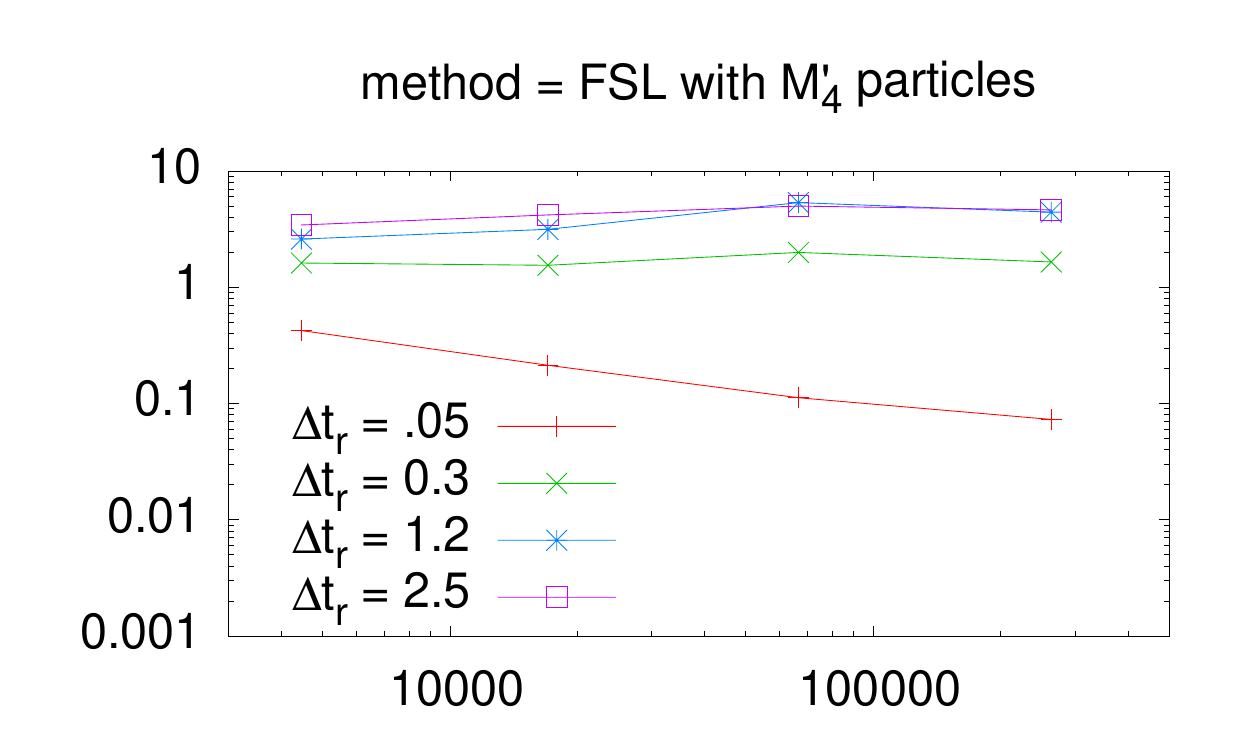}
& % \hspace{5pt} 
\includegraphics[width=0.45\textwidth]{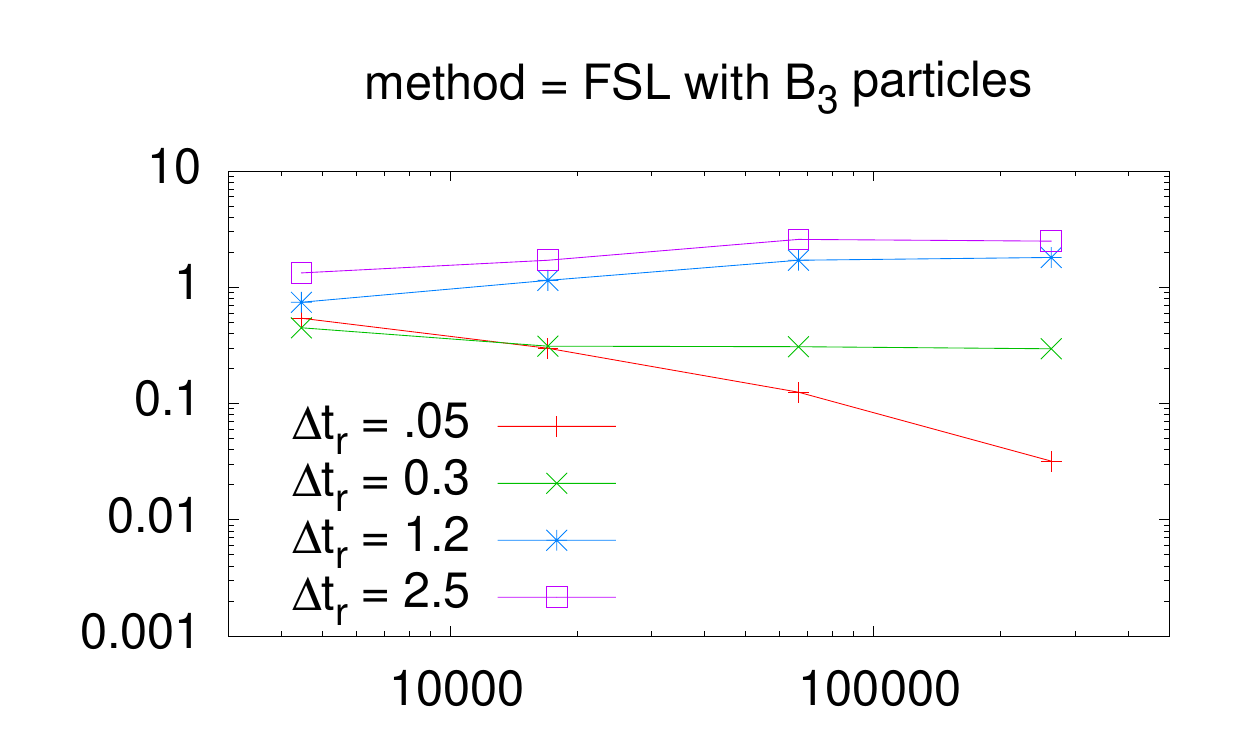}
\vspace{5pt}
\\
\includegraphics[width=0.45\textwidth]{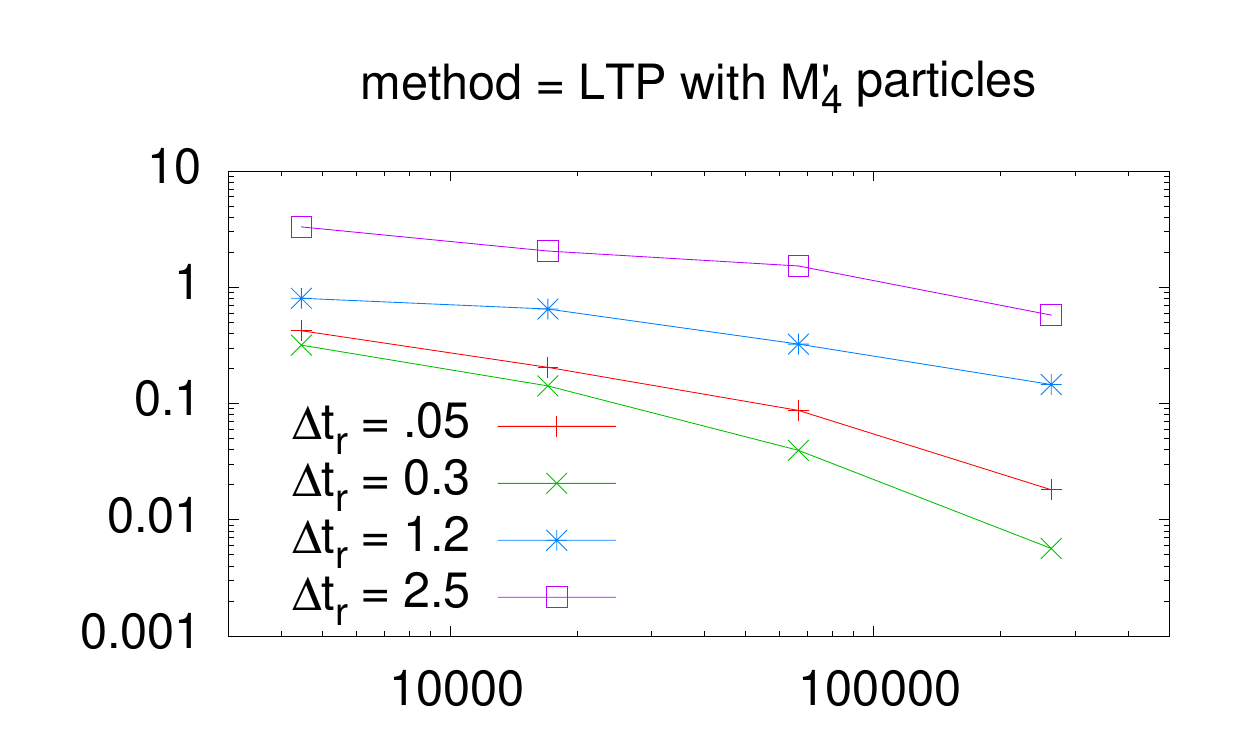}
& % \hspace{5pt} 
\includegraphics[width=0.45\textwidth]{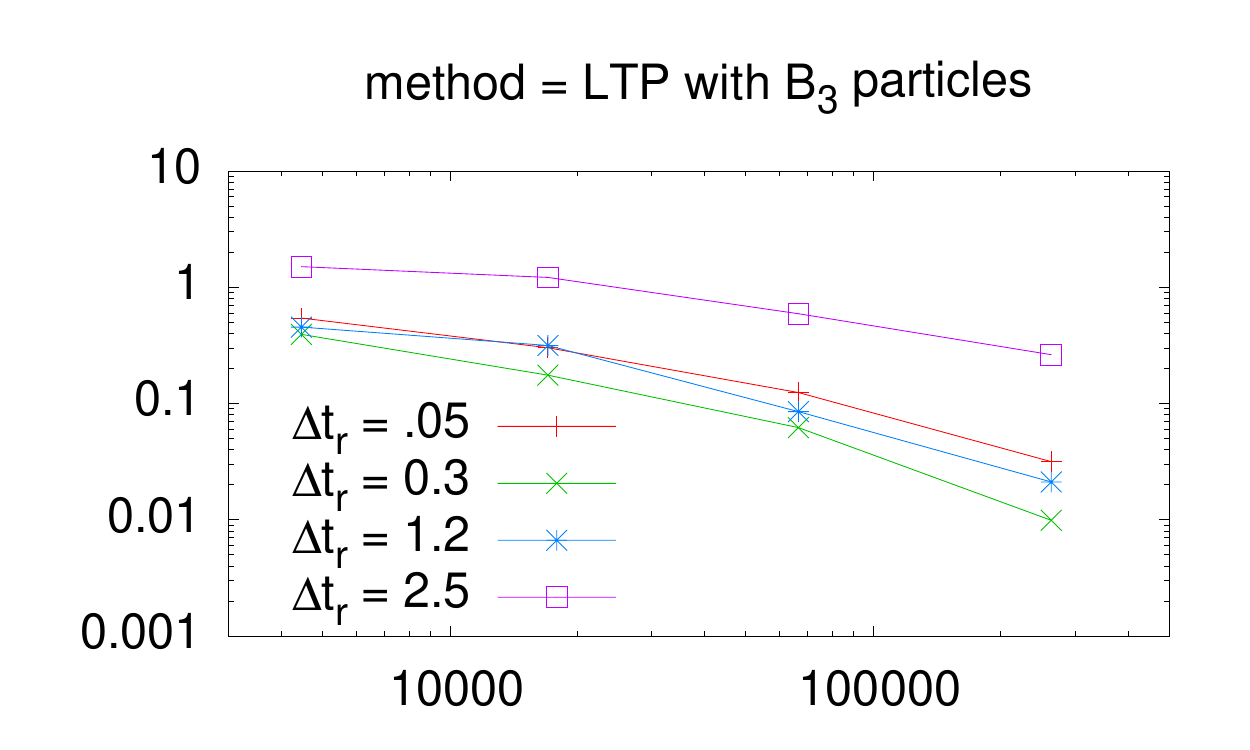}
\vspace{5pt}
\\
\includegraphics[width=0.45\textwidth]{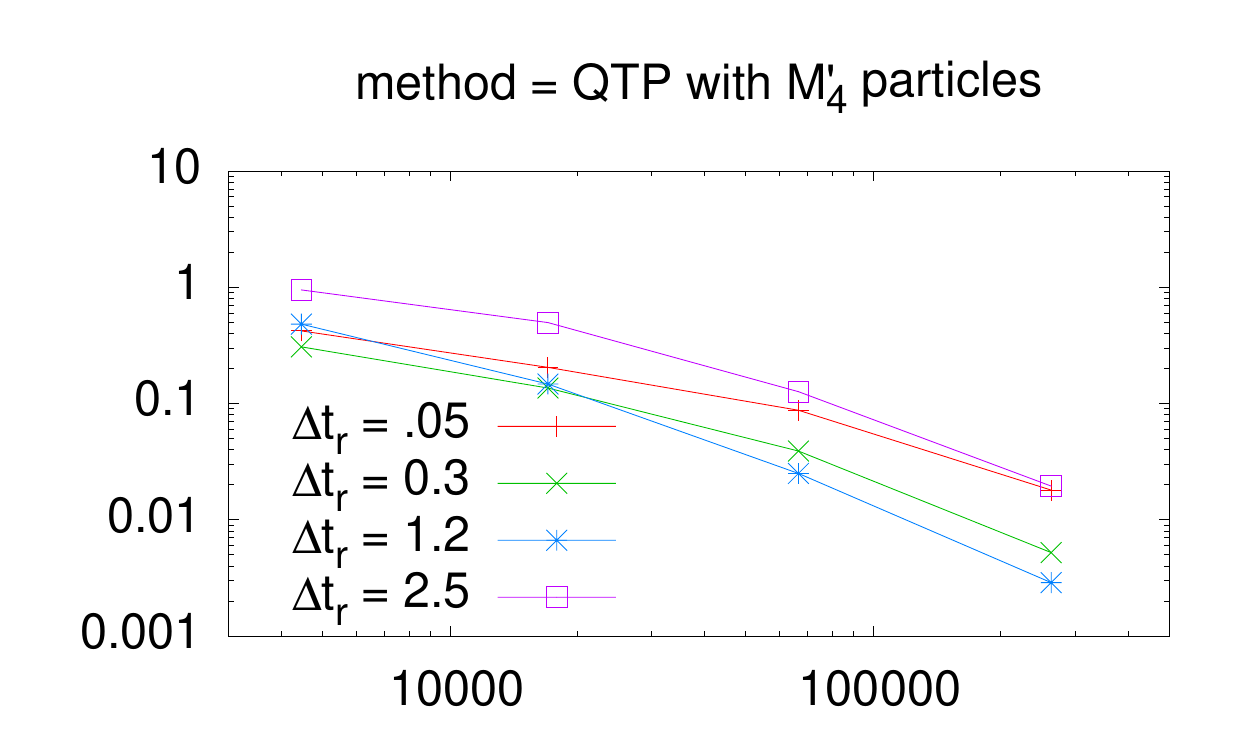}
& % \hspace{5pt} 
\includegraphics[width=0.45\textwidth]{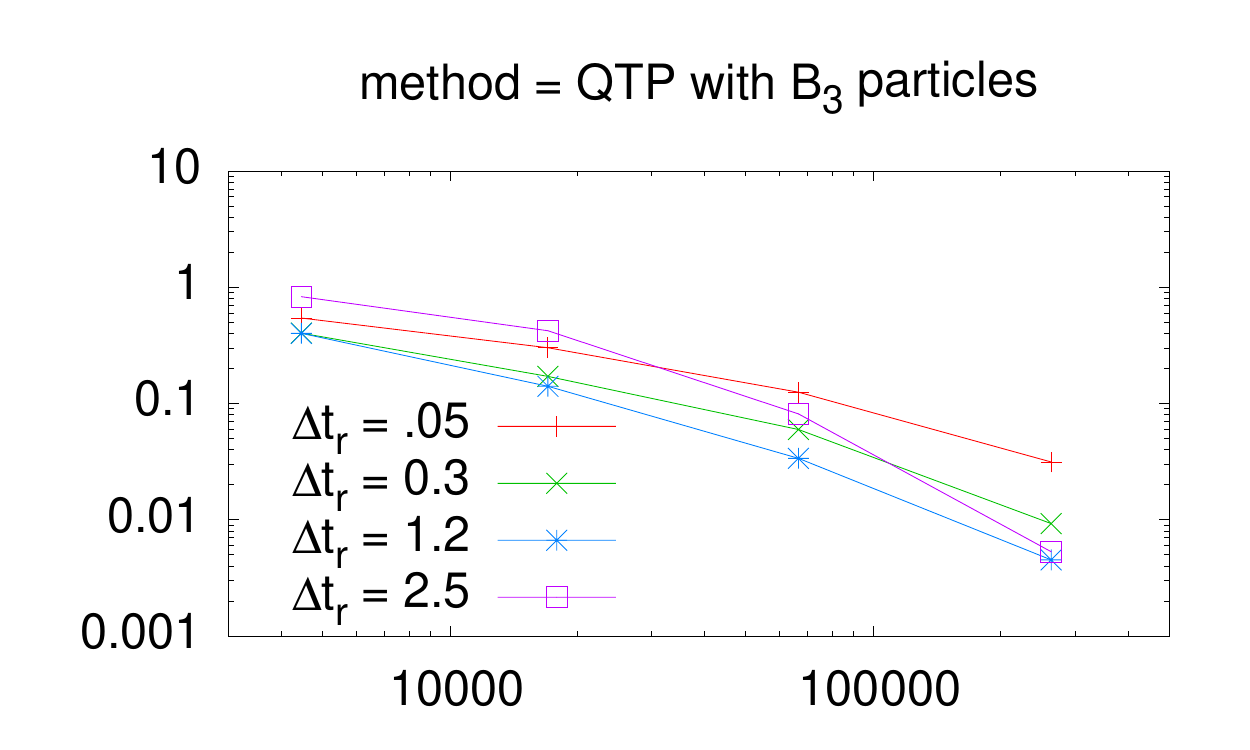}
\end{tabular}
  \caption{(Color) Convergence curves (relative $L^\infty$ errors at $t=T$ vs. average number of active particles)
  for the reversible test case SW-hump defined in Table~\ref{tab:test-cases}, solved with the different methods (see text for details). 
  The first row shows the profile of the exact solution:
  the initial (and final) density $f^0=f(T)$ is on the left, whereas the intermediate solution $f(T/2)$ 
  (with maximum stretching) is on the right.
  }
  \label{fig:cc-SW-hump}
 \end{center}
\end{figure}

\begin{figure} [!ht]
\begin{center}
\begin{tabular}{cc}
\hspace{30pt}
\includegraphics[height=0.3\textwidth]{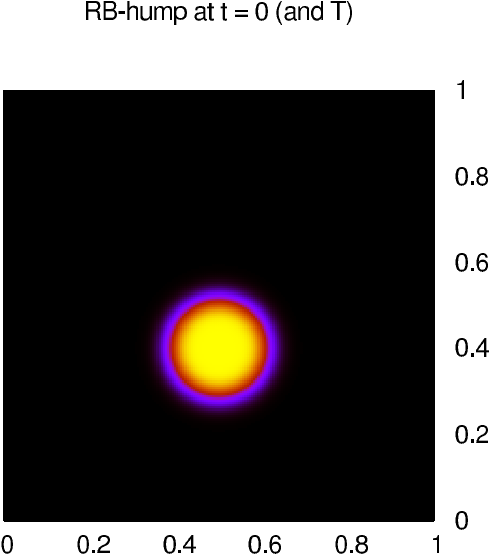}
& % \hspace{5pt} 
\hspace{20pt}
\includegraphics[height=0.3\textwidth]{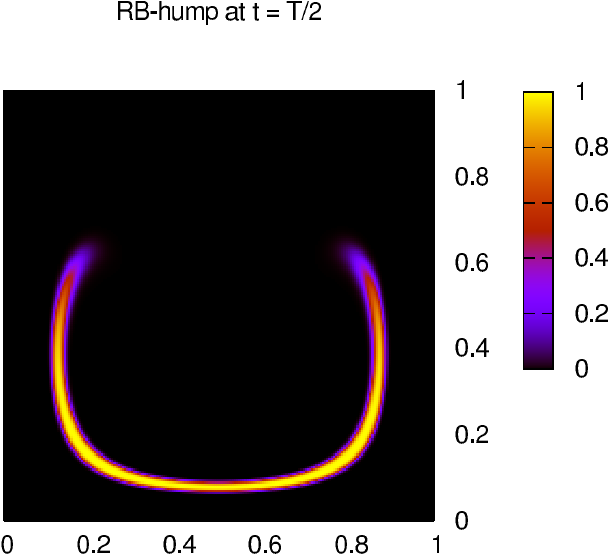}
\vspace{5pt}
\\
\includegraphics[width=0.45\textwidth]{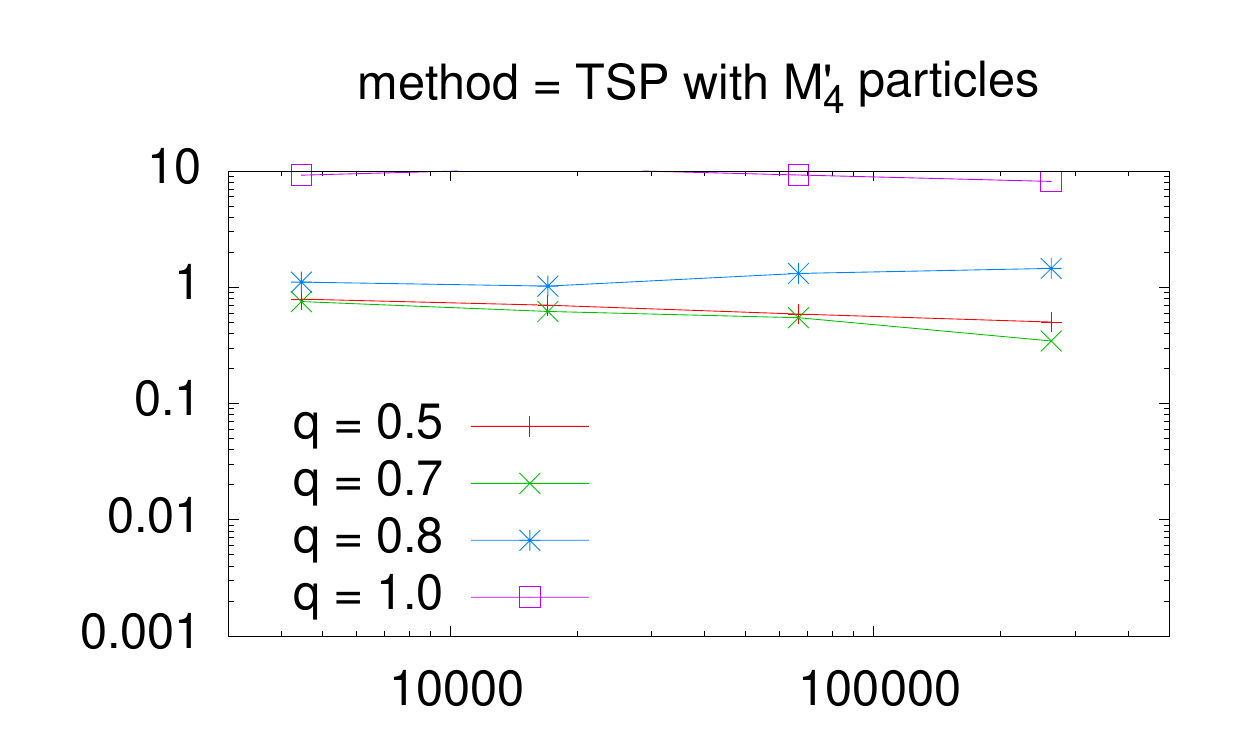}
& % \hspace{5pt} 
\includegraphics[width=0.45\textwidth]{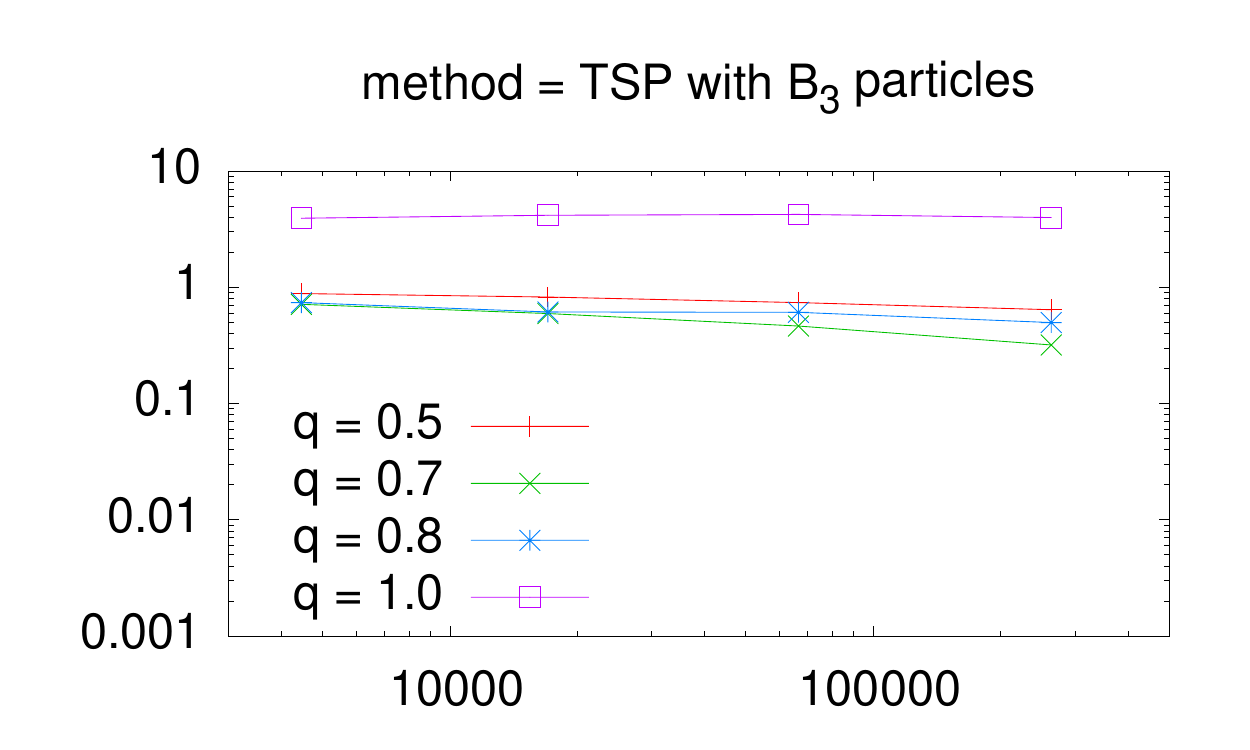}
\vspace{5pt}
\\
\includegraphics[width=0.45\textwidth]{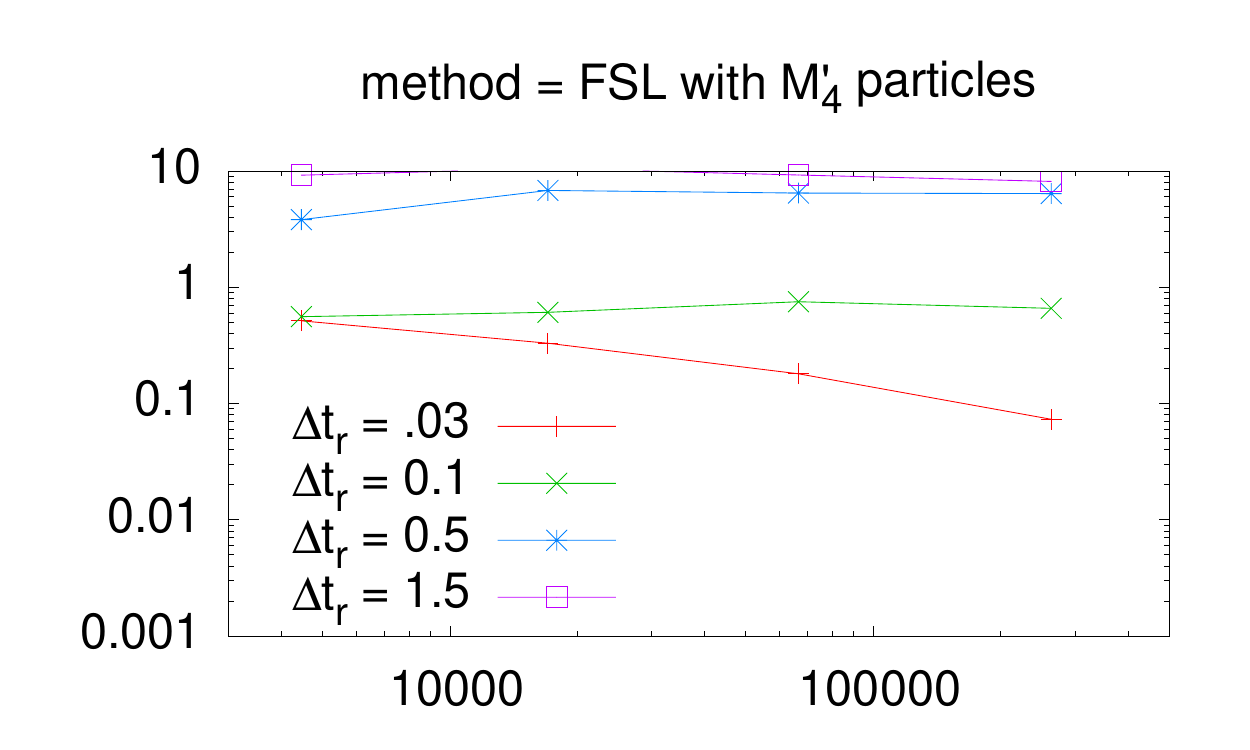}
& % \hspace{5pt} 
\includegraphics[width=0.45\textwidth]{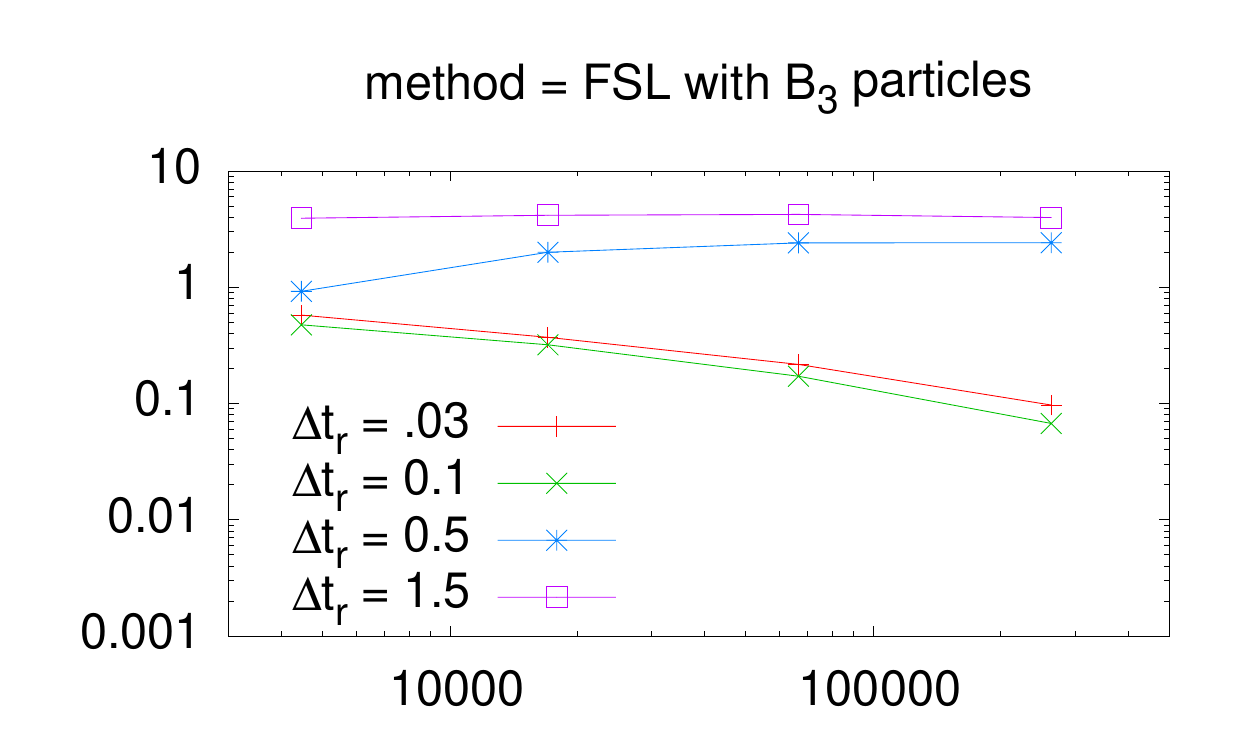}
\vspace{5pt}
\\
\includegraphics[width=0.45\textwidth]{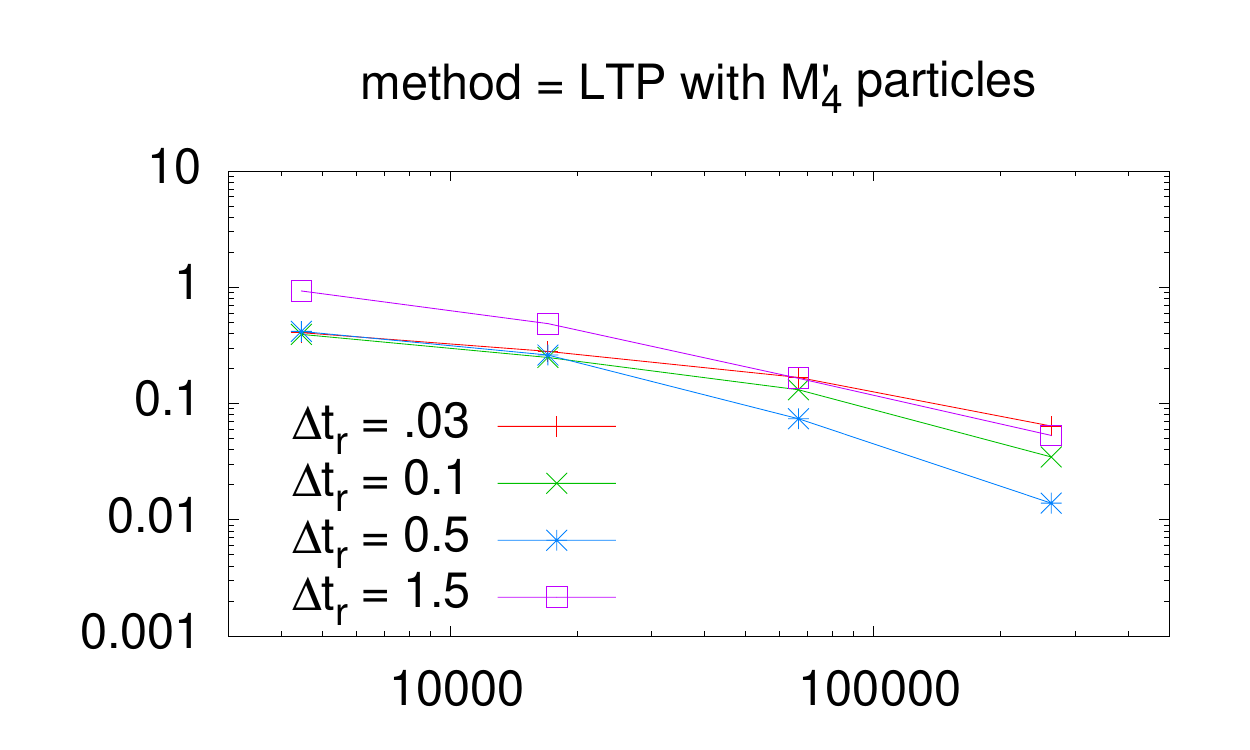}
& % \hspace{5pt} 
\includegraphics[width=0.45\textwidth]{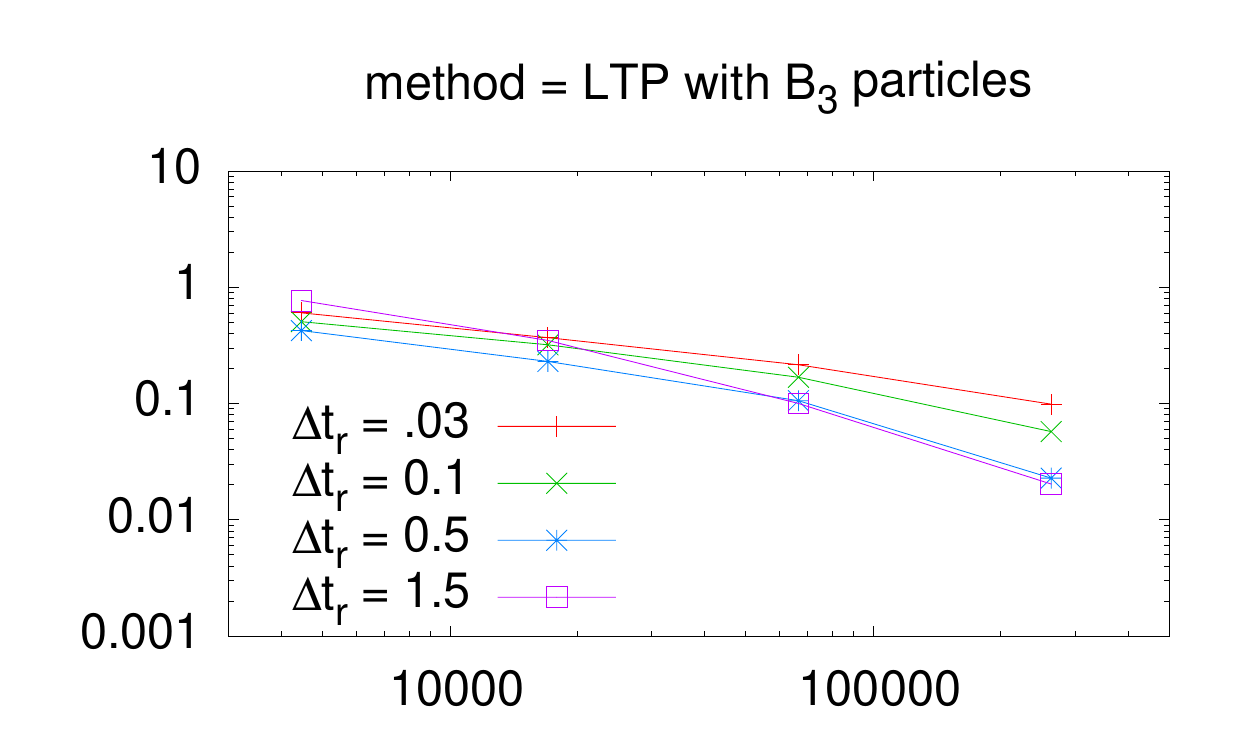}
\vspace{5pt}
\\
\includegraphics[width=0.45\textwidth]{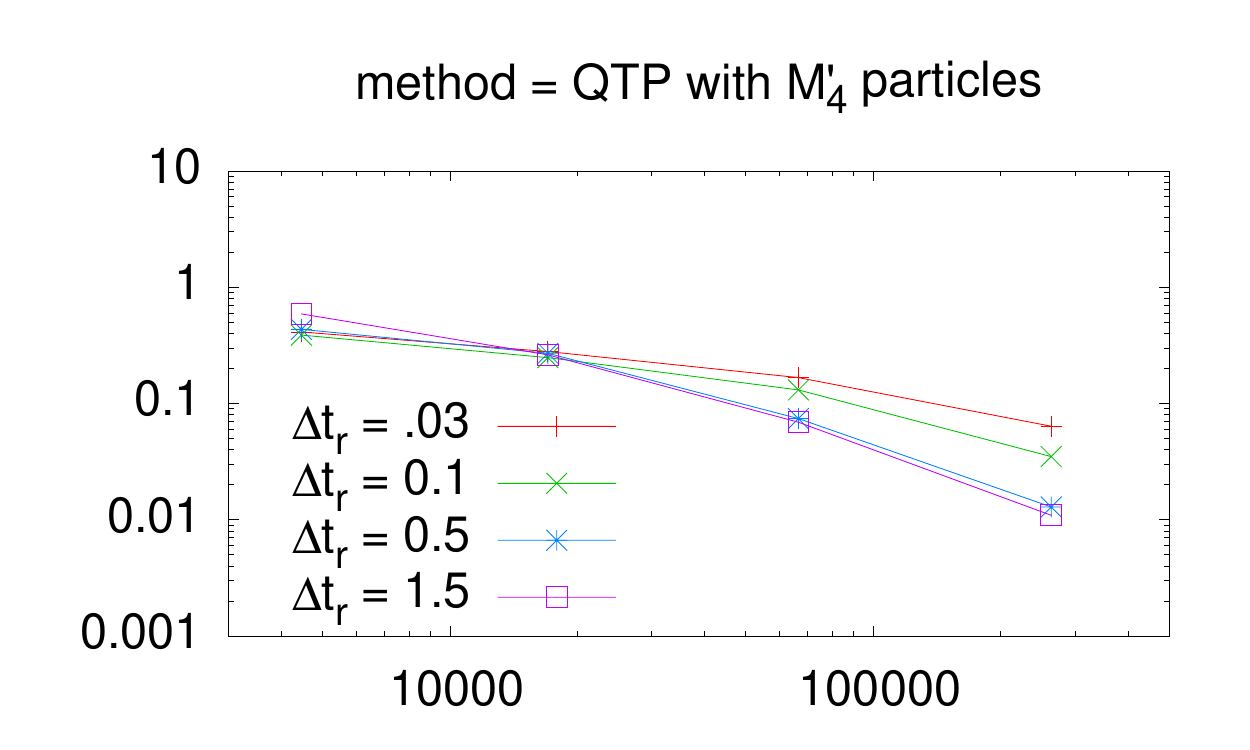}
& % \hspace{5pt} 
\includegraphics[width=0.45\textwidth]{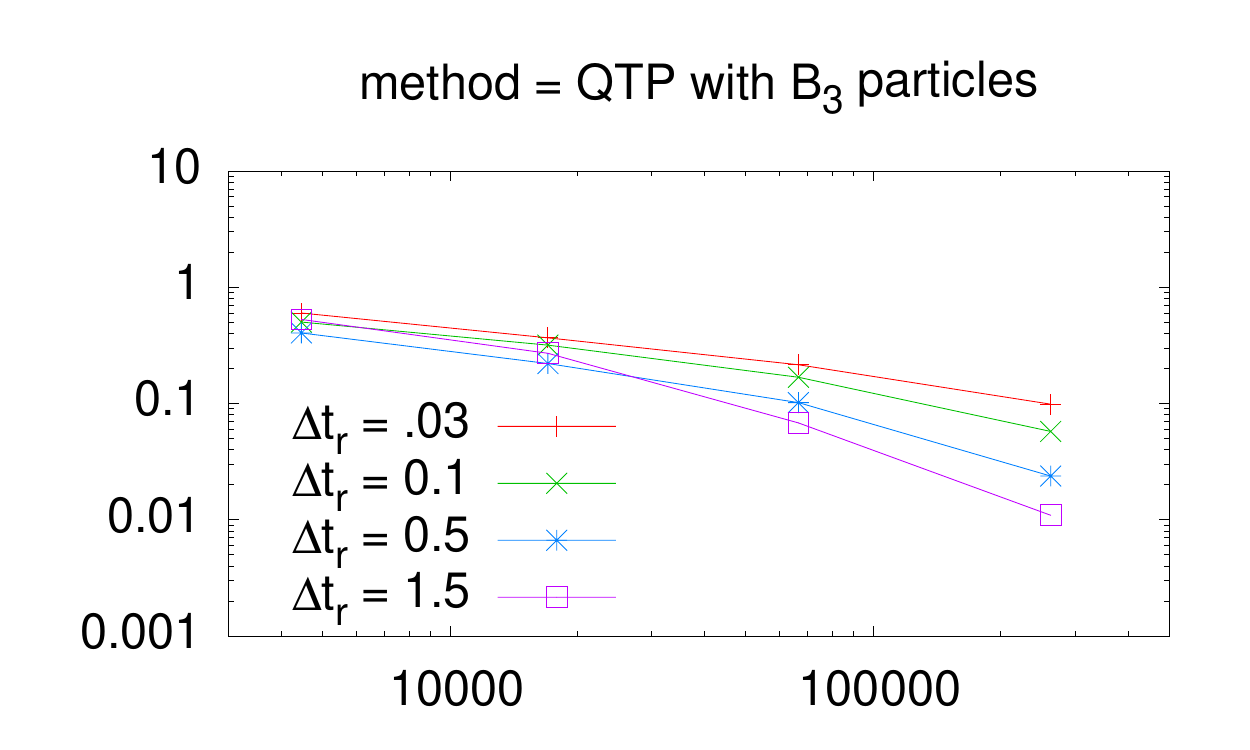}
\end{tabular}
  \caption{(Color) Convergence curves (relative $L^\infty$ errors at $t=T$ vs. average number of active particles)
  for the reversible test case RB-hump defined in Table~\ref{tab:test-cases}, solved with the different methods (see text for details). 
  The first row shows the profile of the exact solution:
  the initial (and final) density $f^0=f(T)$ is on the left, whereas the intermediate solution $f(T/2)$ 
  (with maximum stretching) is on the right.
  }
  \label{fig:cc-RB-hump}
 \end{center}
\end{figure}

\begin{figure} [!ht]
\begin{center}
\begin{tabular}{cc}
\hspace{30pt}
\includegraphics[height=0.3\textwidth]{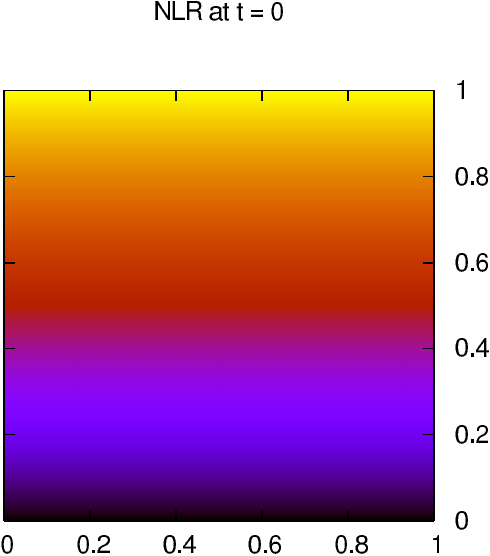}
& % \hspace{5pt} 
\hspace{20pt}
\includegraphics[height=0.3\textwidth]{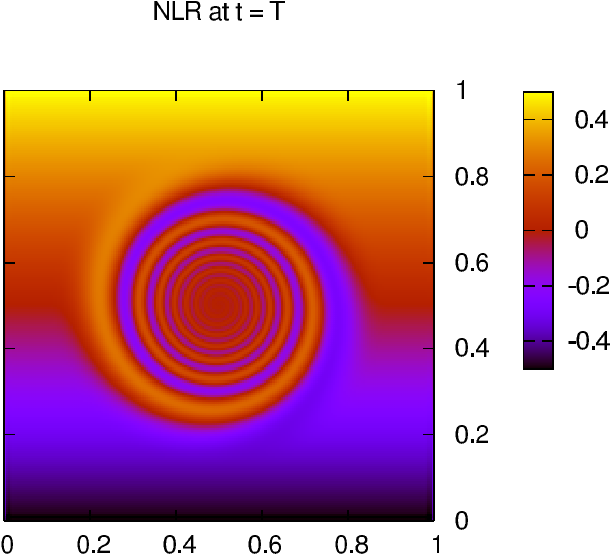}
\vspace{5pt}
\\
\includegraphics[width=0.45\textwidth]{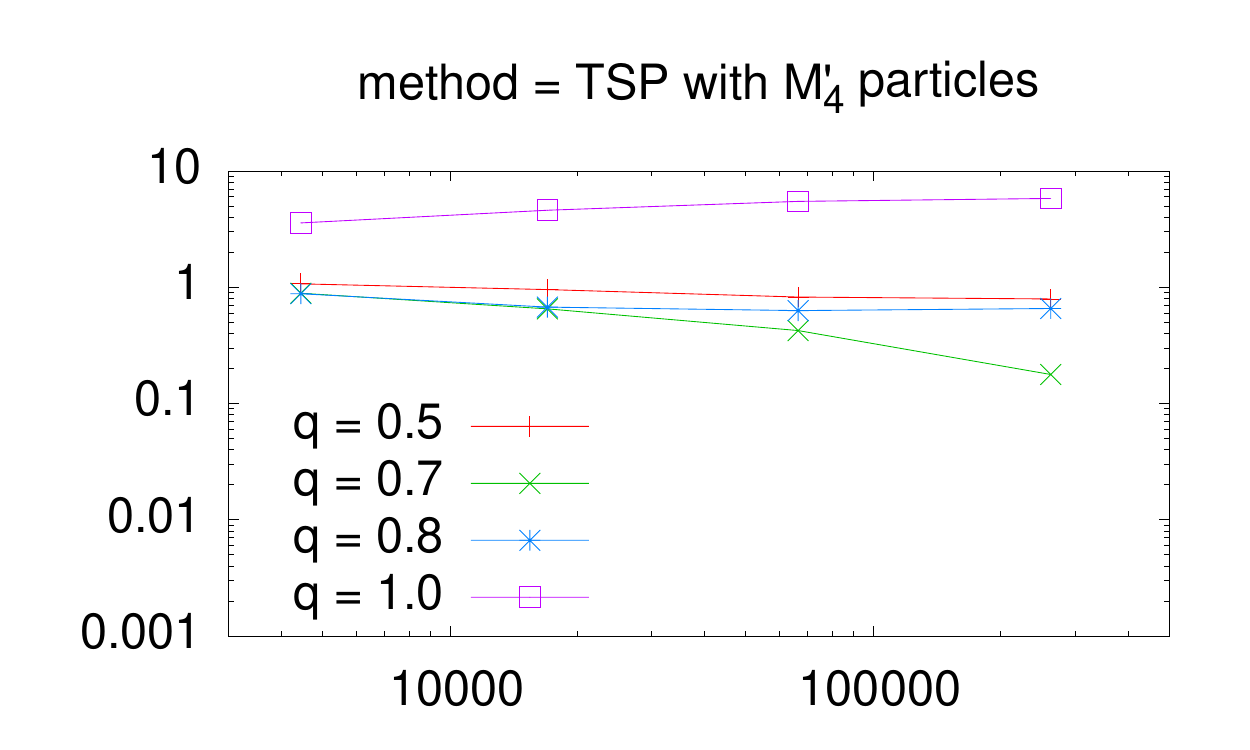}
& % \hspace{5pt} 
\includegraphics[width=0.45\textwidth]{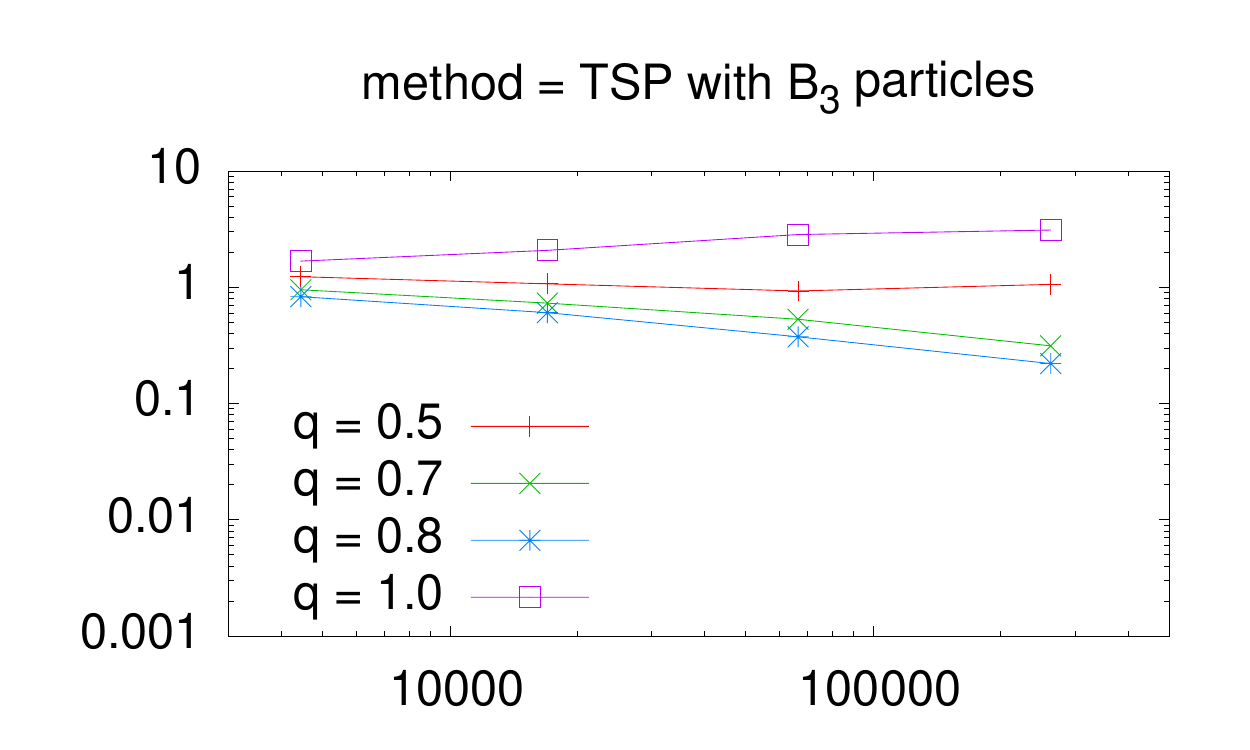}
\vspace{5pt}
\\
\includegraphics[width=0.45\textwidth]{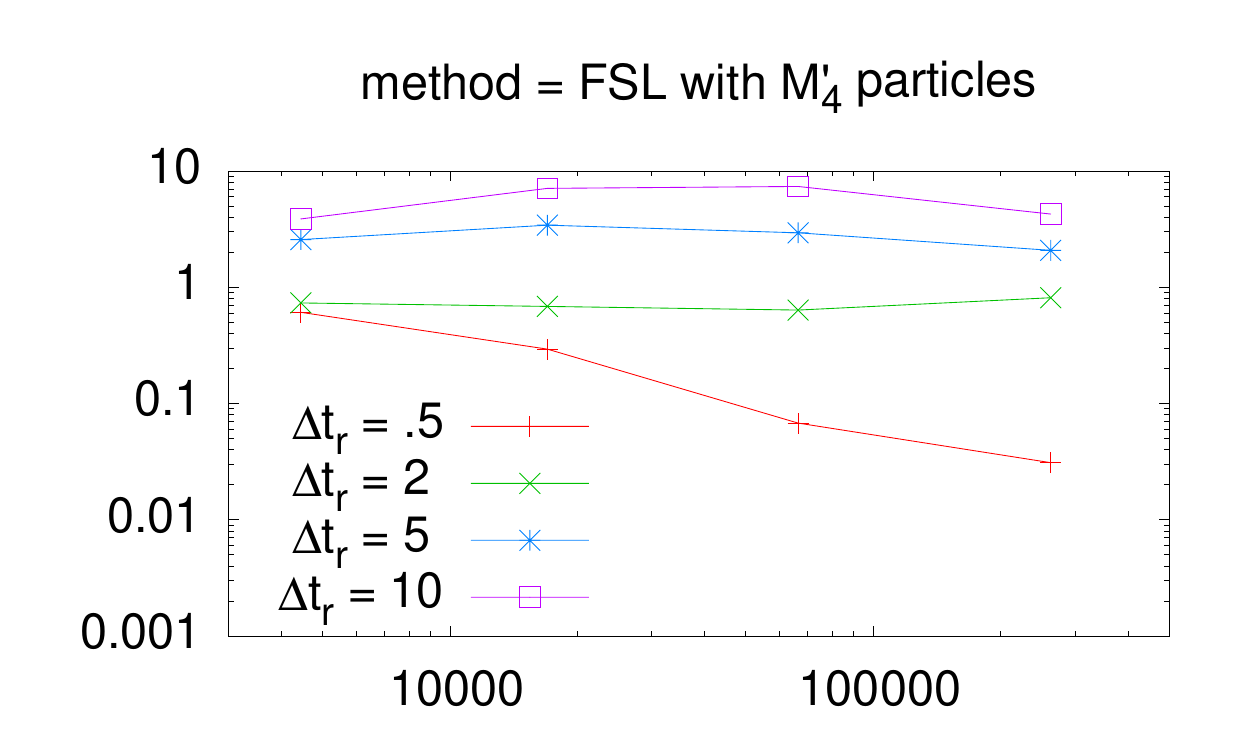}
& % \hspace{5pt} 
\includegraphics[width=0.45\textwidth]{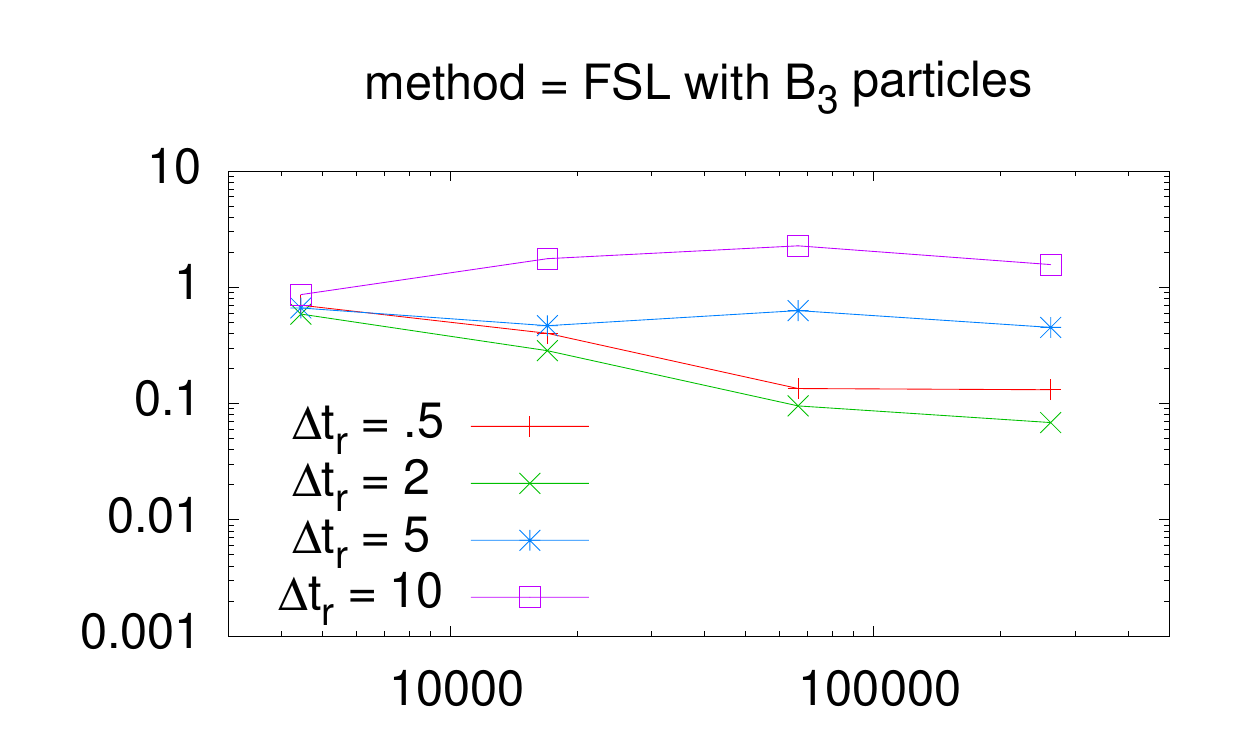}
\vspace{5pt}
\\
\includegraphics[width=0.45\textwidth]{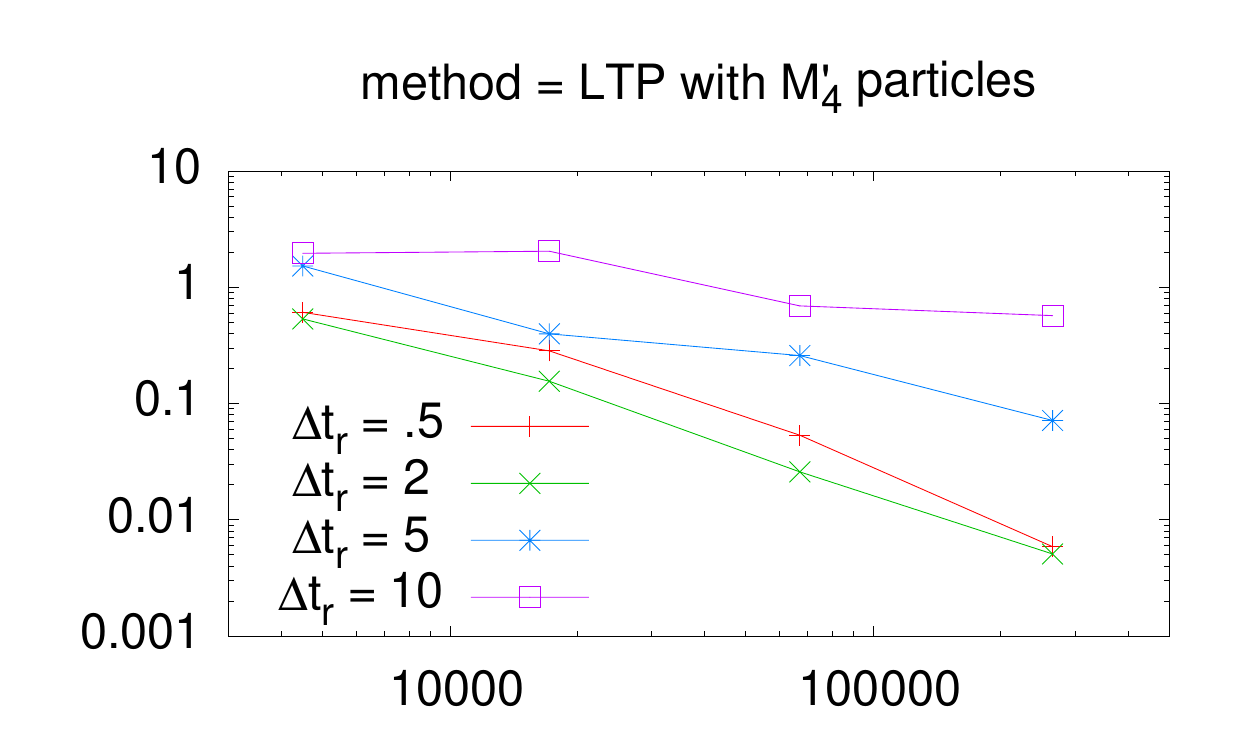}
& % \hspace{5pt} 
\includegraphics[width=0.45\textwidth]{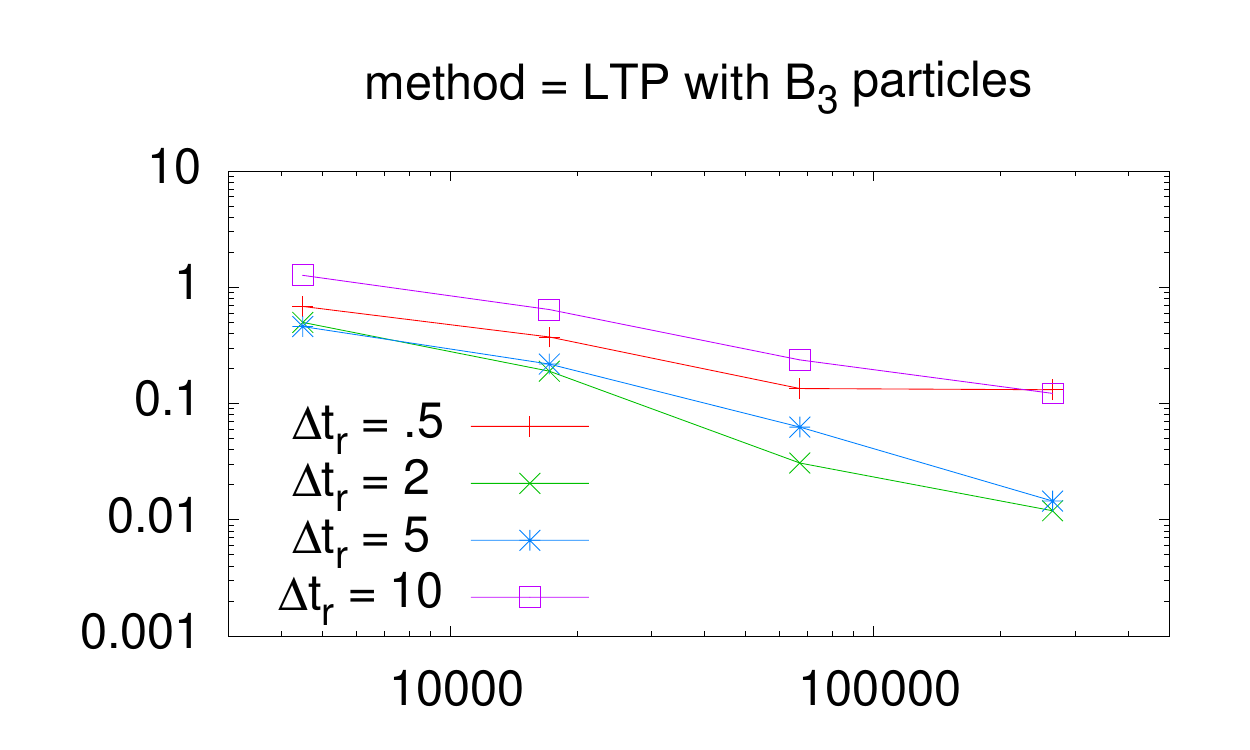}
\vspace{5pt}
\\
\includegraphics[width=0.45\textwidth]{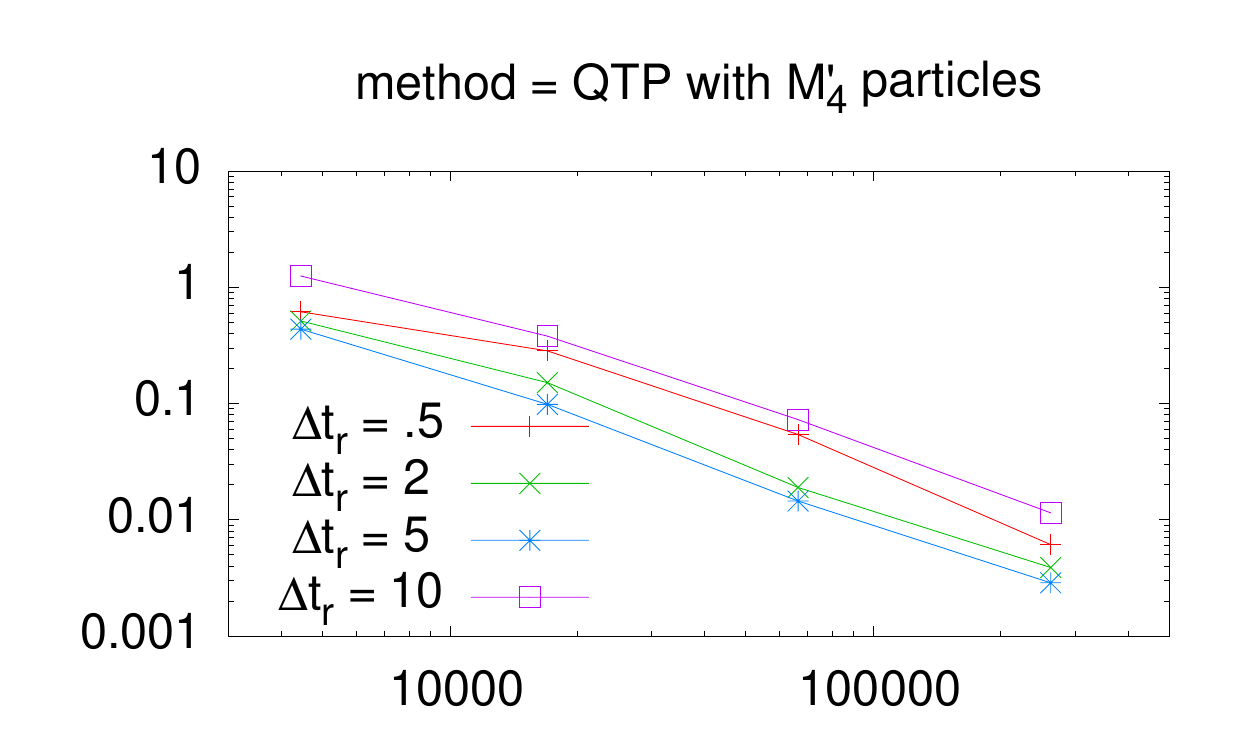}
& % \hspace{5pt} 
\includegraphics[width=0.45\textwidth]{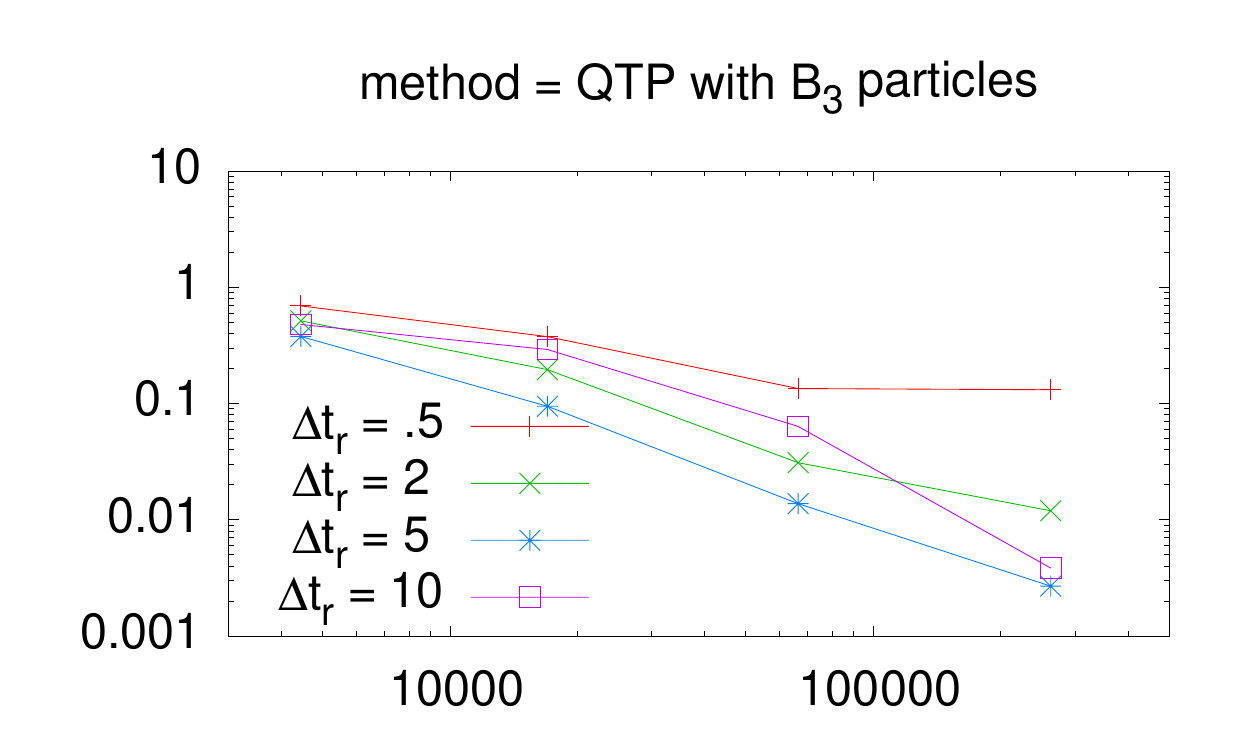}
\end{tabular}
  \caption{(Color) Convergence curves (relative $L^\infty$ errors at $t=T$ vs. average number of active particles)
  for the non-reversible test case NLR defined in Table~\ref{tab:test-cases}, solved with the different methods (see text for details). 
  The first row shows the profile of the exact solution:
  the initial density $f^0$ is on the left and the final solution $f(T)$ 
  (with maximum stretching) is on the right.
  }
  \label{fig:cc-NLR}
 \end{center}
\end{figure}

To measure the convergence properties of the proposed methods, we plot in Figures~\ref{fig:cc-SW-cone}-\ref{fig:cc-NLR} the 
relative $L^\infty$ errors at $t = T$ versus the average number of active particles. Results are shown 
for the four test-cases defined in Table~\ref{tab:test-cases}, and for each case the LTP and QTP methods are compared 
with the standard TSP and FSL methods described in Section~\ref{sec:overview},
using either $M'_4$ or cubic B-spline particles. %, as described in Section~\ref{sec:shapes}.
For the reversible test-cases we always remap the particles at $t=T/2$ when the solutions are stretched most, 
in order to take into account the accuracy of the intermediate approximations in the final measurements.

For every method we also compare different runs obtained by varying their main parameter: with the TSP method we take different
values of the exponent $q$ for which the particle scale $\ve$ behaves like $h^q$ (smaller values of $q$ corresponding 
to more particle overlapping, see Section~\ref{sec:tsp}), and with the remapped particle schemes we test 
different values of the remapping period $\Dtr = \Nr \Dt$.
From these results we make the following observations.
\begin{itemize}
\item[$\bullet$]
    As expected from the theory, the TSP runs (shown on the first rows)
    only converge for values of $q$ smaller that 1, which corresponds to an extended overlapping: the ratio 
    $\ve/h$ must go to $+\infty$ as $h$ goes to 0. We also observe that the numerical convergence is always slow when compared to 
    the remapped particle methods, despite the fact that the $M'_4$ kernel satisfies a moment condition of order 4, 
    see Estimate~\eqref{stand-est}.
\item[$\bullet$]
    In order to converge, the FSL method (second rows) must be run with very short 
    remapping periods, making it somehow closer to an Eulerian method. When remapped every few time steps indeed, the FSL runs 
    exhibit significantly faster convergence rates than the TSP method. However, for fixed values of $\Dtr$ there 
    is always a point where the convergence stops. This amounts to asking for more remappings when
    more particles are used, which in practice has an effect similar to imposing a CFL constraint.
\item[$\bullet$]
    The most striking results with the LTP and QTP runs is that they seem to completely suppress the loss of convergence observed 
    with the FSL method. Moreover, the observed behavior is now radically different: not only does the convergence always hold, 
    it is {\em improved} when the remapping period is increased up to a certain value that is significantly larger than the time step.
    Specifically, we observe that the accuracy of the LTP runs improves for remapping periods as large as 10 to 50 times the time step,
    depending on the test cases. And with the QTP method these ratios go up to values between 30 and 50. 
\item[$\bullet$]
    Broadly speaking, the above hold true for both the $M'_4$ and the $B_3$ particles, the main difference being
    that the former perform better for small remapping periods. 
%    while cubic spline particles give better results for large remapping periods. 
    Given the fact that they are remapped with lower-order but also lower dissipation than the $B_3$ particles, this behavior is rather expected.
\end{itemize}
The remaining sections are devoted to further investigating the influence of the remapping period $\Dtr$ on the accuracy.

%%%%%%%%%
%%%%%%%%%

%%%%%%%%%
%%%%%%%%%

%%%%%%%%%
%%%%%%%%%

  %%%%%%%%%
%%%%%%%%%

%%% %%% %%% %%% %%% %%% %%% %%% %%% %%%
%%% %%% %%% %%% %%% %%% %%% %%% %%% %%%
\subsection{Influence of the remapping period}

\begin{figure} [!ht]
\begin{center}
\begin{tabular}{cc}
\includegraphics[width=0.45\textwidth]{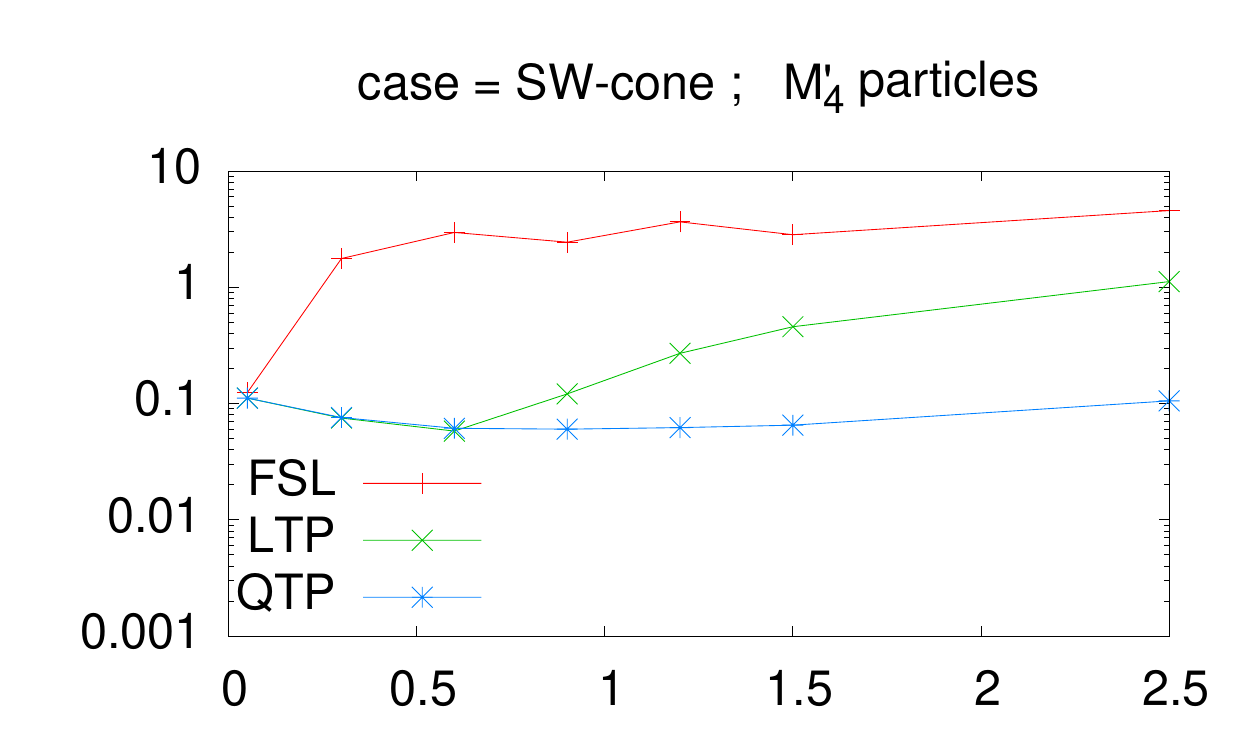}
& % \hspace{5pt} 
\includegraphics[width=0.45\textwidth]{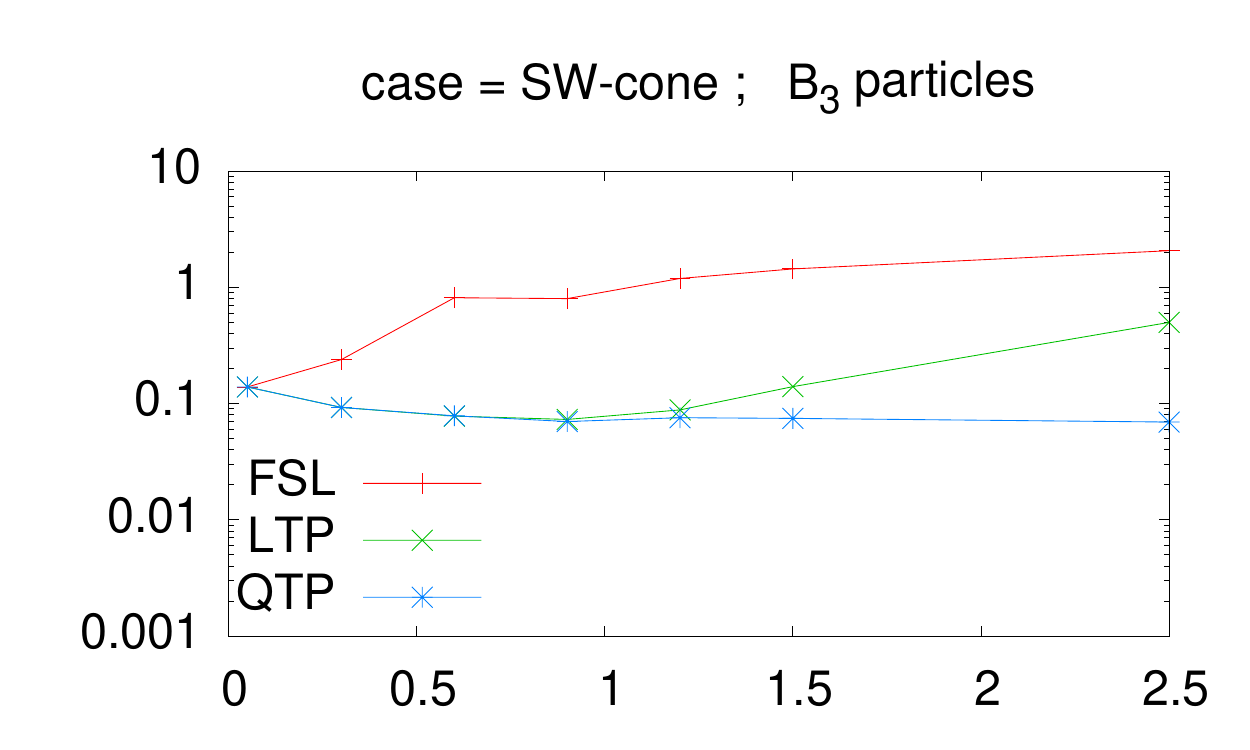}
\vspace{10pt}
\\
\includegraphics[width=0.45\textwidth]{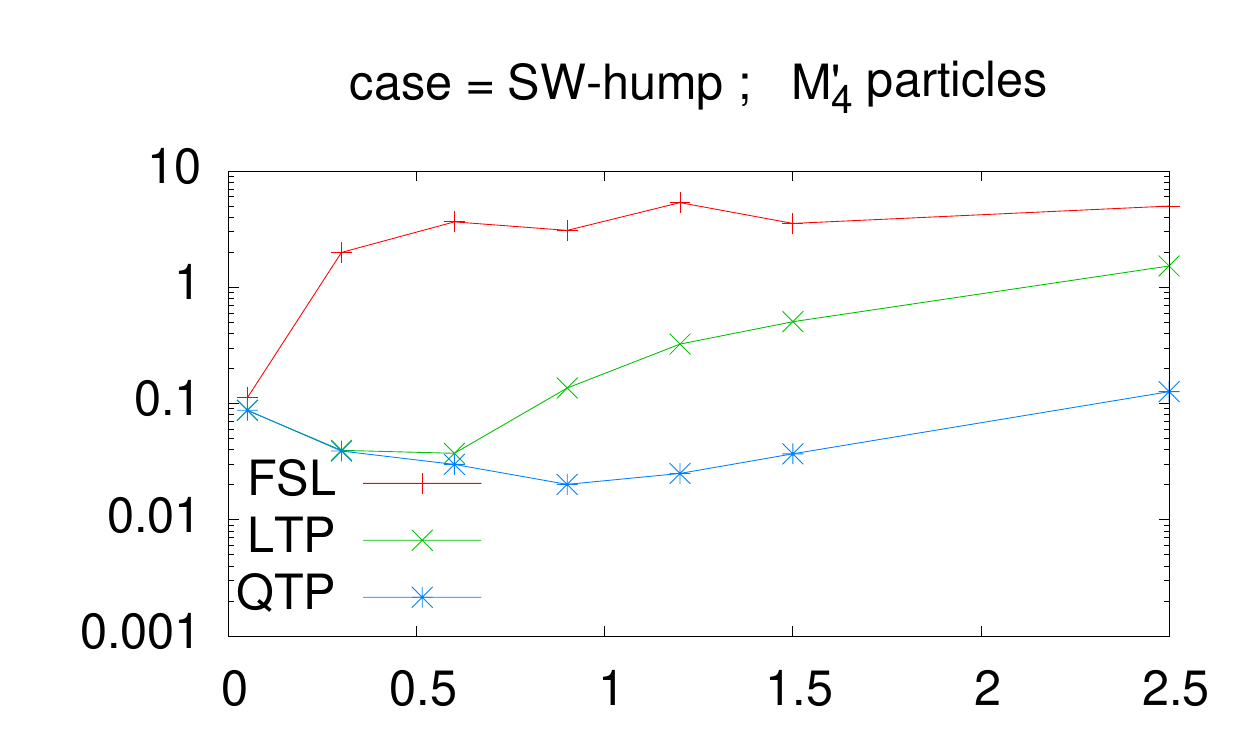}
& % \hspace{5pt} 
\includegraphics[width=0.45\textwidth]{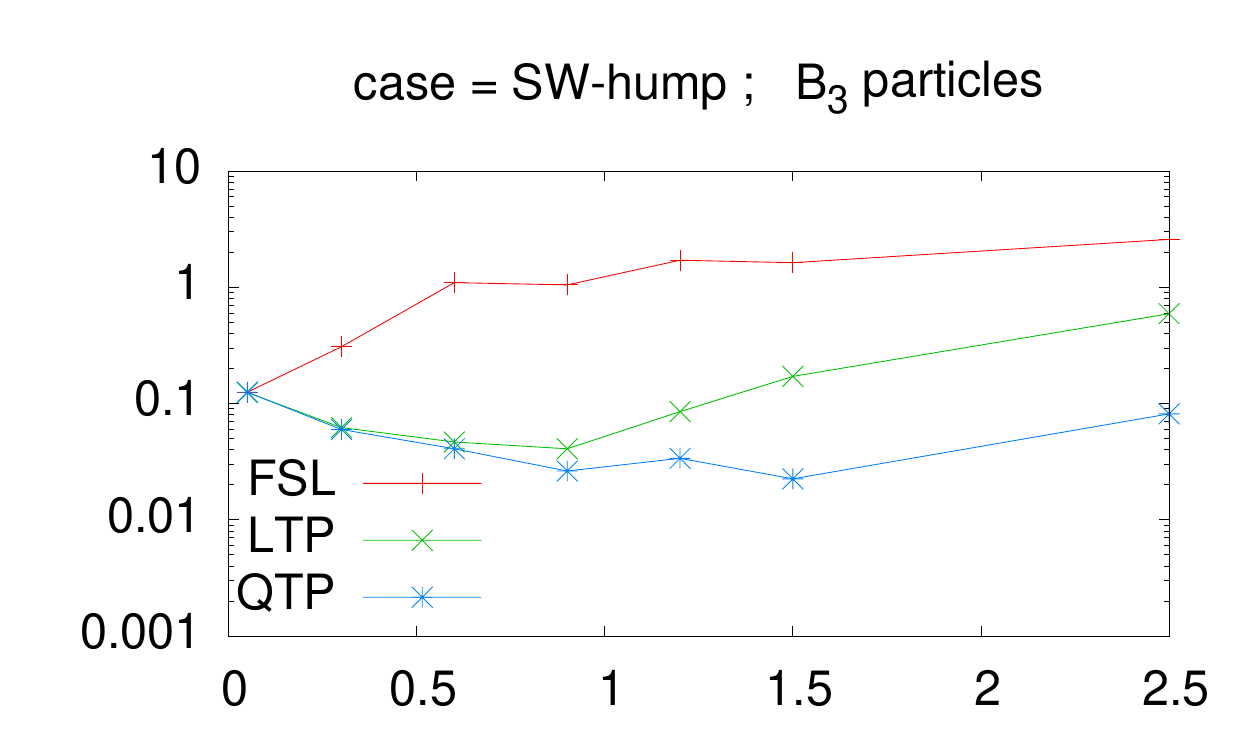}
\vspace{10pt}
\\
\includegraphics[width=0.45\textwidth]{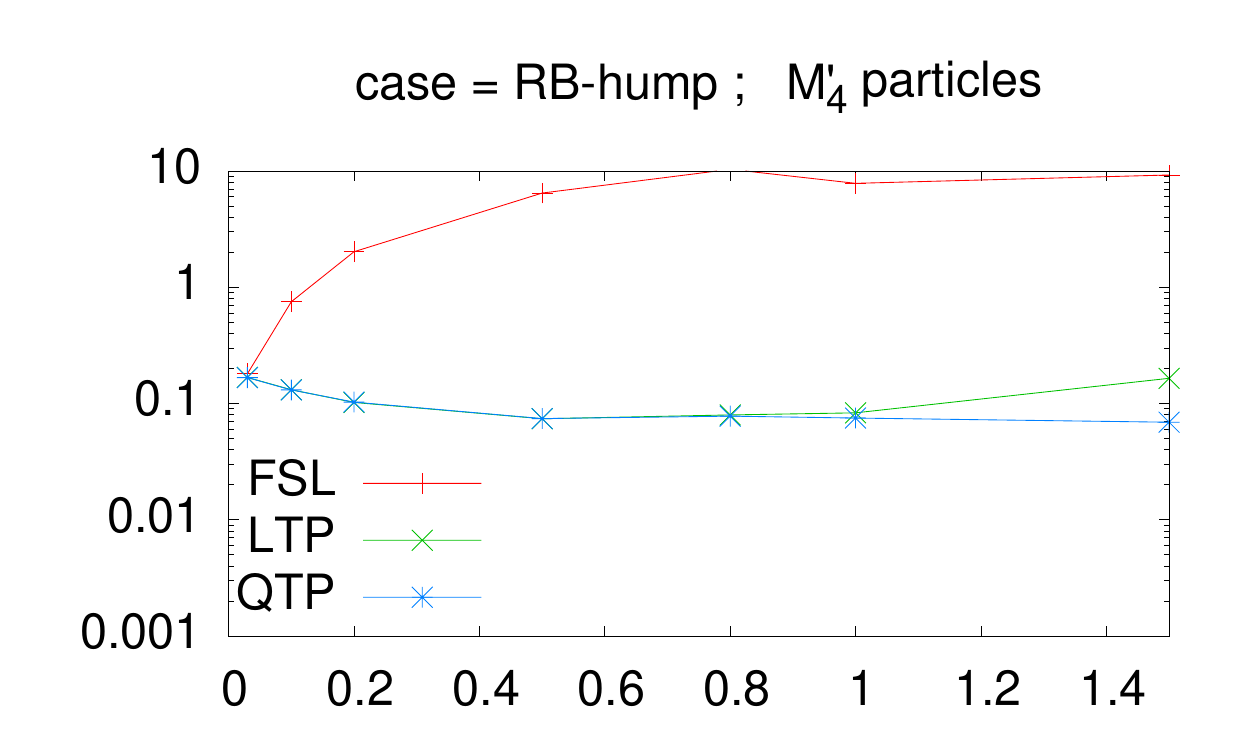}
& % \hspace{5pt} 
\includegraphics[width=0.45\textwidth]{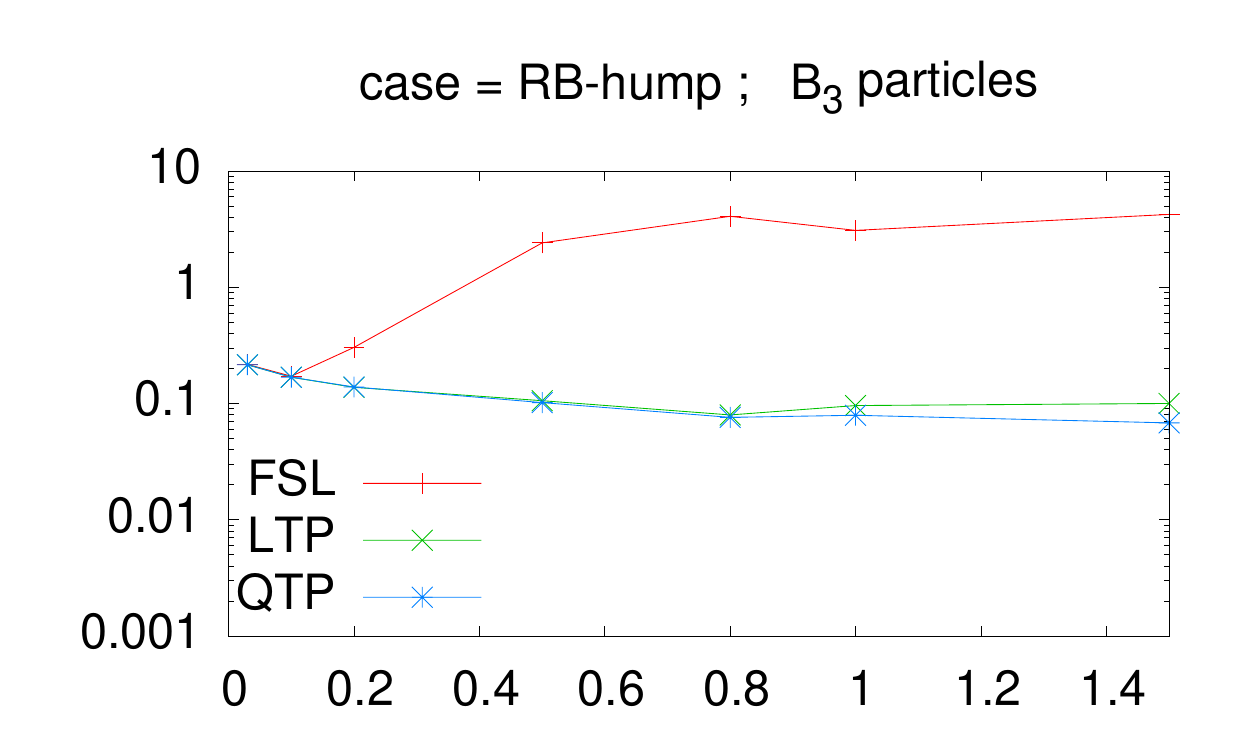}
\vspace{10pt}
\\
\includegraphics[width=0.45\textwidth]{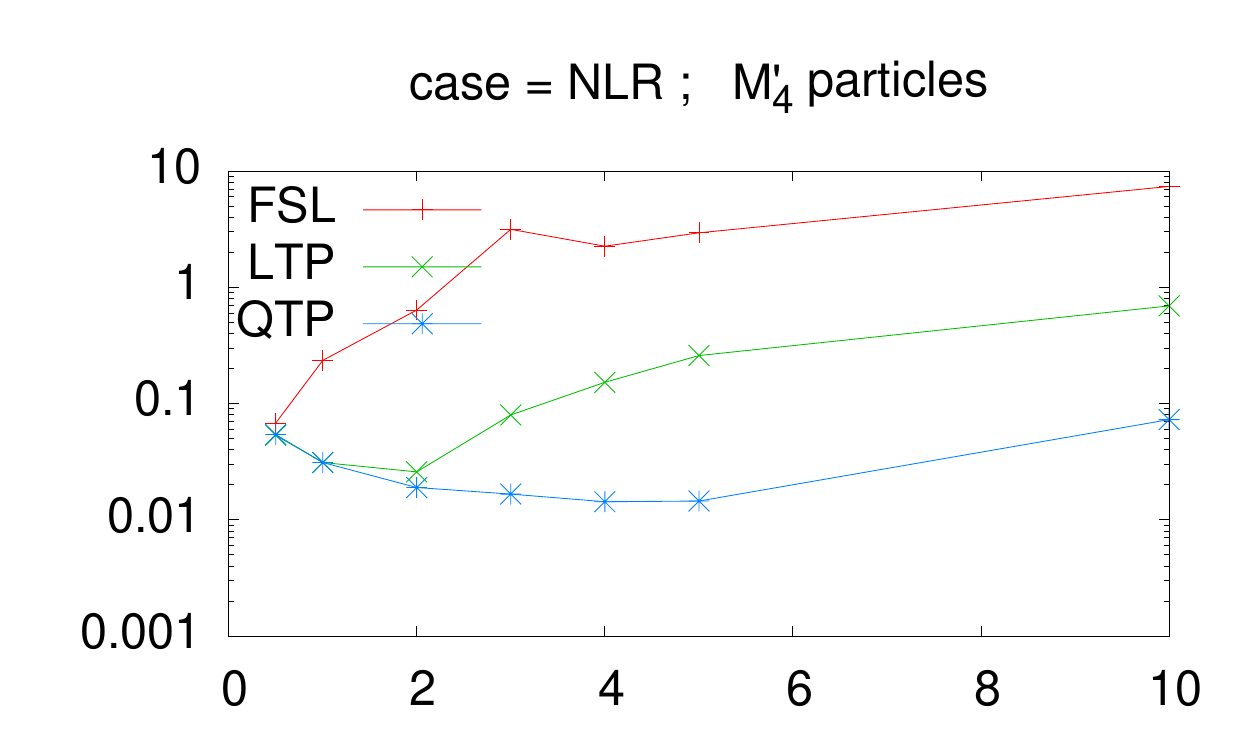}
& % \hspace{5pt} 
\includegraphics[width=0.45\textwidth]{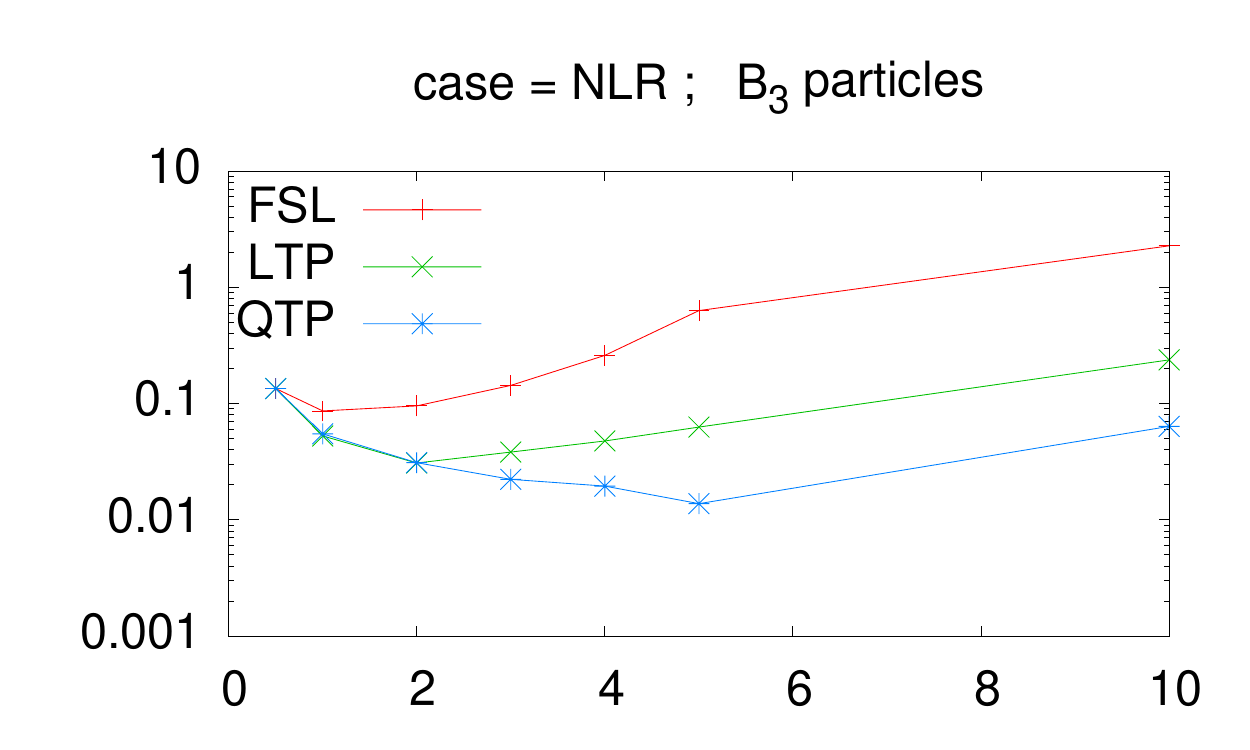}
\end{tabular}
  \caption{(Color) Relative $L^\infty$ errors at $t=T$ vs. remapping period $\Dtr$ for the different test cases
  solved with remapped particle methods of order 0 (FSL), 1 (LTP) and 2 (QTP).
  Here the particles are initialized and remapped on a cartesian grid with $h = 2^{-8}$, which correspond to a maximum of 
  $256 \times 256$ particles (in some cases only a fraction of these particles are activated). 
  Left panels show results obtained with $M'_4$ particles, whereas right panels show results obtained with cubic B-spline particles.
  Qualitatively similar curves have been obtained with runs using finer grids of $512 \times 512$ particles, not shown here.}
  \label{fig:err-Dtr}
 \end{center}
\end{figure}

In order to better understand how the numerical accuracy evolves with respect to the remapping period,
we plot in Figure~\ref{fig:err-Dtr} the final errors obtained with a series of runs using increasing values for $\Dtr$. 
This leads us to the following observations.
\begin{itemize}
\item[$\bullet$]
    With the FSL scheme the particles must be remapped each few time steps, as the accuracy of the method
    quickly deteriorates for increasing remapping periods. This confirms the previous observations.
\item[$\bullet$]     
    By transforming the particles either linearly or quadratically, we obtain a twofold benefit: first, 
    the accuracy is always improved, and sometimes significantly. Second, the optimal accuracy now 
    corresponds to some trade off between small remapping periods where the remappings error 
    dominate the transport errors, and large ones where the opposite occurs.
\item[$\bullet$]
    By transforming the particles quadratically instead of linearly, we observe some further improvements in the accuracy, 
    but the major benefit seems to be a significant gain in robustness with respect to the remapping period.
    In the NLR test-case for instance, the optimal remapping period for the QTP scheme is about five time larger than for the LTP
    scheme, be it with the moderate runs shown in Figure~\ref{fig:err-Dtr}, or with finer ones using $512 \times 512$ particles, 
    not shown here.
\end{itemize}
%    
%%%%%%%%%
%%%%%%%%%

%%% %%% %%% %%% %%% %%% %%% %%% %%% %%%
%%% %%% %%% %%% %%% %%% %%% %%% %%% %%%
\subsection{Dynamic remapping strategies}

In this section we propose and study a practical tool for automatically selecting the time steps
where the particles should be remapped. Although our indicator does not {\em always} perform as well
as the optimal static strategy, it gives very satisfactory results when tested on our different 
benchmark problems, and it is fully local in the sense that it only involves particle-wise 
computations, which best fits parallel solvers. 

To describe it let us denote by $n_i$, $i = 0, \ldots , {R-1}$, the time steps of the initial 
approximation and further remappings, as selected by a given strategy.
Then if $T_h^{n,m} = T_h^{m-1} \cdots T_h^n$ denotes the transport operator acting without remappings 
between the times $t^n$ and $t^m$, the remapped particle scheme that maps $f^0$ to $f^N_h$ reads
(writing $n_R = N$ and $f^0_h = f^0$ for simplicity)
$$
f^N_h = f^{n_R}_h = S^{n_R}_h  f^0_h =T^{n_{R-1},n_R}_h A_h  f^{n_{R-1}}_h = (T^{n_{R-1},n_R}_h A_h) \cdots (T^{n_0,n_1}_h A_h) f^0_h
$$
and the corresponding error $e^N_h := \norm{S^{N}_h-T^{0,N}_\ex}_{L^\infty}$ satisfies 
\begin{align*}
e^N_h &\le \norm{(T^{n_{R-1},n_R}_h - T^{n_{R-1},n_R}_\ex) A_h f^{n_{R-1}}_h}_{L^\infty}
            + \norm{T^{n_{R-1},n_R}_\ex( A_h - I) f^{n_{R-1}}_h}_{L^\infty} 
\\
      & \mspace{200mu}      
            + \norm{T^{n_{R-1},n_R}_\ex( S_h^{n_{R-1}} - T^{0,n_{R-1}}_\ex) f^0}_{L^\infty} 
\\
      &\le \norm{(T^{n_{R-1},n_R}_h - T^{n_{R-1},n_R}_\ex) A_h f^{n_{R-1}}_h}_{L^\infty}
            + \norm{( A_h - I) f^{n_{R-1}}_h}_{L^\infty} 
            + e^{n_{R-1}}_h 
\\
      &\le \sum_{i = 0}^{R-1} 
            \Big(\norm{(T^{n_i,n_{i+1}}_h - T^{n_i,n_{i+1}}_\ex) A_h f^{n_i}_h}_{L^\infty}
            + \norm{( A_h - I) f^{n_i}_h}_{L^\infty} \Big).
\end{align*}
In particular, we see that the global error essentially consists of the transport errors and the remapping ones. 

Our heuristic is then as follows: although estimate~\eqref{Ae-errestim} tells us that the remapping errors
can grow quickly when the smoothness of the solutions deteriorate, in practice we have observed that
they do not depend much on the selected remapping steps. On the contrary, the transport errors often
increase at a comparatively fast rate when the remapping period grows large. And it is easily seen from 
Theorem~\ref{th:conv-f-T1} and \ref{th:conv-f-Tr} that each remapping resets them to an increasing 
function starting at 0, by resetting the corresponding flow to the identity.
Therefore, it seems reasonable to remap the particles when the estimated transport error becomes larger than
the estimated remapping error. Specifically we propose to remap $f^n_h$ when
\begin{equation} \label{rem-crit}
C_{\rm remap} \cE\Big((T^{n_i,n}_h - T^{n_i,n}_\ex) A_h f^{n_i}_h)\Big) 
    \ge
    \cE\Big((A_h -I) f^{n}_h)\Big) 
\end{equation}
where $n_i$ denotes the last remapping step preceding $n$ and $C_{\rm remap}$ is a parameter to be determined
from numerical experiments. It remains to specify local indicators for the transport and remapping errors.

Considering for simplicity the direct approach described in Section~\ref{sec:discrete},
we use Estimates~\eqref{errest-eB-T1} and \eqref{errest-eB*-Tr}-\eqref{errest-Theta*}, respectively,
to derive the indicator
\begin{equation} \label{E-transp-err}
\cE\Big((T^{n_i,n}_h - T^{n_i,n}_\ex) A_h f^{n_i}_h)\Big) := 
    \Big(1+\frac{\hat e^n_{B,(1)}(h)}{h}\Big)^d \frac{\hat e^n_{B,(r)}(h)}{h} \norm{f^{n_i}_h}_{L^\infty}
\end{equation}
with $r = 1$ in the LTP case and 2 in the QTP case.
Here the numerical indicators for the backward flow errors are computed as in \eqref{hat-eB},
using the numerical flows \eqref{numflow-ltp} and \eqref{numflow-qtp}.

As for the remapping error, we rely for simplicity on a first order estimate % in \eqref{Ae-errestim},
\begin{equation}\label{first-order}
\norm{(A_h -I) f^{n}_h}_{L^\infty} \lesssim h \abs{f^n_h}_1
\lessapprox h \sum_{j=1}^d \sup_{k \in \ZZ^d} \abs{\partial_j f^n_h(x^n_k)}
\end{equation}
and to approximate the spatial derivatives we write
%$
%\abs{f^n_h}_1 %= \sum_{j=1}^d \norm{\partial_j f^n_h}_{L^\infty} 
%    \approx \sum_{j=1}^d \sup_{k \in \ZZ^d} \abs{\partial_j f^n_h(x^n_k)}
%$
%we approximate 
$$ %\begin{align*}
\partial_j f^n_h(x^n_k)  \approx \partial_j (f^{n_i}_h \circ B^{n_i,n}_\ex ) (x^n_k) 
%\\
%& \approx 
= \sum_{\ell=1}^d \partial_{\ell} f^{n_i}_h (B^{n_i,n}_\ex(x^n_k)) \partial_j (B^{n_i,n}_\ex)_\ell (x^n_k).
$$ %\end{align*}
We next observe that for the remapped particle scheme the approximations involved in \eqref{ltp} read
$x^n_k \approx F^{n_i,n}_\ex(x^0_k)$ and $D^n_k \approx J_{B^{n_i,n}_\ex}(x^n_k)$, we derive the  
following indicator for the remapping error,
\begin{equation} \label{E-remap-err}
\cE\Big((A_h -I) f^{n}_h)\Big) :=
     h \sum_{j=1}^d \sup_{k \in \ZZ^d} \bigabs{ \sum_{\ell=1}^d \partial_{\ell} f^{n_i}_h (x^0_k) (D^n_k)_{\ell,j}}.
\end{equation}
Here the spatial derivatives of $f^{n_i}_h$ can be estimated with finite differences using the values of 
the numerical density on the structured grid $\{ x^0_k : k \in \ZZ^d\}$, since they are needed in the remapping algorithm. 
In particular, we verify that computing the indicators \eqref{E-transp-err} and \eqref{E-remap-err} does not
require any additional communication between the particles, which is an important property for parallel implementations. 

\begin{figure} [!ht]
\begin{center}
\begin{tabular}{cc}
\includegraphics[width=0.45\textwidth]{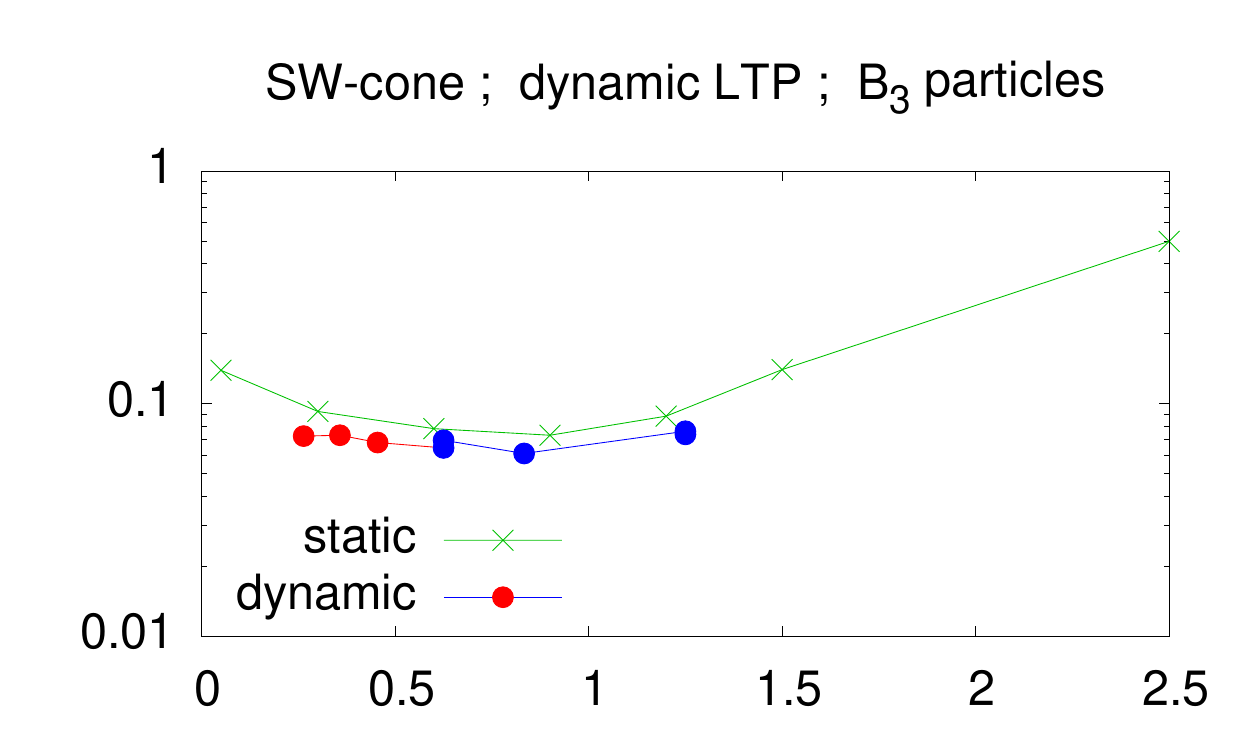}
& % \hspace{5pt} 
\includegraphics[width=0.45\textwidth]{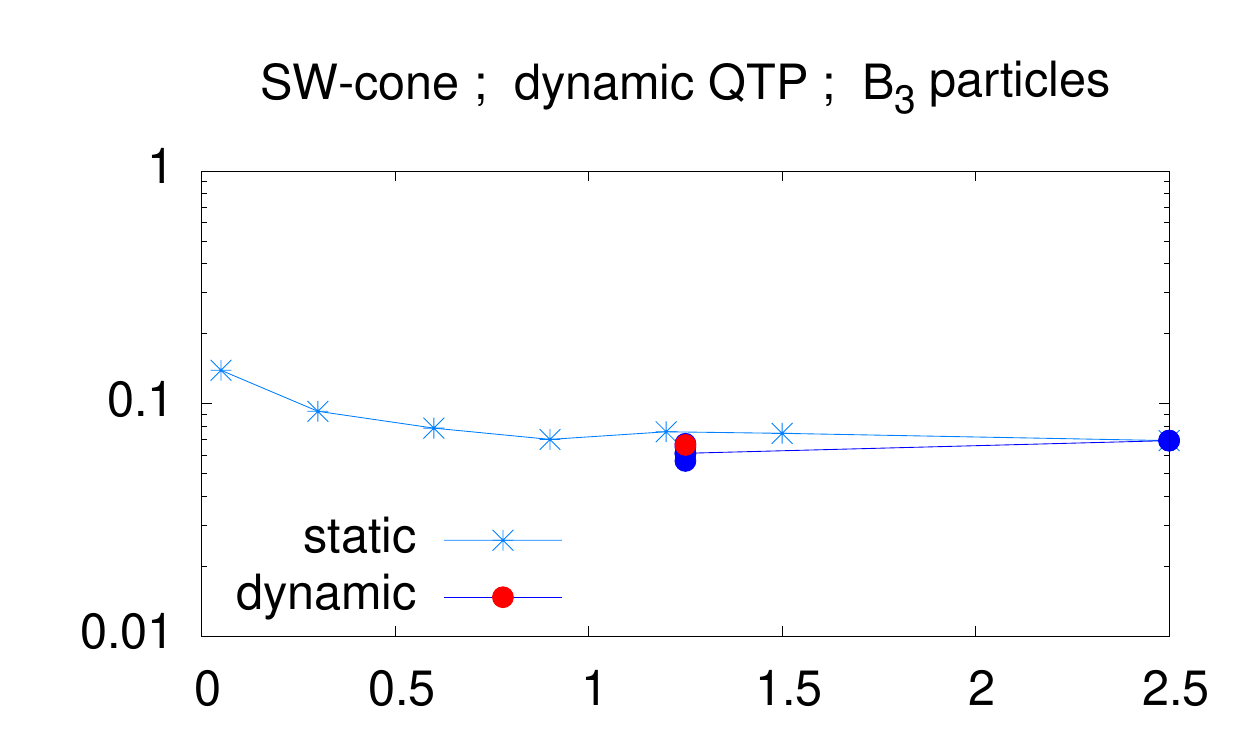}
\vspace{10pt}
\\
\includegraphics[width=0.45\textwidth]{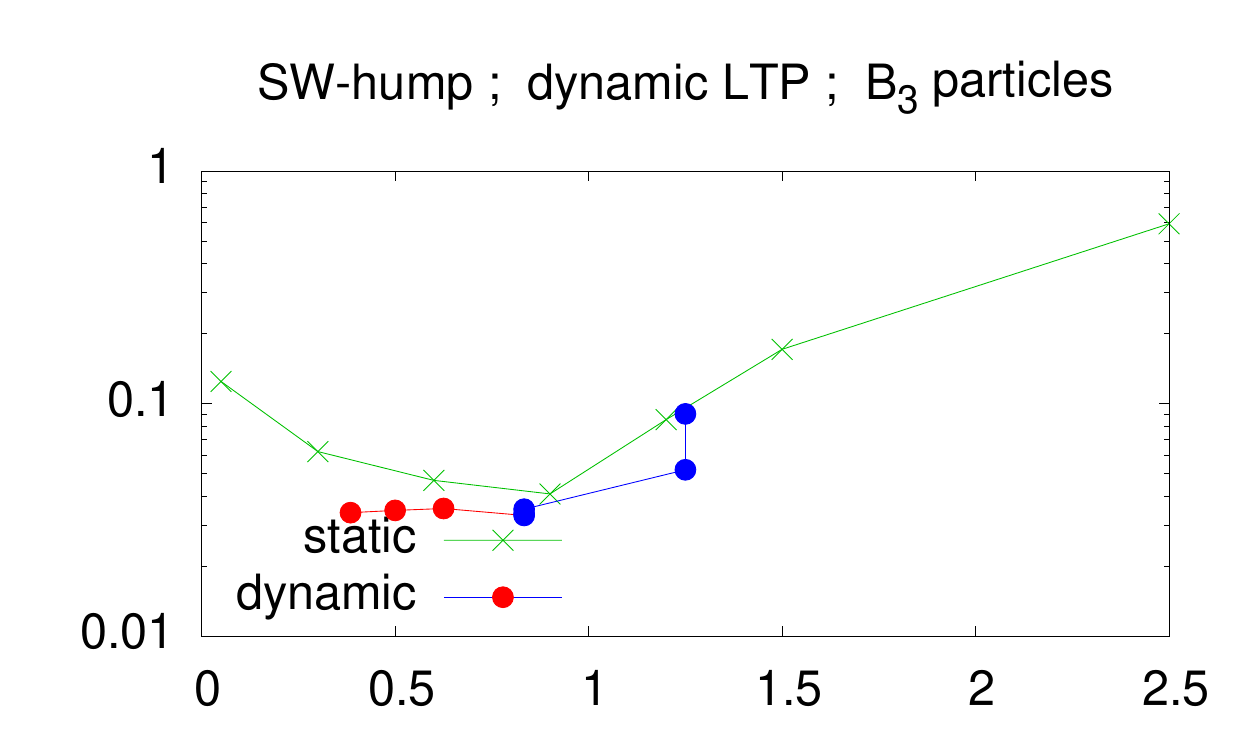}
& % \hspace{5pt} 
\includegraphics[width=0.45\textwidth]{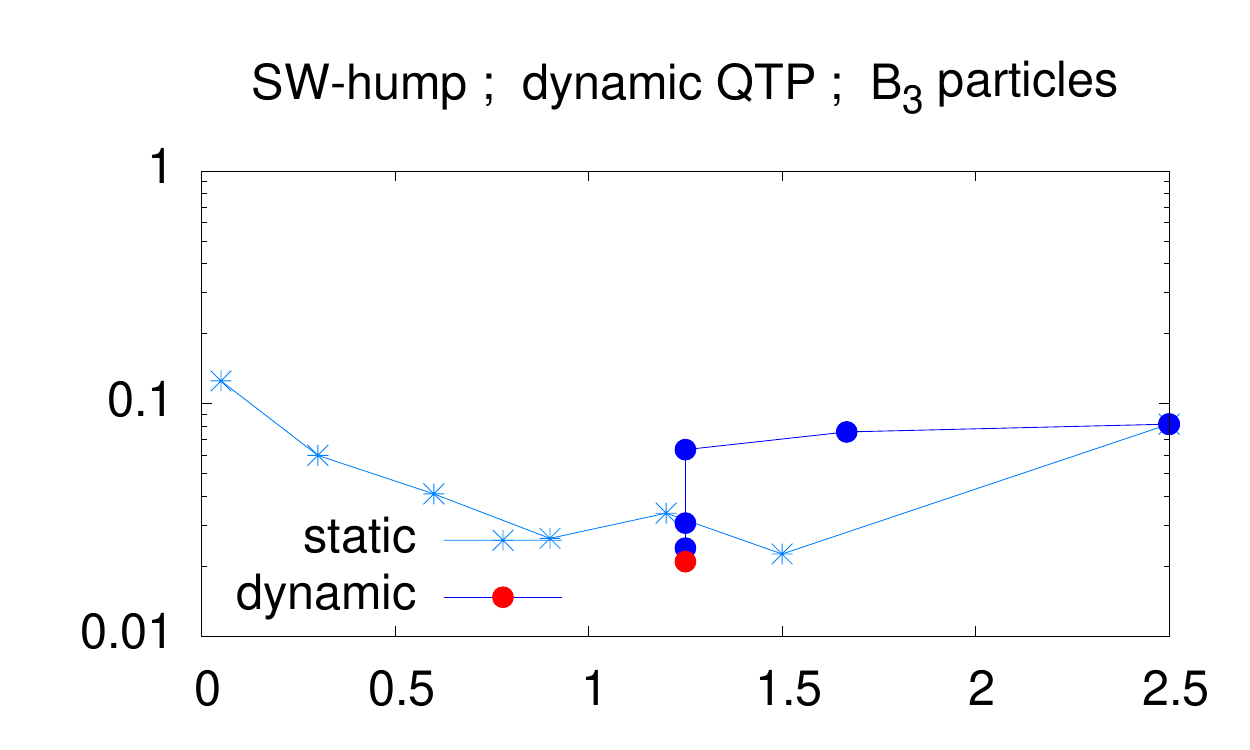}
\vspace{10pt}
\\
\includegraphics[width=0.45\textwidth]{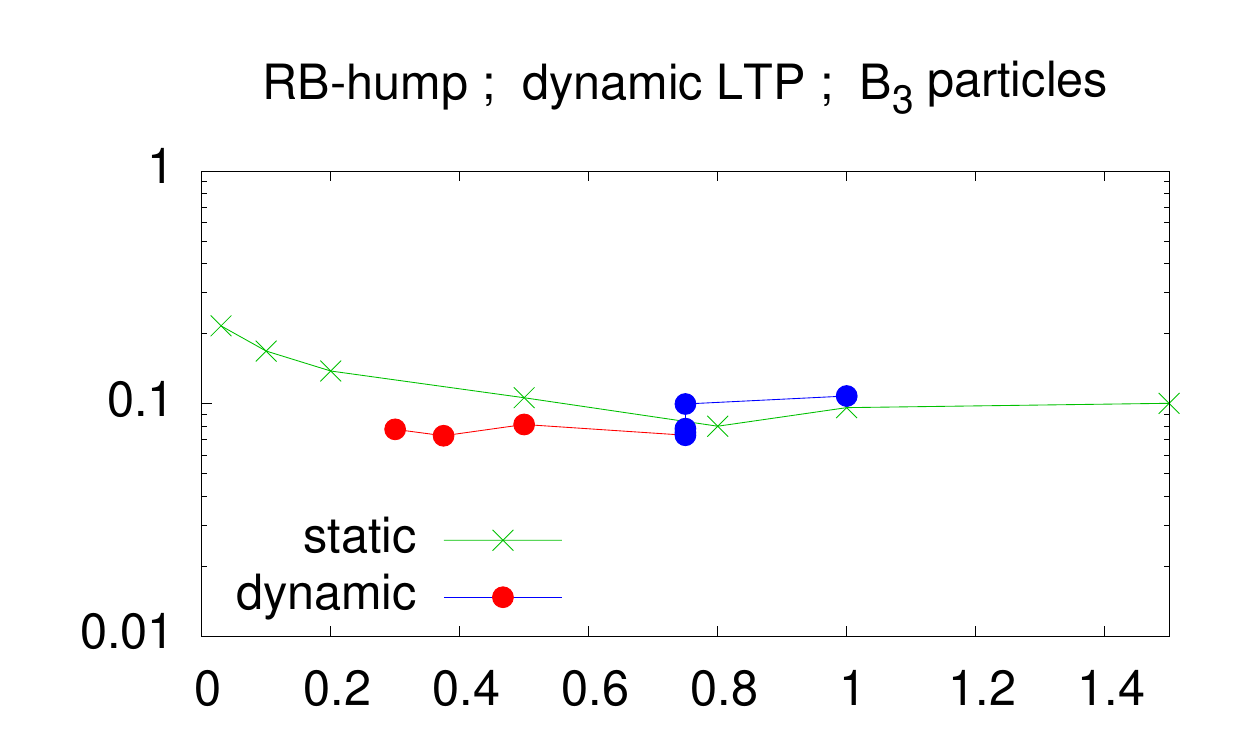}
& % \hspace{5pt} 
\includegraphics[width=0.45\textwidth]{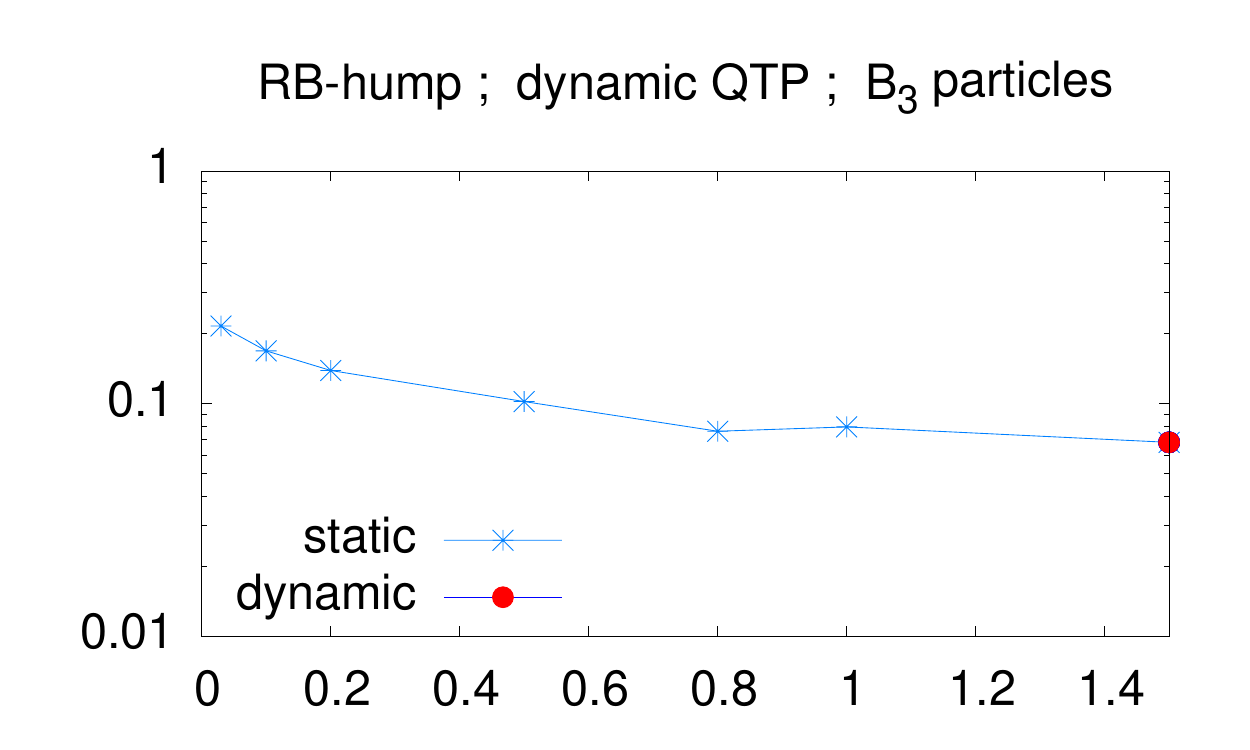}
\vspace{10pt}
\\
\includegraphics[width=0.45\textwidth]{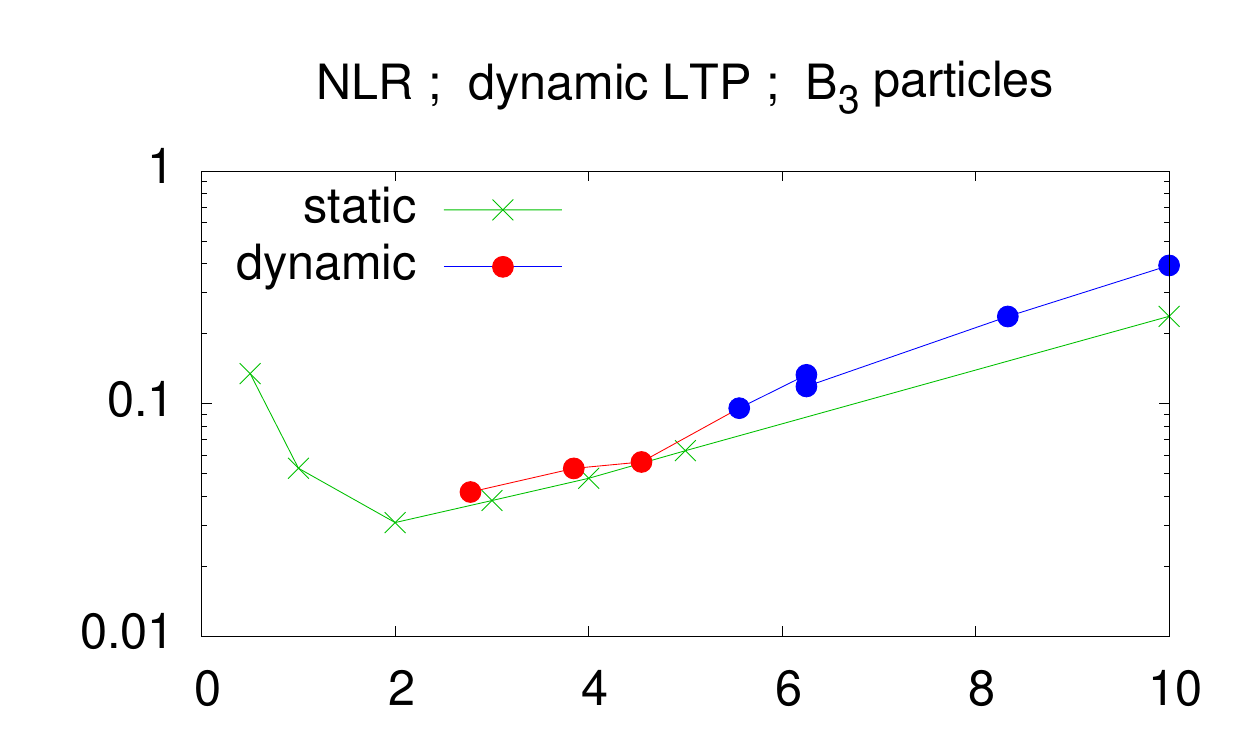}
& % \hspace{5pt} 
\includegraphics[width=0.45\textwidth]{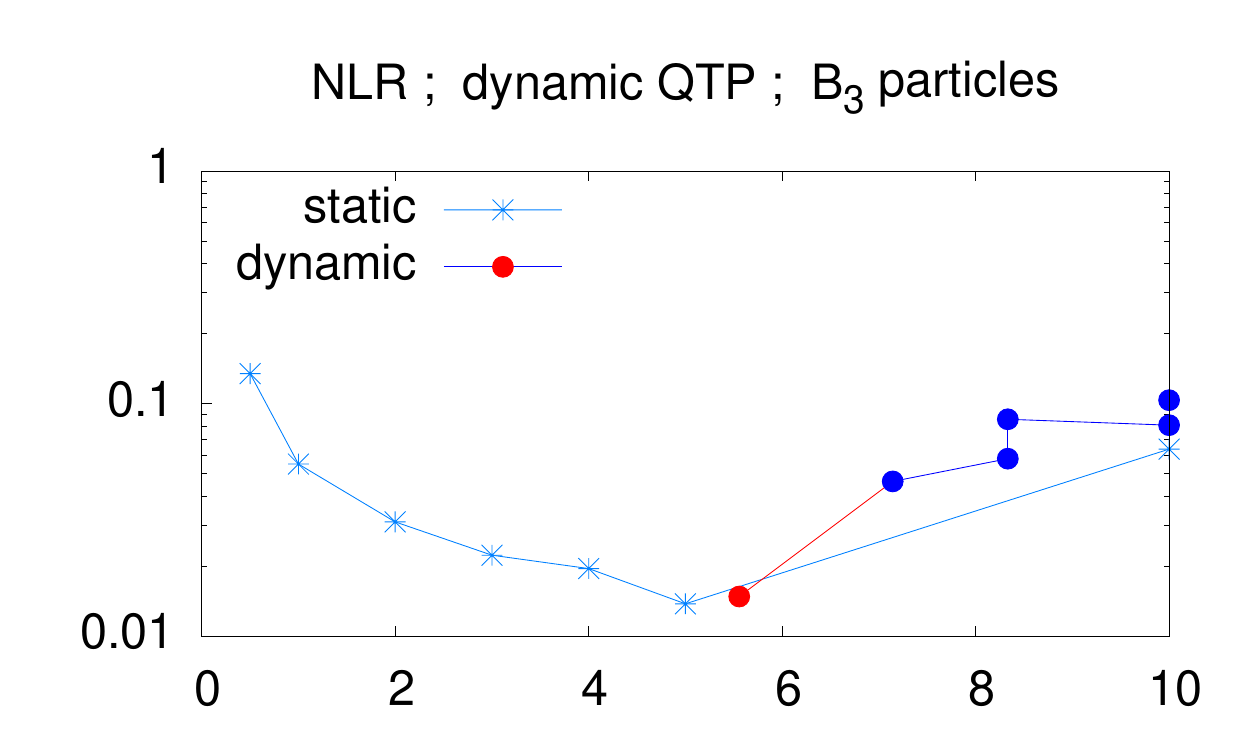}
\end{tabular}
  \caption{(Color) Performance of the dynamic remapping strategy applied to the LTP and QTP method.
    Here the relative $L^\infty$ errors at $t=T$ are plotted vs. the average remapping periods, 
    for the four test cases defined in Table~\ref{tab:test-cases} (every run is obtained with
    $256\times 256$ particles).
    In the static runs the abscissa represents the constant remapping period $\Dtr$ as in 
    Figure~\ref{fig:err-Dtr}.
    In the dynamic runs it corresponds to the ratio $T/R$ where $R$ is the 
    number of remappings (initialization included) selected by the criterion \eqref{rem-crit} with the 
    error indicators \eqref{E-transp-err}, \eqref{E-remap-err}.
    Finally the different points in the dynamic runs correspond to different values of the 
    constant in \eqref{rem-crit}. For the LTP runs the red (resp. blue) points correspond to 
    values larger (resp. smaller) than $C_{\rm remap} = 1$. 
    For the QTP runs the threshold value is $C_{\rm remap} = 5$.
    }
  \label{fig:dyn-rem}
 \end{center}
\end{figure}

In Figure~\ref{fig:dyn-rem} we plot the results obtained with the above dynamic criterion \eqref{rem-crit} and compare them with 
those obtained with a static strategy. To do so, for each case we plot both error curves using 
as $x$-axis an average remapping period : for the static runs it is the constant remapping period (hence the curves correspond
to those already shown in Figure~\ref{fig:err-Dtr}) and for the dynamic runs it is defined as $T/(1+R)$ where $R$ is the 
number of dynamic remappings (initialization included). 

Finally, the different points in the dynamic runs correspond to different values of the 
constant in \eqref{rem-crit}. For the LTP runs the red (resp. blue) points correspond to values larger (resp. smaller)
than $C_{\rm remap} = 1$. For the QTP runs the threshold value is $C_{\rm remap} = 5$.

From the results shown in Figure~\ref{fig:dyn-rem} we may draw a positive conclusion: indeed in every case but one, our dynamic
remapping strategy achieves the same level of accuracy as the best static run (in some cases it is even more accurate), while using
about the same number of remappings. And in the NLR case where the dynamic strategy fails to reach that optimal accuracy 
one could argue that the solutions are very smooth, hence by only considering the first order error estimate \eqref{first-order}
we obtain a remapping error indicator that is overly pessimistic. One would expect better results with a numerical indicator based 
on a second order estimate.

%%%%%%%%%
%%%%%%%%%

\section{Conclusion}

We have introduced a new class of particle methods with polynomial deformations
of arbitrary degree $r$, and for transport problems with underlying smooth characteristic flow
we have established their $L^\infty$ convergence with order $r$. These estimates come with uniform
bounds on the particle overlapping, and no particle remapping is required.

For practical applications we have described and tested fully discrete implementations 
of the first and second order cases. %, where particle shapes are deformed linearly and quadratically. 
An important feature of the resulting linearly-transformed particle (LTP) and quadratically-transformed particle 
(QTP) schemes is that they only involve pointwise evaluations 
of the forward characteristic flow. Since virtually every particle code contains routines that compute accurate
approximations of this flow to push the particles forward in time, this feature should facilitate their writing 
in existing codes, with small additional programming cost.

By testing the proposed LTP and QTP methods on a series of 2d test problems we have demonstrated both their 
enhanced accuracy and robustness with respect to the remapping period. In particular, we have observed that with the LTP 
scheme the optimal remapping periods could be as large as about 10 to 50 times what they are for standard (fixed-shape)
remapped particle methods, and with the QTP scheme this ratio would often reach values beyond 30.

We thus believe that this approach has a few potential advantages over standard smoothed ($\ve$-convoluted or remapped) 
particle methods: first, the scale $\ve$ of the particle shape functions can be taken equal to that of the initialization grid, 
which makes the overlapping of the particles uniformly bounded with respect to the resolution of the method. 
Second, the method is more accurate and much more robust with respect to particle remappings, leading to less numerical dissipation
-- two properties already observed in a recent work \cite{CP.S.F.G.L.2013.arxiv} where an LTPIC scheme has been used to solve Vlasov-Poisson plasmas. 
In parallel codes where inter-particle communications are expensive, fewer remappings also imply lower computational costs.

\begin{acknowledgements}
The author thanks Eric Sonnendr\"ucker, Albert Cohen, Jean-Marie Mirebeau and Jean Roux for valuable discussions during 
the early stages of this research.
\end{acknowledgements}

% BibTeX users please use one of
%\bibliographystyle{spbasic}      % basic style, author-year citations
\bibliographystyle{spmpsci}      % mathematics and physical sciences
\bibliography{campos-spmws.bib}   % name your BibTeX data base

% Non-BibTeX users please use
%\begin{thebibliography}{}
%%
%% and use \bibitem to create references. Consult the Instructions
%% for authors for reference list style.
%%
%\bibitem{RefJ}
%% Format for Journal Reference
%Author, Article title, Journal, Volume, page numbers (year)
%% Format for books
%\bibitem{RefB}
%Author, Book title, page numbers. Publisher, place (year)
%% etc
%\end{thebibliography}

\end{document}